\documentclass[11pt, a4paper]{article}
\usepackage{amsmath, amsthm, amssymb, eucal}
\usepackage{palatcm}
\newcommand{\partentry}[1]{\addtocontents{toc}
{\small\bfseries#1\hfill\thepage\par}}
\makeatletter
\def\@part[#1]#2{%
    \ifnum \c@secnumdepth >\m@ne
      \refstepcounter{part}
      \partentry{\protect\makebox[2em][l]{\thepart}#1}
\else
      \partentry{#1}
    \fi
    {\noindent\normalfont\Large\bfseries\thepart\hspace{1em}#2\par}
    \nobreak
    \vskip 3ex
    \@afterheading}
\def\@spart#1{%
    {\noindent\normalfont\Large\bfseries #1\par} 
     \nobreak
     \vskip 3ex
     \@afterheading}
\renewcommand\section{\@startsection{section}{1}{\z@}
{-3.5ex \@plus -1ex \@minus -.2ex}
{2ex \@plus.2ex}
{\large\bfseries}}
\renewcommand\subsection{
\@ifstar{\setcounter{subsection}{\value{equation}}
\@startsection{subsection}{2}{\z@}
                          {1.75ex \@plus.5ex \@minus.2ex}%
                           {-.4em} 
			{\itshape}*}
{\setcounter{subsection}{\value{equation}}
\stepcounter{equation}
\@startsection{subsection}{2}{\z@}
                          {1.75ex \@plus.5ex \@minus.2ex}%
                           {-.4em} 
{\itshape}}}
\def\@seccntformat#1{\@ifundefined{#1@cntformat}%
{\csname the#1\endcsname\quad} 
{\csname #1@cntformat\endcsname}}
\def\section@cntformat{\thesection.~}
\def\subsection@cntformat{(\thesubsection)\ }
\renewcommand*\l@section{\mdseries\small\@dottedtocline{1}{1.5em}{2em}}
\makeatother

\numberwithin{equation}{section}
\theoremstyle{plain}
\newtheorem{maintheorem}{Theorem}
\swapnumbers
\newtheorem{theorem}[equation]{Theorem}
\newtheorem{corollary}[equation]{Corollary}
\newtheorem{lemma}[equation]{Lemma}
\newtheorem{proposition}[equation]{Proposition}
\theoremstyle{definition}
\newtheorem{definition}[equation]{Definition}
\theoremstyle{remark}
\newtheorem{remark}[equation]{Remark}

\newcommand{\cliff}{\mathrm{Cliff}}
\newcommand{\fock}{\mathbf{F}}
\newcommand{\hilb}{\mathbf{H}}
\newcommand{\lilb}{\mathbf{L}}
\newcommand{\philb}{\mathbf{H}'}

\newcommand{\spin}{\mathbf{S}}
\newcommand{\pspin}{\mathbf{S}'}

\newcommand{\naff}{N_\mathrm{aff}}
\newcommand{\nafff}{N_\mathrm{aff}}
\newcommand{\waff}{{W_\mathrm{aff}}}
\newcommand{\QB}{Q_I}
\newcommand{\rS}{\sqrt[r]{S^1}}
\newcommand{\vep}{\varepsilon}
\newcommand{\mi}{\mathrm{i}}
\newcommand{\ul}{\underline}
\newcommand{\ol}{\overline}
\newcommand{\plg}{L'\frg}
\newcommand{\hatplg}{\widehat{L}'\frg}
\newcommand{\tildeplg}{\widetilde{L}'\frg}
\newcommand{\plge}{L'_\vep\frg}
\newcommand{\hatplge}{\widehat{L}'_\vep \frg}
\newcommand{\wtlat}{\Lambda}
\newcommand{\kstar}{K^*}
\newcommand{\dirac}{\textup{\mbox{D\hspace{-0.6em}\raisebox{.1ex}%
{/}\hspace{.07em}}}}

\newcommand{\ofrac}[2]{\genfrac{}{}{0pt}{}{#1}{#2}}
\newcommand{\ad}{\mathrm{ad}}
\newcommand{\Ad}{\mathrm{Ad}}
\newcommand{\cA}{\mathcal{A}}

\newcommand{\cN}{\mathcal{N}}
\newcommand{\cO}{\mathcal{O}}

\newcommand{\cV}{\mathcal{V}}
\newcommand{\cU}{\mathcal{U}}
\newcommand{\fra}{\mathfrak{a}}
\newcommand{\frg}{\mathfrak{g}}
\newcommand{\frh}{\mathfrak{h}}
\newcommand{\frn}{\mathfrak{n}}
\newcommand{\frt}{\mathfrak{t}}
\newcommand{\frz}{\mathfrak{z}}

\newcommand{\frN}{\mathfrak{N}}
\newcommand{\frO}{\mathfrak{O}}
\newcommand{\frP}{\mathfrak{P}}
\newcommand{\frU}{\mathfrak{U}}

\newcommand{\bC}{\mathbb{C}}
\newcommand{\bP}{\mathbb{P}}
\newcommand{\bR}{\mathbb{R}}
\newcommand{\bT}{\mathbb{T}}
\newcommand{\bZ}{\mathbb{Z}}

\newcommand{\bfv}{\mathbf{v}}

\begin{document}                                                                                                   
\title{\textbf{Loop Groups and twisted $K$-theory III}}
\author{Daniel S.~Freed \and Michael J.~Hopkins \and Constantin Teleman} 
\date{10 November 2005} 
\maketitle

\begin{quote}
\abstract{\noindent This is the third paper of a series relating the 
equivariant twisted $K$-theory of a compact Lie group $G$ to the 
``Verlinde space'' of isomorphism classes of projective lowest-weight representations of the loop groups. Here, we treat arbitrary compact 
Lie groups. In addition, we discuss the relation to semi-infinite 
cohomology, the fusion product of Conformal Field theory, the r\^ole 
of energy and the topological Peter-Weyl theorem.}
\end{quote}
\vskip 1cm

\part*{Introduction}
In~\cite{fht1, fht2} the twisted equivariant $K$-theory of a compact 
Lie group was described in terms of positive energy representations of 
its loop group. There, we assumed that the group was connected, with 
torsion-free fundamental group. Here, we remove those restrictions; 
we also relax the constraints on the twisting, assuming only its
\textit {regularity}. Additional constraints allow the introduction 
of an \textit{energy} operator, matching the rotation of loops, and 
lead to the \textit{positive energy representations} relevant to 
conformal field theory. Finer restrictions on the twisting lead to 
a structure of 2-dimensional topological field theory, the ``Verlinde 
TFT" \cite{fht1}. This is not discussed here, but we do prove two of 
the key underlying results: we identify the fusion product with the 
topological cup-product, and equate the topologically constructed TFT 
bilinear form with the duality pairing between irreducible representations 
at opposite levels.

Capturing the Verlinde ring topologically lets us revisit, via twisted 
$K$-theory, some constructions on representations that were hitherto 
assumed to rely on the algebraic geometry of loop groups. Thus, 
restriction to and induction from the maximal torus in twisted 
$K$-theory recover of \textit{semi-infinite restriction} and \textit
{induction} of Feigin and Frenkel \cite{ff} on representations.
The \textit{energy} operator comes from the natural circle action on 
the quotient stack of $G$, under its own conjugation action. The numerator 
in the character formula can be obtained by dualising the Gysin inclusion 
of the identity in $G$. Next, the cup-product action of $R(G)$ on 
$K^\tau_G(G)$ corresponds to the \textit{fusion} of Conformal Field 
Theory, defined via holomorphic induction. Finally, we discuss the 
Borel-Weil theorem for the ``annular" flag variety of a product of 
two copies of the loop group, interpreted now as a topological Peter-Weyl 
theorem. This last result can be interpreted as a computation of the 
TFT bilinear form mentioned earlier, but in addition, it can be 
further extended to an index theorem for generalised flag varieties of 
loop groups, in which twisted $K$-theory provides the topological side. 
We refer to \cite[\S8]{fht3} for a verification of this result in the 
special case of connected groups with free $\pi_1$, and to \cite{tel2} 
for further developments concerning higher twistings of $K$-theory.

The paper is organised as follows. Chapter I states the main theorems 
and describes the requisite technical specifications. Two examples are 
discussed in Chapter II: the first relates our theorem in the case of 
a torus to the classical spectral flow of a family of Dirac operators, 
while the second recalls the Dirac family associated to a compact group
\cite{fht2}, whose loop group analogue is the ``non-abelian spectral 
flow" implementing our isomorphism. Chapter III computes the twisted $K$-theory 
$K^\tau_G(G)$ topologically, by reduction to the maximal torus and its 
normaliser in $G$. Chapter IV reviews the theory of loop groups and 
their lowest-weight representations; the classification of irreducibles 
in \S\ref{repclassifsect} reproduces the basis for $K^\tau_G(G)$ constructed 
in Chapter III. The Dirac family in Chapter V assigns a twisted $K$-class 
to any (admissible) representation of the loop group, and this is shown 
to recover the isomorphism already established by our classification. 
Chapter VI gives the topological interpretation of some known constructions 
on loop group representations as discussed above. Appendix A reviews 
the diagram automorphisms of simple Lie algebras and  relates our 
definitions and notation with those in Kac \cite{kac}. 

\section*{Acknowledgements.} The authors would like to thank Graeme Segal 
for many useful conversations. We also thank G.~Landweber for detailed 
comments on the early version of the manuscript and U.~Bunke for discussions 
on the G\"ottingen seminar notes \cite{bunk}. During the course of this 
work, the first author was partially supported by NSF grants DMS-0072675 
and DMS-0305505, the second partially by NSF grants DMS-9803428 and 
DMS-0306519, and the third partially by NSF grant DMS-0072675.  We also 
thank the KITP of Santa Barbara (NSF Grant PHY99-07949) and the Aspen Centre 
for Physics for hosting their "topology and physics" programs, where 
various sections of this paper were revised and completed.

\section*{Index of Notation} 

\begin{tabular}{ll}
\textit{Groups} & \\
$G$, $G_1$ & Compact Lie group and its identity component\\
$T$, $N$ & Maximal torus and its normaliser in $G$ \\
$W$, $W_1$ & Weyl groups $N/T$, $N\cap G_1/T$ of $G$ and $\frg$\\
$G(f)$, $N(f)$ & Centralisers in $G$ and $N$ of the connected 
	component of $f$ (\S\ref{topktorus})\\
$\frg$,  $\frt$, $\frg_{\bC}$, $\frt_{\bC}$ & 
	Lie algebras and their complexifications\\
$\frn\subset\frg_{\bC}$ & Sub-algebra spanned by
	the positive root vectors\\ 
$\langle\ |\ \rangle$, $\{\xi_a\}$ & Basic inner product on $\frg$ 
	(when semi-simple); orthonormal basis (\S\ref{finitedirac})\\ 
$\rho, \theta \in \frt^*$ & Half-sum of positive roots,
	highest root (\S\ref{finitedirac})\\
$h^\vee $ & (for simple $\frg$) Dual Coxeter number $\rho\theta + 1$ 
	(\S\ref{finitedirac})\\ 
$\vep$; $\ul\frg$, $\ul\frt$ & Diagram automorphism of $\frg$; 
	$\vep$-invariant sub-algebras (\S\ref{topkcompact})\\ 
$\ul {W}$; $\ul T$ &	Weyl group of $\ul {\frg}$; torus $\exp(\ul\frt)$
	(\S\ref{topkcompact})\\
$\ul {\rho}, \ul{\theta} \in \ul{\frt}^*$ & Half-sum of positive roots 
	in $\ul \frg$, highest $\ul\frg$-weight of ${\frg}/{\ul \frg}$ 
	(\S\ref{twistedaff})\\
$R$, $R^\vee $ & Root and co-root lattices\\ 
$\wtlat$, $\ul\wtlat$; $\ul{\wtlat}^\tau$ & Weight lattices of $T$ and 
	$\ul T$; lattice of $\tau$-affine weights \\
\end{tabular}

\noindent
\begin{tabular}{ll}
\textit{Loop Groups} & \\
$LG$, $L_fG$ & Smooth loop group and twisted loop group 
	(\S\ref{stat})\\
$LG^\tau$ &  Central extension by $\mathbb{T}$ with cocycle $\tau$\\
$L\frg$, $L_f \frg$  & Smooth loop Lie algebras \\
$\plg$, $L'G$ & Laurent polynomial Lie algebra (\S\ref{affnot}), 
	loop group (\S\ref{fusion})\\
$\naff = \Gamma_fN$ & Group of (possibly $f$-twisted) geodesic loops in $N$ 
	(\S\ref{topktorus}) \\
$\waff(\frg,f)$, $\waff$ & $f$-twisted affine Weyl group of $\frg$, 
	extended affine Weyl group $\nafff/\ul T$ \\
$\fra,\ul\fra$ & (simple $\frg$) Alcove of dominant $\xi\in\frt,\ul\frt$ 
	with $\theta(\xi) \le 1$, resp. $\ul{\theta}(\xi)\le1/r$ \\
$\tau\cdot\ul\fra^*\subset\ul\frt$ & Product of the centre of $\ul\frg$ and 
	the $[\tau]$-scaled alcoves on simple factors (\S\ref{repclassifsect})\\
\end{tabular}

\vspace{0.5cm}
\noindent
\begin{tabular}{ll}
\textit{Twistings} & \\
$\tau$; $[\tau]$ & 2-cocycle on $LG$; level in $H_G^3(G;\bZ)$ 
	(\S\ref{topregsect})\\
$\kappa^\tau$ & Linear map $H_1(\ul T)\to H^2(B\ul T)$ defined from 
	$[\tau]$ (\S\ref{topktorus})\\
$\sigma$, $\ul\sigma$ & $LG$-cocycle of the Spin modules for $L\frg$, 
	$L_f\frg$ (\S\ref{morita})\\
$\sigma(\frt)$, $\sigma(\ul\frt)$ & $W$-cocycle for the spinors on $\frt$ 
	and $\ul\frt$ (\S\ref{topktorus})\\ 
$\tau'$, $\tau''$ & Twisting for the $\waff$-action on $\ul{\wtlat}^\tau$; 
	shifted twisting $\tau'-\sigma(\ul\frt)$ (\S\ref{topktorus})\\
\end{tabular}

\section*{}
\begin{minipage}[t]{12cm}
\tableofcontents
\end{minipage}

\part{Statements}
\label{stat}

Throughout the paper, cohomology and $K$-theory have integer coefficients, 
if no others are specified. $K$-theory has \textit{compact supports}; however, 
for proper actions of non-compact groups, or for stacks in general, this 
refers to the quotient space. For a twisting $\tau$ on $X$, the twisted 
$K$-theory will be denoted $K^\tau(X)$. This is a $\bZ/2$-graded group, 
whose two components are denoted $K^{\tau+0}(X), K^{\tau+1}(X)$. For a
central extension $G^\tau$ of $G$, the Grothendieck group of $\tau$-projective
representations is denoted by $R^\tau(G)$; it is a module over the 
representation ring $R(G)$.\footnote{Note that, when $\tau$ is graded, 
this module can have an \textit{odd} component, cf.~\S\ref{untwisted}.}

\section{Main theorems}
\subsection{Simply connected case.} \label{singlespec} 
The single most important special case of our result concerns a simple, 
simply connected compact Lie group $G$. Central extensions of its smooth 
loop group $LG$ by the circle group $\bT$ are classified by their \textit
{level}, the Chern class $c_1\in H^2(LG) = \bZ$ of the underlying circle 
bundle. These extensions are equivariant under loop rotation. Among the 
projective representations of $LG$ with fixed level $k$ are the \textit
{positive energy} ones: they are those which admit an intertwining action 
of the group of loop rotations, with spectrum bounded below. Working up 
to infinitesimal equivalence, as is customary with non-compact groups, 
these representations are semi-simple, with finitely many irreducibles, 
all of them unitarisable \cite{psloop}. The free abelian group on irreducible 
isomorphism classes is called the \textit{Verlinde ring of $G$ at level 
$k$}; the multiplication is the \textit{fusion} of Conformal Field Theory. 
We denote this ring by $R^k(LG)$, by analogy with the representation ring 
$R(G)$ of $G$. 

\begin{maintheorem} 
If $k+h^\vee >0$, $R^k(LG)$ is isomorphic to the twisted K--theory 
$K_G^{[k+h^\vee]+d}(G)$. 
\end{maintheorem}
\noindent Here, $G$ acts on itself by conjugation, $h^\vee$ is the \textit
{dual Coxeter number} of $G$, $[k+h^\vee]$ is interpreted as a \textit
{twisting class} in $H_G^3(G)\cong \bZ$ for equivariant $K$-theory, while 
$d:=\dim G$ is a degree-shift: the two sides are supported in degrees $0$, 
resp.\ $d \mod 2$. The ring structure on $K$-theory is the convolution 
(Pontryagin) product. The isomorphism is established by realising both 
sides as quotient rings of $R(G)$, via \textit{holomorphic induction} 
on the loop group side, and via the Thom push-forward from the identity 
in $G$, on the $K$ side.

\subsection{General groups.}\label{conditions}
The isomorphism between the two sides and the relation between level 
and twisting cannot be described so concisely for general compact Lie 
groups. This is due to the presence of torsion in the group $H^3$ of
twistings, to an additional type of twistings classified by $H^1_G(G;
\bZ/2)$, related to gradings of the loop group, and to the fact that the 
two sides need not\footnote{A simple statement can be given when $G$ is 
connected and $\pi_1(G)$ is free \cite[\S6]{fht3}, precisely because  
both sides are quotients of $R(G)$.} be quotients of $R(G)$. For a 
construction of the map via a correspondence induced by conjugacy classes, 
we refer to \cite{freedicm} (see also \S\ref{inductorb} here). Ignoring 
the difficulties for a moment, there still arises a natural isomorphism 
between the twisted equivariant $K$-groups of $G$ and those of the 
category of positive energy representations at a \textit{shifted} level, 
provided that:
\begin{enumerate} 
\itemsep0ex 
\item we use $\bZ/2$-\textit{graded} representations; 
\item we choose a central extension of $LG$ which is equivariant
under loop rotation; 
\item the cocycle of the extension satisfies a positivity 
condition. 
\end{enumerate}
The energy operator cannot be defined without (ii), and without (iii), 
representations of positive energy do not exist. While (iii) is 
merely a question of choosing the correct sign, an obstruction to 
equivariance in (ii) is the absence of symmetry in the level (see 
\S\ref{energysect}). This can only happen for tori --- whose loop groups, 
ironically, have a simple representation theory. 

This formulation is unsatisfactory in several respects. The loop group 
side involves the \textit{energy}, with no counterpart in $K_G(G)$: 
instead, a rotation-equivariant version of the latter will be more 
relevant. There is also the positivity restriction, whereas the 
topological side is well-behaved for \textit{regular} twistings (\S\ref
{techdef}). There is, finally, the unexplained ``dual Coxeter" shift. 

We now formulate the most canonical statement. This need not be the most 
useful one (see Theorems 3 and 5 instead). However, it has the virtue of 
explaining the shift between level and twisting, as the projective cocycle 
of the positive energy spinors on $LG$. Gradings in (i), if not originally 
present in the twisting $\tau$, are also imposed upon us by the spinors 
whenever the Ad-representation of $G$ does not spin. 

\subsection{Untwisted loop groups.}\hspace{-.4em}\footnote
{\textit{Twisting} for loop groups (\S\ref{twistedgroups}) and for 
$K$-theory mean different things, but both uses are well-entrenched.}
\label{untwisted} 
Let $G$ be any compact Lie group and $LG^\tau$ a smooth $\bT$-central 
extension of its loop group. We allow $LG$ to carry a \textit{grading}, 
or homomorphism to $\bZ/2$; this is classified by an element of 
$H^1_G(G_1;\bZ/2)$ and is notationally incorporated into $\tau$. An 
Ad-invariant $L^2$ norm on $L\frg$ defines the graded Clifford 
algebra\footnote{This algebra should really be based on \textit{half-forms} 
on the circle, which carry a natural bilinear form.} $\cliff (L\frg^*)$, 
generated by odd elements $\psi(\mu)$, $\mu \in L\frg^*$, with 
relations $\psi(\mu)^2= \|\mu\|^2$.

A \textit{$\tau$-representation} of $LG$ is a graded representation of 
$LG^\tau$ on which the central circle acts by the natural character. We 
are interested in complex, graded $\tau$-representations of the crossed 
product $LG\ltimes \cliff (L\frg^*)$, with respect to the co-adjoint action. 
Graded modules for $\cliff(L\frg^*)$ can be viewed as \textit{$b$-projective representations} of the odd vector space $\psi(L\frg^*)$, where $b$ is 
the $L^2$ inner product, so we are considering $(\tau,b)$-representations 
of the \textit{graded super-group} $LG_s := LG\ltimes\psi(L\frg^*)$. 
Subject to a \textit{regularity restriction} on $\tau$, an \textit
{admissibility condition} on representations will ensure their complete 
reducibility (\S\ref{techdef}).

A \textit{super-symmetry} of a graded representation is an odd automorphism
squaring to $1$. Let $R^{\tau+0}(LG_s)$ be the free $\bZ$-module of graded 
admissible representations, modulo super-symmetric ones, and $R^{\tau+1}
(LG_s)$ that of representations with a super-symmetry, modulo those carrying 
a second super-symmetry anti-commuting with the first. These should be 
regarded as the $LG_s$-equivariant $K^\tau$-groups of a point. The reader 
should note that defining $K$-theory for graded algebras a delicate matter 
in general \cite{black}; the shortcut above, also used in \cite[\S4]
{fht3}, relies on the semi-simplicity of the relevant categories of modules. 

Since $K^\tau_G(G)$ is a $K_G(G)$-module, it carries in particular an
action of the representation ring $R(G)$. \textit{Fusion} with 
$G$-representations defines an $R(G)$-module structure on $R^\tau(LG_s)$; 
the definition is somewhat involved, and we must postpone it until \S\ref
{fusion}. Here is our main result.

\begin{maintheorem} 
For regular $\tau$, there is a natural isomorphism of (graded) $R(G)
$-modules $R^\tau(LG_s) \cong K_G^\tau(G)$, wherein $K$-classes arise 
by coupling the Dirac operator family of Chapter \ref{diracfam} to 
admissible $LG_s$-modules.
\end{maintheorem}

\begin{remark}
For twistings that are suitably \textit{transgressed from $BG$}, both 
sides carry isomorphic Frobenius ring structures. The portion of the 
product structure that exists for any regular twisting is discussed 
in \S\ref{fusion}. A geometric construction of the duality pairing 
is described in \S\ref{peterweylsection}.
\end{remark}

\subsection{Twisted loop groups.}\label{twistedgroups} 
When $G$ is disconnected, there are twisted counterparts of these notions. 
For any $f\in G$, the \textit{twisted loop group} $L_fG$ of smooth maps 
$\gamma: \bR \to G$ satisfying $\gamma(t+2\pi) = f\gamma(t)f^{-1}$ 
depends, up to isomorphism, only on the conjugacy class in $\pi_0G$ of 
the component $fG_1$ of $f$. Let $[fG_1]\subset G$ denote the union of 
conjugates of $fG_1$; the topological side of the theorem is $K_G^\tau
\left([fG_1]\right)$, while the representation side involves the 
admissible representations of $L_fG\ltimes \psi(L_f\frg^*)$. 
 
\subsection{Removing the spinors.}\label{morita}
A lowest-weight spin module $\spin$ for $\cliff(L\frg^*)$ (see \S\ref 
{admissibility}) carries an intertwining projective action of the loop 
group $LG$. Denoting by $\sigma$ ($\ul\sigma$, in the $f$-twisted case) 
the projective cocycle of this action and by $\ul d$ the dimension of 
the centraliser $G^f$, a \textit{Morita isomorphism}
\begin{equation}
R^\tau(L_fG_s)\cong R^{\tau-\ul\sigma-\ul d}(L_fG)
\end{equation}
results from the fact that an admissible, graded $\tau$-module of 
$L_fG\ltimes\cliff(L_f\frg^*)$ has the form $\hilb\otimes\spin$, for a 
suitable $(\tau-\ul\sigma)$-representation $\hilb$ of $L_fG$, unique up to 
canonical isomorphism. Note in particular the dimension shift by $\ul d$, 
from the parity of the Clifford algebra. We obtain the following 
reformulation of Theorem~2. 
\begin{maintheorem} 
For regular $\tau$, there is a natural isomorphism 
$K_G^\tau\left([fG_1]\right)\cong R^{\tau-\ul\sigma-\ul d}(L_fG)$. 
\end{maintheorem}
\noindent The loop group may well acquire a grading from the spinor twist 
$\ul\sigma$, even if none was present in $\tau$; if so, $R^{\tau-\ul\sigma}
(L_fG)$ is built from graded representations, as in \S\ref{untwisted}. 

\subsection{Classifying representations.} 
\label{repclasstatement} In proving the theorems, we compute both sides 
of the isomorphism in Theorem~3. More precisely, we compute the twisted 
$K$-theory by reduction to the torus and the Weyl group, and produce an answer 
which agrees with the classification of admissible representations in terms 
of their lowest weights. In fact, twisted $K$-theory allows for an attractive formulation of the lowest-weight classification for disconnected (loop) 
groups, as follows. 

Choose a maximal torus $T\subset G$ which, along with a dominant chamber, 
is stable under $f$-conjugation. (Such tori always exists, see Proposition 
\ref{diagaut}.) Recall that the \textit{extended affine Weyl group} $\waff$ 
for $L_fG$ is $\pi_0$ of $L_fN$, the group of $f$-twisted loops in the 
normaliser $N$ of $T$. Let $\ul T\subset T$ denote the subtorus centralised 
by $f$, and $\ul{\wtlat}^\tau$ the set of its $\tau$-affine weights. The 
conjugation action of $L_fN$ on $\ul T$ descends to an action of $\waff$ on 
$\ul{\wtlat}^\tau$, which preserves the subset $\ul{\wtlat}^\tau_{reg}$ of 
\emph{regular} weights. A tautological twisting $\tau'$ is defined for this 
action, because every weight defines a $\bT$-central extension of its 
centraliser in $\waff$ (see \S\ref{morewaff} for details). Finally, after 
projection to the finite Weyl group $W=N/T$, $\waff$ also acts on the Lie 
algebra $\ul\frt$ of $\ul T$.
\begin{maintheorem}
The category of graded, admissible $\tau$-representations of $L_fG\ltimes 
\cliff(L_f\frg^*)$ is equivalent to that of $\tau'$-twisted $\waff\ltimes
\cliff(\ul{\frt})$-modules on $\ul{\wtlat}^\tau_{reg}$. 
\end{maintheorem}  

\noindent It follows that the corresponding $K$-groups agree. 
We reduce Theorem 4 in \S\ref{repclassifsect} to the well-known cases 
of simply connected compact groups and tori. 

Computing both sides is a poor explanation for a natural isomorphism, 
and indeed we improve upon this in Chapter V by producing a map from 
representations to $K$-classes using families of Dirac operators. The 
construction bypasses Theorem~4 and ties in beautifully with Kirillov's 
\textit{orbit method}, recovering the co-adjoint orbit and line bundle
that correspond to an irreducible representation. Another offshoot of 
this construction emerges in relation with the \textit{semi-infinite 
cohomology} of Feigin and Frenkel \cite{ff}, for which we give a 
topological model (Theorem \ref{seminfres}): for integrable representations, 
the Euler characteristic of semi-infinite $L\frn$-cohomology becomes the 
restriction from $G$ to $T$, on the $K$-theory side. While this can also 
be checked by computing both sides, our Dirac family gives a more natural 
proof, providing the same rigid model for both. 

\subsection{Loop rotation.}\label{looprot}
Assume now that the extension $LG^\tau$ carries a lifting of the loop 
rotation action on $LG$. It is useful to allow \textit{fractional} 
lifts, that is, actions on $LG^\tau$ of a \textit{finite cover} $\bT$ 
of the loop rotation circle; such a lift always exists when $G$ is 
semi-simple (Remark \ref{ssfine}). If so, admissible $\tau$-representations 
carry an intertwining, semi-simple action of this new $\bT$.\footnote{A 
further positivity condition (\S\ref{poserg}) on $\tau$ ensures that the 
spectrum of this action is bounded below, and the real infinitesimal 
generator of the intertwining action is then called the \textit{energy}.} 
Schur's lemma implies that the action is unique up to an overall shift 
on any irreducible representation.

In this favourable situation, we can incorporate the loop rotation into 
our results. The requisite object on the topological side is the \textit
{quotient stack} of the space $\cA$ of $\frg$-valued smooth connections 
on the circle by the semi-direct product $\bT\ltimes LG$, the loop group 
acting by gauge transformations and $\bT$ by loop rotation. We denote 
the twisted $K$-theory of this stack by $K^\tau_\bT(G_G)$. This notation, 
while abusive, emphasises that the $\bT$-action makes it into an 
$R(\bT)$-module; its \textit{fibre over $1$} is the quotient by the 
augmentation ideal of $R(\bT)$. The following formulation, while awkward, 
has the virtue of being canonical; there is no natural isomorphism of 
$K_\bT^\tau(G_G)$ with $K_G^\tau(G)\otimes R(\bT)$. 
\begin{maintheorem}
If the regular twisting $\tau$ is rotation-equivariant, $K_\bT^\tau
(G_G)$ is isomorphic to the $R^\tau$-group of graded, admissible, 
representations of $\bT \ltimes LG_s$ (cf.~\S\ref{untwisted}). It is a 
free module over $R(\bT)$, and its fibre over $1$ is isomorphic to $
K_G^\tau(G)$.
\end{maintheorem}
\noindent A noteworthy complement to Theorem 5 is that $K_\bT^\tau
(G_G)$ contains the Kac numerator formula for $LG^\tau$-representations: 
see \S\ref{kachar}. It would be helpful to understand this as a twisted 
Chern character, just as the the Kac numerator at $q=1$ is the Chern 
character for $K_G^\tau(G)$ \cite{fht3}.


\section{Technical definitions}\label{techdef}
In this section, we describe our regularity conditions on the central 
extension $LG^\tau$ and define the class of admissible representations. 
There is a \textit{topological} and an \textit{analytical} component 
to regularity.  

\subsection{Topological regularity.}\label{topregsect} 
The central extension $LG^\tau$ has a characteristic class $[\tau]\in 
H_G^3(G_1)$, the \textit{level}. It is an equivariant version of the 
\textit{Dixmier-Douady invariant} of a gerbe, and arises from the 
connecting arrow in the exponential sequence for group cohomology 
with smooth circle coefficients, $H_{LG}^2(\bT) \to H_{LG}^3(\bZ)$: 
the last group is purely topological, and equals $H^3(BLG)\cong H_G^3(G_1)$. 
When $\frg$ is semi-simple, the smooth-cochain group cohomology $H_{LG}^2
(\bR)$ vanishes \cite[XIV]{psloop}, and $[\tau]$ then determines the 
central extension $LG^\tau$, up to isomorphism. In any case, restricting 
to a maximal torus $T\subset G$ and writing $H_T^2$ for $H^2(BT)$, we 
obtain a class in 
\[
H_T^3(T)= H^1(T)\otimes H_T^2 \oplus H^3(T). 
\]
For classes arising from central extensions, it turns out that the $H^3(T)$ 
component vanishes. In view of the isomorphism $H^1(T)\cong H^2_T$, we make 
the following

\begin{definition}\label{topreg}
We call $\tau$ \emph{topologically regular} iff $[\tau]$ defines 
a non-singular bilinear form on $H_1(T)$. 
\end{definition}

\noindent For a twisted loop group $L_fG$, topological regularity is 
detected instead by the $f$-invariant sub-torus $\ul T \subset T$ in an 
$f$-stable maximal torus $T$ as in \S\ref{repclasstatement}. Restricting 
$[\tau]$ there leads to a bilinear form on $H_1(\ul T)$, and regularity 
refers to the latter. In the next section, we will see how the bilinear 
form captures the commutation in $LT^\tau$ of the constant loops $T$ with 
the group of components $\pi_1 T$.

\subsection{Analytic regularity.} \label{anregsect}
This condition, which holds in the standard examples, concerns the 
centre $\frz\subset\frg$, and ensures that the topologically invisible 
summand $L\frz/\frz$ does not affect the classification of representations 
of $LG^\tau$. Split $L\frz$ into the constants $\frz$ and their 
$L^2$-complement $V$, and observe that $LG$ is the semi-direct product of 
the normal subgroup $\exp(V)$ by the subgroup $\Gamma$ of loops $\gamma$ 
whose velocity $d\gamma\cdot\gamma^{-1}$ has constant $\frz$-projection. 
Because the action of $\Gamma$ on $V$ factors through the \textit{finite} 
group $\pi_0G$, invariant central extensions of $\exp(V)$ have a preferred 
continuation to $LG$. 

\begin{definition}\label{anreg}
$\tau$ is \emph{analytically regular} iff it is the sum of an extension 
of $\Gamma$ and a Heisenberg extension of $\exp(V)$, and, moreover,  
the Heisenberg cocycle $\omega:\Lambda^2V \to \mi\bR$ has the form 
$\omega(\xi, \eta) = b(S\xi,\eta)$, for some skew-adjoint Fredholm 
operator $\mi S$ on $V$.
\end{definition} 

\noindent The standard example\footnote{This is the only possibility 
for $\mathrm{Diff}(S^1)$-equivariant extensions \cite[VIII]{psloop}.} 
has $S= \mi d/dt$, an unbounded operator, so we really ask that $S/
(1+\sqrt{S^*S})$ should be Fredholm. We need to tame $\omega$ for the 
Dirac constructions in Chapter \ref{diracfam}. For twisted loop groups, 
the analytic constraints refer to ${L_f\frz}\left/{\frz^f}\right.$.
  
\subsection{Linear splittings.} \label{linsplit}
Restricted to any simple summand in $\frg$, every extension class is a 
multiple of the \textit{basic} one in \S\ref{basicext}, and is detected by 
the level $[\tau]$. However, the extension cocycle $\omega: \Lambda^2 
L\frg \to \mi\bR$ depends on a linear splitting of the extension 
\begin{equation}\label{exactLie}
0 \to \mi\bR \to L\frg^\tau \to L\frg \to 0.
\end{equation}
For the unique $\frg$-invariant splitting, $S$ is a multiple of $\mi d/dt$. 
Preferred splittings for the twisted loop groups also exist; they are 
discussed in \S\ref{twistedaff}. Hence, subject to topological regularity, 
and using the preferred splittings, the second part of Condition \eqref
{anreg} holds for the entire Lie algebra. Varying the splitting by a \textit{representable} linear map $L\frg\to\mi\bR$, that is, one of 
the form $\eta\mapsto \omega(\xi,\eta)$, changes $S$ by an \textit
{inner} derivation. We assume now that such a splitting has been chosen.

\begin{remark}\label{clarif}\begin{trivlist}\itemsep0ex
\item(i) $S$ must vanish on $\frz$, because the latter exponentiates to a 
torus, over which any $\bT$-extension is trivial. The Heisenberg condition 
allows $\ker S \cap L\frz$ to be no larger. Combining this with the 
discussion of simple summands shows that, for regular $\tau$, $\ker S$ 
is the Lie algebra of a full-rank compact subgroup of $LG$. This is 
the constant copy of $G$, for the standard splitting.
\item(ii) Assuming rotation-equivariance (\S\ref{looprot}), $S$ 
commutes with $d/dt$ and $\Gamma$ with $V$; this justifies our first 
analytic constraint in \eqref{anreg}. However, the constraint truly 
needed to classify admissible representations is weaker: conjugation by 
$\Gamma$ should implement a \textit{representable} change in any splitting 
of (\ref{exactLie}) over $V$. If $\Ad\gamma$ changes the splitting by 
$\eta\mapsto \omega(\xi,\eta)$ for some ($\gamma$-dependent) $\xi\in V$, 
then the alternate copy of $\Gamma$ in $L_fG$, which replaces $\gamma$ 
by $e^{-\xi}\gamma$, decomposes the latter as a semi-direct product of 
$\exp(V)^\tau$ by a central extension of the new $\Gamma$. The new 
$\Gamma$-extension may differ from the original, but has the same 
topological level.
\end{trivlist}\end{remark}

\subsection{Lowest-weight representations.}\label{loweightrep} 
The semi-positive spectral projection of $S$ is an $\omega$-isotropic 
sub-algebra $\frP\subset L\frg_\bC$; we call it the \textit{positive 
polarisation}. The strictly positive part $\frU \subset \frP$ is a Lie 
ideal, and $\ker S\otimes\bC$ is isomorphic to $\frP/\frU$. A linear 
splitting in \S\ref{linsplit} restricts to a \textit{Lie algebra} 
splitting over $\frP$. A \textit{lowest weight} $\tau$-representation of 
$L\frg$ is one generated by an irreducible module of $\ker S$, which is 
killed by the lifted copy of $\ol{\frU}$ in \eqref{exactLie}. 

The lowest-weight condition depends on $S$ and on the splitting of \eqref
{exactLie} over the centre $\frz$. However, if we insist on \textit
{integrability} of the representation to the identity component of the 
loop group (see \S\ref{integrable}), lowest-weight modules are 
irreducible, unitarisable, and their Hilbert space completions are 
unchanged under a representable variation of that splitting.

\subsection{Admissible representations.}\label{admissibility} 
A projective representation of $LG$ is called \textit{admissible} if it 
decomposes as a finite-multiplicity sum of Hilbert space completed 
lowest-weight representations of the Lie algebra. Assuming topological 
regularity, any integrable lowest-weight representation of $L\frg$ 
exponentiates to an action of the identity component of $LG$ on the 
Hilbert space completion. This then induces an admissible representation 
of the full loop group. Moreover, at fixed level, there are finitely 
many irreducibles, up to isomorphism; see \S\ref{repclassifsect}.

There is a similar notion of lowest-weight and admissibility for 
$\cliff(L\frg^*)$-modules, using the same polarisation. (Note that
$\frU$ is $b$-isotropic). As in the finite-dimensional case, there are 
one or two isomorphism classes of lowest-weight representations, according 
to whether $\dim\frg$ is even or odd, and they are irreducible. The 
numbers are switched if we ask for \textit{graded} representations; 
any of the graded irreducibles is called a spin module. The $K$-theory 
of graded, admissible $\cliff(L\frg^*)$-modules (as in \S\ref{untwisted}) 
is $\bZ$, in degrees $\dim G \pmod{2}$. The two spin modules, in the 
even case, differ by parity-reversal, and represent opposite generators 
of $K^0$. (In the odd case, two opposite generators come from the 
two choices of a super-symmetry on the irreducible spin module.)

\begin{remark}
The algebraic approach to representations starts from the Laurent 
polynomial loop algebra $\plg$ and the finite-multiplicity sums of 
integrable lowest-weight modules of $\plg \ltimes \cliff(\plg^*)$. 
These are the Harish-Chandra modules underlying our admissible 
representations. However, as our Dirac construction of $K$-classes 
involves the smooth loop group and its unitary representations, 
we must work more analytically.
\end{remark}

\part{Two examples}
We recall from \cite{fht2} two examples relevant to the construction 
of the Dirac operator families in Chapter~IV, which relate 
representations to $K$-theory classes. The first concerns the group 
$LT$ of loops in a torus; the second is a finite-dimensional Dirac family, 
which leads to an interpretation of our theorem as an infinite-dimensional 
Thom isomorphism.

\section{Spectral flow over a torus} 
\label{spectorus}
\subsection{The circle \cite{atpatsing}.} 
Let $\dirac:=d/d\theta$ be the one-dimensional Dirac operator on the 
complex Hilbert space $\lilb := L^2(S^1;\bC)$, acting as $\mi n$ on 
the Fourier mode $e^{\mi n\theta}$. For any $\xi\in\bR$, the modified 
operator $\dirac_\xi := \dirac + \mi\xi$ has the same eigenvectors, 
but with shifted spectrum $\mi(n+\xi)$. Let $M: \lilb \to \lilb$ 
be the operator of multiplication by $e^{\mi\theta}$. The relation 
$M^{-1}\dirac_\xi M=\dirac_{\xi +1}$ shows that the family 
$\dirac_\xi$, parametrised by $\xi \in \bR$, descends to a family 
of operators on the Hilbert bundle $\bR \times_\bZ \lilb$ over 
$\bR/\bZ$ ($M$ generates the $\bZ$-action on $\lilb$).

Following the spectral decomposition of $\dirac_\xi$, we find that 
one eigenvector crosses over from the negative to the positive imaginary 
spectrum as $\xi$ passes an integer value. Thus, the dimension of the 
positive spectral projection, although infinite, changes by $1$ as we 
travel once around the circle $\bR/\bZ$. This property of the family 
$\dirac_\xi$ is invariant under continuous deformations and captures 
the following topological invariant. Recall \cite{atsing} that the 
interesting component\footnote{The components of essentially positive 
and essentially negative Fredholm operators are contractible.} of the 
space $\mathrm{Fred^{sa}}$ of skew-adjoint Fredholm operators on $\lilb$ 
has the homotopy type of the \textit{small unitary group} $U(\infty )$; 
in particular, $\pi_1\mathrm{Fred^{sa}} = \bZ$. Weak contractibility of 
the \textit{big} unitary group allows us to trivialise our Hilbert bundle 
on $\bR/\bZ$, uniquely up to homotopy; so our family defines a map from 
the circle to $\mathrm{Fred^{sa}}$, up to homotopy. This map detects a 
generator of $\pi_1\mathrm{Fred^{sa}}$. 

\subsection{Generalisation to a torus.} A metric on the Lie algebra 
$\frt$ of a torus $T$ defines the Clifford algebra $\cliff(\frt^*)$, 
generated by the dual $\frt^*$ of $\frt$. Denote by $\psi(\mu)$ the 
Clifford action of $\mu\in\frt^*$ on a complex, graded, irreducible spin
module $\spin(\frt) = \spin^+(\frt)\oplus\spin^-(\frt)$ \cite{abs}. Let 
$\lilb^\pm = L^2(T)\otimes\spin^\pm(\frt)$, denote by $\dirac$ the 
Dirac operator $\sum_a \partial/\partial\theta^a\otimes\psi^a$ on 
$\lilb := \lilb^+ \oplus\lilb^-$, and consider the family of operators 
parametrised by $\mu\in\frt^*$,
\[
\dirac_\mu =\dirac + \mi\psi(\mu): \lilb^+\to\lilb^-. 
\]
Let $\Pi= (2\pi)^{-1}\log 1$ be the \textit{integer lattice} in $\frt$, 
isomorphic to $\pi_1T$. For a weight $\lambda\in\Pi^*:= \mathrm{Hom}
(\Pi;\bZ)$, let $M_\lambda:\lilb\to\lilb$ be the operator of multiplication 
by the associated character $e^{\mi\lambda}: t\mapsto t^{\mi\lambda}$. 
The relation $M_{-\lambda}\circ\dirac_\mu\circ M_{\lambda} = \dirac_{\mu+
\lambda}$ shows that $\dirac_\mu$ descends to a family of fibre-wise 
operators on the Hilbert bundle $\frt^*\times_{\Pi^*}\lilb$ over the dual 
torus $T^* := \frt^*/\Pi^*$. Here, $\Pi^*$ acts on $\frt^*$ by translation 
and on $\lilb$ via the $M$. As before, contractibility of the unitary group 
leads to a continuous family of Fredholm operators over $T^*$. When $\ell 
:= \dim\frt$ is odd, we choose a self-adjoint volume form $\varpi\in
\cliff^1(\frt^*)$. This commutes with all the $\psi^\bullet$ and converts 
$\dirac_\mu$ to a skew-adjoint family $\varpi\cdot\dirac_\mu$ of operators 
acting on $\lilb^+$. Thus, in every case, we obtain a class in $K^\ell(T^*)$. 
This is the $K$-theoretic volume form; more precisely, it is a Fredholm 
model for the Thom push-forward of the identity in $T^*$.

\subsection{Representations and twisted $K$-classes.}\label{fock} 
Relating this construction to our concerns requires a bit more structure, 
in the form of a linear map $\tau:\Pi\to\Pi^* $ (not related to the metric). 
A \textit{central extension} $\Gamma^\tau$ of the product $\Gamma := 
\Pi\times T$ by the circle group $\bT$ is defined by the commutation 
rule
\begin{equation}\label{commute}
ptp^{-1} = t\cdot t^{\mi\tau(p)} \qquad p\in\Pi, t\in T\; \text{ and }\;
t^{\mi\tau(p)}\in\bT. 
\end{equation}
The group $\Gamma^\tau$ has a unitary representation on $L^2(T)$, with 
$T$ acting by translation and $\Pi$ by the $M_{\tau(p)}$'s. If $\tau$ 
has full rank, $L^2(T)$ splits into a finite sum of irreducible $^\tau
\Gamma$-representations $\fock_{[\lambda]}$, each of them comprising 
the weight spaces of $T$ in a fixed residue class $[\lambda ]\in {\Pi^*}
/\tau(\Pi)$. Moreover, these are all the unitary $\tau$-irreducibles 
of $\Gamma$, up to isomorphism. (This will be shown in \S\ref
{repclassifsect}.) 

Now, $\tau$ also induces a map $T \to T^*$, where-under the pull-back of 
$\lilb$ splits, according to the splitting of $L^2(T)$ into the 
$\fock_{[\lambda]}$. Each component carries the lifted Dirac family 
$\dirac_\xi := \dirac + \mi\psi(\tau(\xi))$, descending to a spectral 
flow family over $T$. Except at the single value $\exp(\tau^{-1}
[-\lambda])\in T$ of the parameter, $\dirac_\bullet$ is invertible
on the fibres $\fock_{[\lambda]}\otimes\spin$.
 
All families $\fock_{[\lambda]}\otimes\spin$ have the same image in
$K^\ell(T)$, but this problem is cured by remembering the $T$-action, 
as follows. Instead of viewing the $\dirac_\bullet$ as families over $T$, 
we interpret them as $^\tau\Gamma$-equivariant Fredholm families parametrised 
by $\frt$. Now, $\frt$ is a principal $\Pi$-bundle over the torus $T$, 
equivariant for the trivial action of $T$ on both, and the central 
extension $^\tau\Gamma$ defines a \textit{twisting} for the $T$-equivariant 
$K$-theory of $T$ \cite{fht1}. Classes in $K_T^{\tau+0}(T)$ are then described 
by $\Gamma$-equivariant families of Fredholm operators, parametrised by 
$\frt$, on $\tau$-projective unitary representations of $\Gamma$; twisted 
$K^1$-classes are represented by skew-adjoint families. Thus, our 
families $\dirac_\xi: \fock_{[\lambda]} \otimes \spin^+ \to \fock_
{[\lambda]} \otimes \spin^-$ (respectively $\varpi \cdot\dirac_\xi$ on 
$\fock_{[\lambda]} \otimes\spin^+$ in odd dimensions) give classes in 
$K_T^{\tau+\ell}(T)$. A special case of our main theorem asserts that, when 
$\tau$ is regular, these classes form a $\bZ$-basis of the twisted $K$-groups 
in dimension $\ell\pmod 2$, while the other $K$-groups vanish.

\begin{remark} \label{integration}
The inverse map from $K_T^{\tau+\ell}(T)$ to representations of $\Gamma^\tau$ 
ought to be ``integration over $\frt$" from $K_{\Gamma}^{\tau+\ell}(\frt)$ 
to $R^{\tau+0}(\Gamma)$. This is consistent with our interpretation of our 
main theorem as an infinite-dimensional Thom isomorphism, on the space
of connections over the circle (Chapter V). However, we only know how 
to define the last group in terms of $C^*$-algebras. 
\end{remark}

\subsection{Direct image interpretation.}\label{dirimg} 
Here, we give a topological meaning for the family $(\dirac_\bullet,
\lilb)$; this will be used in \S\ref{peterweylsection}. We claim it 
represents the image of the unit class $[1]$ under the Gysin map 
\[
p_*: K^0(T) \to K^{\tau-\ell}_T(T).
\]
To define $p_*$, we must trivialise the lifted twisting $p^*\tau$. Recall 
that the twisting $\tau$ for the (trivial) $T$-action on $T$ is the 
groupoid defined from the action of $\Gamma^\tau$ on $\frt$. The matching 
model for $p^*\tau$ on $T=\frt/\Pi$ comes from the restricted extension 
$\Pi^\tau$, and this is trivialised by its construction \eqref{commute}. 
The class $[1]$ then corresponds to the trivial line bundle on $\frt$ with 
trivial $\Pi$-action.

We now give an equivalent, but more concrete model for $p_*$. Replace 
$K^*(T)$ by $K^*_T(T\times T)$, where $T$ translates the second factor; 
the projection $P$ to the first factor replaces $p$. If we represent 
$T$ by the $\Gamma$-action groupoid on $\frt\times T$, where $\Pi$ and 
$T$ act by translation on $\frt$, resp.\ $T$, then the twisting $P^*\tau$ 
is represented by the action of $\Gamma^\tau$ on $\frt\times T$.

Call $\cO(\tau)$ the trivial line bundle on $\frt\times T$, but with the 
translation action of $T$ and with $\Pi$-action via the operators $M_{\tau
(\bullet)}$. The two assemble to a $\tau$-action of $\Gamma$, 
so $\cO(\tau)$ gives a class in $K^{\tau+0}_T(T\times T)$. We claim that 
this is the image of $[1]$ under the trivialisation of $p^*\tau$. Indeed, 
our model for $p^*\tau$ as the action of $\Pi^\tau$ on $\frt$ maps to the 
model for $P^*\tau$ by inclusion at $\frt\times\{1\}$; thereunder, 
$\cO(\tau)$ restricts to the trivial bundle with trivial $\Pi$-action.

The Gysin image $P_*[1]$ is now represented by any $\Gamma^\tau$-invariant 
family of Dirac operators on the fibres of $P$, and $(\dirac_\bullet,
\lilb)$ is an example of this.

\subsection{Relation to the loop group $LT$.}\label{loopconnect}
Decompose $LT$ as $\Gamma \times\exp(V)$, where $V=L\frt\ominus \frt$. 
Central extensions of $\exp(V)$ by the circle group $\bT$ are classified 
by skew 2-forms $\omega$ on $V$. We choose a regular such form, in the 
sense of \S\ref{anreg}, together with a positive isotropic subspace 
$\frU \subset V_\bC$. There exists then, up to isomorphism, a unique 
irreducible, unitary projective \textit{Fock representation} $\fock$ of 
$\exp(V)$ which contains a vector annihilated by $\ol {\frU}$. The sum 
of $\Gamma^\tau$ and our extension of $\exp(V)$ is a $\bT$-central 
extension $LT^\tau$ of $LT$, whose irreducible admissible representations 
are isomorphic to the $\fock_{[\lambda]}\otimes\fock$. Our construction 
assigns to each of these a class in $K_T^{\tau+\ell}(T)$.

We will extend this construction and resulting correspondence between 
$LG$-representations and twisted $K$ classes to arbitrary compact groups 
$G$. Observe, by factoring out the space of based loops, that $\Gamma
^\tau$-equivariant objects over $\frt$ are in natural correspondence 
to $LT^\tau$-equivariant ones over the space $\cA$ of $\frt$-valued 
connection forms on the circle, for the gauge action; and it is in this 
form that our construction of the Dirac spectral flow generalises. The 
explicit removal of the Fock factor $\fock$ has no counterpart for 
non-abelian groups, and the same effect is achieved instead by coupling 
the Dirac operator to the spinors on $L\frt/\frt$.
 
\section{A finite-dimensional Dirac family}
\label{finitedirac}
We now recall from \cite{fht2} the finite-dimensional version of our 
construction of twisted $K$-classes from loop group representations 
(Chapter V). For simplicity, we take $G$ to be simple and simply 
connected. Choosing a dominant Weyl chamber in $\frt$ defines the 
nilpotent algebra $\frn$ spanned by positive root vectors, the \textit
{highest root} $\theta$ and the \textit{Weyl vector} $\rho$, the 
half-sum of the positive roots. Roots and weights live in $\frt^*$, a 
weight $\lambda$ defines the character $e^{\mi\lambda}:T\to\bT$, 
sending $e^{\,\xi}\in T$ to $e^{\mi\lambda(\xi)}$.  

The \textit{basic} invariant bilinear form $\langle\ |\ \rangle$ on 
$\frg$ is normalised so that the long roots have square-length 2. Define 
the \textit{structure constants} $f_{ab}^c$ by $[\xi_a,\xi_b] = f_{ab}^c
\xi_c$, in an orthonormal basis $\left\{\xi_a\right\}$ of $\frg$ with 
respect to this bilinear form.\footnote{We use the Einstein summation 
convention, but will also use the metric to raise or lower indexes as 
necessary, when no conflict arises.} Note that $f_{bc}^af_{ad}^c = 
2h^\vee \delta _{bd}$, where $h^\vee = \rho\theta + 1$ is the dual 
Coxeter number. Let $\cliff(\frg^*)$ be the Clifford algebra generated 
by elements $\psi^a$ dual to the basis $\xi_a$, satisfying 
$\psi^a\psi^b + \psi ^b\psi^a = 2\delta ^{ab}$, and let $\spin=\spin^+ 
\oplus \spin^-$ be a graded, irreducible complex module for it. This 
is unique up to isomorphism and (if $\dim\frg$ is even) up to parity switch. There is a unique action of $\frg$ on $\spin$ 
compatible with the adjoint action on $\cliff(\frg^*)$; the action of 
$\xi_a$ can be expressed in terms of Clifford generators as
\[
\sigma_a = -\frac{1}{4}\: f_{bc}^a \cdot \psi^b\psi^c.
\]
It follows from the Weyl character formula that $\spin$ is a sum of $2^{\lceil\dim \frt/2\rceil}$ copies of the irreducible representation 
$V_{-\rho}$ of $\frg$ of lowest weight $(-\rho)$. The lowest-weight 
space is a graded $\cliff(\frt^*)$-module; for dimensional reasons,  
it is irreducible.

\subsection{The Dirac operator.} Having trivialised the Clifford and 
Spinor bundles over $G$ by left translation, consider the following 
operator on spinors, called by Kostant \cite{kost1} the ``cubic Dirac 
operator'':
\begin{equation} \label{findirac}
\dirac = 
	R_a\otimes \psi ^a + \frac{1}{3}\:\sigma _a\cdot \psi^a=
	R_a\otimes\psi ^a - \frac{1}{12}\:f_{abc}\psi ^a
	\psi ^b\psi ^c,
\end{equation}
where $R_a$ denotes the right translation action of $\xi_a$ on functions.
Let also $T_a=R_a+\sigma_a$ be the total right translation action of
$\xi_a$ on smooth spinors. 

\begin{proposition} \label{diracrels}
 $\left[\dirac,\psi^b \right]= 2 T_b$;\ \ $\left[\dirac,T_b \right]=0$.
\end{proposition}

\begin{proof} The second identity expresses the right-invariance of
the operator, while the first one follows by direct computation:

\begin{minipage}[b]{14.5cm}
\begin{eqnarray*}
\left[ \dirac,\psi ^b \right] &=& 
	R_a \otimes \left[\psi^a, \psi ^b \right] + 
\frac{\sigma_a} {3}\left[ {\psi ^a, \psi ^b} \right] - 
\frac{1}{3}\left[ {\sigma_a,\psi^b} \right]\cdot \psi^a \\
&=& 2 R_b + \frac{2}{3}\:\sigma_b - \frac{1}{3}\: 
f_{ca}^b\cdot \psi^c\psi^a\\
&=& 2(R_b + \sigma_b).
\end{eqnarray*}
\end{minipage}
\end{proof}

\subsection{The Laplacian.} The Peter-Weyl theorem decomposes $L^2(G;
\spin)$ as $\bigoplus_\lambda V_{-\lambda}^* \otimes V_{-\lambda} 
\otimes \spin$, where the sum ranges over the dominant weights $\lambda$ 
of $\frg$. Left translation acts on the left, $R_a$ on the middle 
and $\sigma_a$ on the right factor. Hence, $\dirac$ acts on the two 
right factors alone. As a consequence of \eqref{diracrels}, the Dirac 
Laplacian $\dirac^2$ commutes with the operators $T_\bullet$ and 
$\psi^\bullet$. As these generate $V_{-\lambda}\otimes \spin$ from 
its $-(\lambda+\rho )$-weight space, $\dirac^2$ is determined 
from its action there. To understand this action, rewrite 
$\dirac$ in a root basis of $\frg$,  

\begin{equation} \label{rootexpress}
\dirac=R_j\otimes \psi^j + \frac{1}{3}\:\sigma_j \psi^j + 
R_\alpha\otimes \psi^{-\alpha } + R_{-\alpha}\otimes \psi^\alpha + 
\frac{1}{3}\left(\sigma_\alpha\psi^{-\alpha } + \sigma_{-\alpha}
\psi^\alpha \right),
\end{equation}
where the $j$'s label a basis of $\frt$ and $\alpha$ ranges over the 
positive roots. The commutation relation $\left[\sigma_{-\alpha}, 
\psi^\alpha\right] = \psi(-2\mi\rho)$, where summation over $\alpha$ 
has been implied, converts \eqref{rootexpress} to
\[
\dirac = R_j\otimes \psi^j + \frac{1}{3}\:\sigma_j\psi^j 
- \frac{2\mi}{3}\:\psi(\rho) + R_\alpha \otimes \psi^{-\alpha} 
+ R_{-\alpha}\otimes \psi^\alpha + \frac{1}{3} 
	\left({\sigma_\alpha\psi^{-\alpha} 
	+ \psi^\alpha\sigma_{-\alpha}}\right),
\] 
and the vanishing of all $\alpha$-terms on the lowest weight space 
leads to the following 

\begin{proposition} \label{diracsquare}
\begin{trivlist}\itemsep -.1ex
\item{\upshape (i)} $\dirac= -\mi\psi (\lambda +\rho )$ on the
$-(\lambda+\rho )$-weight space of $V_{-\lambda}\otimes \spin$. 
\item{\upshape (ii)} $\dirac^2=-(\lambda +\rho )^2$ on 
$V_{-\lambda}\otimes \spin$. \qed
\end{trivlist}
\end{proposition}

\subsection{The Dirac family.} Consider now the family $\dirac_\mu:= 
\dirac +\mi\psi(\mu)$, parametrised by $\mu\in\frg^*$. Conjugation by a 
suitable group element brings $\mu$ into the \textit{dominant chamber} 
of $\frt^*$. From \eqref{diracsquare}, we obtain the following 
relations, where $\langle T|\mu \rangle$ represents the contraction of 
$\mu$ with $T \in \frg^* \otimes\mathrm{End} \left(V\otimes\spin\right)$, 
in the basic bilinear form (the calculation is left to the reader). 

\begin{corollary} 
\begin{trivlist}\itemsep-.15ex
\item{\upshape (i)} $\dirac_\mu = \mi\psi (\mu-\lambda-\rho)$
on the lowest weight space of $V_{-\lambda}\otimes\spin$.
\item{\upshape (ii)} $\dirac_\mu^2=-(\lambda +\rho -\mu)^2 + 2\mi\langle
T|\mu \rangle - 2 \langle \lambda +\rho |\mu \rangle$. \qed
\end{trivlist}
\end{corollary} 

\subsection{The kernel.} Because $\mi\langle T|\mu \rangle 
\le \langle \lambda +\rho|\mu \rangle$, with equality only on the 
$-(\lambda+\rho)$-weight space, $\dirac_\mu$ is invertible on $V_
{-\lambda} \otimes\spin$, except when $\mu$ is in the co-adjoint orbit 
$\frO$ of $(\lambda +\rho)$. In that case, the kernel at $\mu \in \frg^*$ 
is that very weight space, with respect to the Cartan sub-algebra 
$\frt_\mu$ and dominant chamber \textit{defined by the regular element} 
$\mu$. This is the lowest-weight line of $V_{-\lambda}$ tensored with the 
lowest-weight space of $\spin$, and is an irreducible module for the 
Clifford algebra generated by the normal space $\frt_\mu^*$ to $\frO$ at 
$\mu$. More precisely, the kernels over $\frO$ assemble to the normal 
spinor bundle to $\frO\subset\frg^*$, twisted by the natural line bundle 
$\cO(-\lambda -\rho)$. Finally, at a nearby point $\mu +\nu$, with $\nu 
\in \frt_\mu ^*$, $\dirac_{\mu +\nu}$ acts on $\ker(\dirac_\mu)$ as 
$\mi\psi(\nu)$. 

\subsection{Topological interpretation.} The family of operators $\dirac_
\mu$ on $V_{-\lambda}\otimes \spin$ is a compactly supported $K$-cocycle 
on $\frg^*$, equivariant for the co-adjoint action of $G$. As before, 
when $\dim \frg$ is odd, we use the volume form $\varpi$ to produce the 
skew-adjoint family $\varpi \dirac_\mu$, which represents a class in 
$K^1_G$. Our computation of the kernel identifies these classes with the 
Thom classes of $\frO\subset\frg^*$, with coefficients in the natural 
line bundle $\cO(-\lambda-\rho)$. Sending $V_{-\lambda}$ to this class 
defines a linear map 

\begin{equation} \label{orbitthom}
R(G) \to K_G^{\dim\frg}(\frg^*).
\end{equation}

There is another way to identify this map. Deform $\dirac_\mu$ to 
$\mi\psi(\mu)$ via the (compactly supported Fredholm) family $\vep
\cdot\dirac +\mi\psi(\mu)$. At $\vep=0$ we obtain the standard Thom class 
of the origin in $\frg^*$, coupled to $V_{-\lambda}$. Therefore, our 
construction is an alternative rigid implementation of the Thom 
isomorphism $K_G^0(0) \cong K_G^{\dim\frg} (\frg^*)$.

The inverse isomorphism is the push-forward from $\frg^*$ to a point. In 
view of our discussion, this expresses $V_{-\lambda}$ as the Dirac index 
of $\cO(-\lambda -\rho)$ over $\frO$, leading to (the Dirac index version 
of) the Borel-Weil-Bott theorem. The affine analogue of the Thom 
isomorphism (\ref{orbitthom}) is Theorem 3, equating the module of 
admissible projective representations with a twisted $K_G(G)$.

\subsection{Application to Dirac induction.} \label{dirind}
For later use, we record here 
the following proposition; when combined with the Thom isomorphisms and 
the resulting twists, it gives the correct version of Dirac induction for 
\textit{any compact Lie group} $G$ (not necessarily connected). Let 
$N\subset G$ be the normaliser of the maximal torus $T$. We have a 
restriction map $K_G(\frg^*)\to K_N(\frt^*)$ and an ``induction" 
$K_N(\frt^*) \to K_G(\frg^*)$ (Thom push-forward from $\frt^*$ to 
$\frg^*$, followed by Dirac induction from $N$ to $G$). 

\begin{proposition}\label{diracbwb}
The composition $K_G(\frg^*)\to K_N(\frt^*) \to K_G(\frg^*)$ is 
the identity.
\end{proposition}
\begin{proof}
Express the middle term as $K_G(G\times_N\frt^*)$, with the left action 
of $G$ on the induced space. The map from $G\times_N\frt^*$ to $\frg^*$ 
sends $(g,\mu)$ to $g\mu g^{-1}$. Since $N$ meets every component of $G$ 
(Prop.~\ref{diagaut}), this map is a diffeomorphism over regular 
points. Every class in $K_G(\frg^*)$ is the Thom push-forward of a 
class $[V]\in K_G(0)$. Deforming this to $\dirac + \mi\psi(\mu)$ 
leads to a class supported on a \textit{regular} orbit; \textit{a 
fortiori}, our composition is the identity on such classes, hence on 
the entire $K$-group.
\end{proof}

\part{Computation of twisted $K_G(G)$}
In this chapter, we compute the twisted $K$-theory $K^\tau_G(G)$ by 
topological methods, for arbitrary compact Lie groups $G$ and regular 
twistings $\tau$. A key step is the reduction to the maximal torus, 
Proposition~\ref{split}. Our answer takes the form of a twisted 
$K$-theory of the set of regular affine weights at level $\tau$, 
equivariant under the \textit{extended affine Weyl group} (\S\ref
{simplewaff}, \S\ref{affaction}). This action has finite quotient 
and finite stabilisers, and the $K^\tau$-theory is a free abelian group 
of finite rank. 

For a detailed discussion of foundational questions on twisted $K$-theory, 
we refer to \cite{fht1} and the references therein.

\section{A ``Mackey decomposition" lemma}
\label{general}

We recall from \cite{fht1} the following construction and generalisation 
of Lemma~2.14 in \cite{fht3}, which will be a key step in our computation of 
$K_G(G)$. This is a topological form of the Mackey decomposition of 
irreducible representations of a group, restricted to a normal subgroup; 
the analogy will be particularly relevant in \S\ref{normal}.

\subsection{Construction.} \label{constr}
Let $H$ be a compact group, acting on a compact Hausdorff space 
$X$, $\tau$ a twisting for $H$-equivariant $K$-theory, $M \subseteq H$ a 
normal subgroup acting trivially on $X$. The following data can be 
extracted from this:
\begin{enumerate}
\itemsep.1ex
\item an $H$-equivariant family, parametrised by $X$, of $\bT$-central 
extensions $M^\tau$ of $M$; 
\item an $H/M$-equivariant covering space $p:Y \to X$, whose fibres label 
the isomorphism classes of irreducible, $\tau$-projective representations 
of $M$; 
\item an $H$-equivariant, tautological projective bundle $\bP R\to Y$, 
whose fibre $\bP R_y$ at $y\in Y$ is the projectivised $\tau$-representation 
of $M$ labelled by $y$;
\item a class $[R] \in K^{\bP R}_H(Y)$, represented by $R$;
\item a twisting $\tau'$ for the $H/M$-equivariant $K$-theory of $Y$, 
and an isomorphism of $H$-equivariant twistings $\tau' \cong p^*\tau - \bP R$.
\end{enumerate}

\noindent Items (iii) and (v) are only defined up to canonical isomorphism. 
Note that, if $M^\tau$ is abelian, as will be the case in our application,
then $\bP R = Y$, which can be taken to represent the zero twisting. 
However, $[R]$ is \textit{not} the identity class $[1]$, because of the 
non-trivial $M$-action on the fibres.

\begin{lemma}[Key Lemma] \label{keylemma}
The twisted K-theories $K^{\tau'}_{H/M}(Y)$ and $K^\tau_H(X)$  
are naturally isomorphic. \qed
\end{lemma}

\noindent Recall that the isomorphism is induced by the composition below:
\begin{equation}\label{chain}
K_{H/M}^{\tau'}(Y) \longrightarrow K_{H}^{\tau'}(Y)
\cong K_{H}^{p^*\tau-\bP R}(Y) \xrightarrow{\otimes[R]\:} 
K_{H}^{p^*\tau}(Y) \xrightarrow{\:p_!\:} K^\tau_H(X).
\end{equation}
The inverse map lifts a class from $K^\tau_H(X)$ to $K_H^{p^*\tau}(Y)$, 
tensors with the dual class $[R^\vee]\in K^{-\bP R}_H(Y)$ and, finally, 
extracts the $M$-invariant part. The last step requires a bit of care, 
as it involves a Morita isomorphism associated to two different 
projective bundle models for the same twisting.

\section{Computation when the identity component is a torus}
\label{topktorus}
To ensure consistency of notation when the identity component $G_1$ is a 
torus $T$, we write $N$ for $G$ and $W$ for $\pi_0N$. Denoting, for any 
$f\in N$, by $N(f)$ the stabiliser in $N$ of the component $fT$, we can 
decompose $K^\tau_N(N)$ as a sum over representatives $f\in N$ of the 
conjugacy classes in $W$:
\begin{equation} \label{decomp}
K^\tau_N(N) \cong \bigoplus\nolimits_f K^\tau_{N(f)}\left(fT\right).
\end{equation}

\subsection{The identity component.} \label{idcomp}
With $H=N$ and $M=X=T$ in construction \ref{constr}, Lemma \ref
{keylemma} gives $K^\tau_N(T)= K_W^{\tau'}(Y)$. It is easy to 
describe the bundle $p:Y\to T$. A twisting class $[\tau] \in H_N^3(T)$ 
restricts to $H_T^3(T)$, hence to $H^1(T)\otimes H^2_T$, and contraction 
with the first factor gives a map $\kappa^\tau: H_1(T)\to H^2_T$. This gives 
a translation action of $\Pi :=\pi_1T$ on the set $\wtlat^\tau$ of $\tau
$-affine weights of $T$, and $Y$ is the associated bundle $\frt\times_\Pi 
\wtlat^\tau$. If $\kappa^\tau$ is injective, as per our regularity condition 
\eqref{topreg}, $Y$ is a union of copies of $\frt$, labelled by $\wtlat
^\tau\left/{\kappa^\tau(\Pi)}\right.$, and integration along $\frt$ gives 
\[
K_W^{\tau'}(Y) = K_W^{\tau'-\sigma(\frt)-\dim T}
	\left(\wtlat^\tau\left/{\kappa^\tau (\Pi)}\right.\right), 
\]
where the down-shift $\sigma(\frt)$ in the twisting is defined by a 
$W$-equivariant Thom class of $\frt$, represented by a choice of spinors 
$\spin(\frt)$ with projective $W$-action.  

\subsection{Affine Weyl action.} \label{simplewaff}
We restate this by observing that the class $[\tau]\in H^3_N(T)$ has a 
``leading term" in $H^1_W(T;H^2_T)$, with respect to the Hochschild-Serre 
spectral sequence $E_2^{p,q}=H^p_W(T;H^q_T) \Rightarrow H^{p+q}_N(T)$. 
This term captures the $W$-action on the covering $Y$ of $T$, but, more 
importantly, defines an affine action on $\wtlat^\tau$ of the \textit
{extended affine Weyl group} $W\ltimes\Pi$, extending the action of 
$\Pi$. Comparing orbits and stabilisers gives an equivalence of 
categories of equivariant bundles, and hence an isomorphism of the 
desired form,
\[
K_W^{\tau'-\sigma(\frt)}\left({\wtlat^\tau}\left/
	{\kappa^\tau (\Pi)}\right.\right) = K_{W \ltimes \Pi}
		^{\tau'-\sigma(\frt)}\left({\wtlat^\tau} \right).
\]

\subsection{A general component.} \label{affaction} 
Let now $T^f$ be the $T$-centraliser of $f\in N$ and $\ul T$ its identity 
component. Then, $fT$ is a homogeneous space, with discrete isotropy, for 
the combined action of ${N(f)}\left/{\ul T}\right.$ by conjugation and 
of $\ul\frt :=\frt^f$ by translation. We thus have an $N(f)\ltimes
\ul\frt$-isomorphism 
\begin{equation}\label{homog}
fT \cong \left.\left[\left(N(f)\left/{\ul T}\right.
\right) \ltimes \ul\frt\right]\right/\waff\,; 
\end{equation}
the stabiliser $\waff$ of $f$ is expressed, by projection to 
$N(f)/\ul T$, as a group extension 
\begin{equation}\label{wfext}
1 \to \ul\Pi \to \waff \to \widetilde{W}^f\to 1, 
\end{equation}
where $\ul\Pi:=\pi_1\ul T$ and $\widetilde{W}^f:= [N(f)/\ul T]^f$ is itself 
an extension of $W^f$ by the finite group $[T/\ul T]^f$: 
\[
1 \to [T/\ul T]^f \to \widetilde{W}^f \to W^f \to 1.
\] 
Exactness on the right follows from the vanishing of $H^1_{\langle f 
\rangle}(T/\ul T)$; that, in turn, follows from the absence of $f
$-invariants in $\pi_1(T/\ul T)$.

With $X=fT$, $H=N(f)$ and $M = \ul T$ in \eqref{constr}, an $N(f)
$-equivariant twisting $\tau$ defines a covering space $Y\to fT$, with 
fibres the sets $\ul{\wtlat}^\tau$ of $\tau$-affine weights of $\ul T$. 
Via \eqref{homog}, this cover is associated to an affine action of 
$\waff$ on $\ul{\wtlat}^\tau$, which is classified by the leading 
component of $[\tau]\in H^3_{N(f)}(fT)$ in 
\begin{equation}\label{leading}
H^1_{N(f)/\ul T}(fT;H^2_{\ul T}) = H^1_\waff(H^2_{\ul T}).
\end{equation}

\begin{theorem} \label{toruskcompute}
\begin{trivlist}\itemsep0ex
\item \textup{(i)} We have a natural isomorphism $K^\tau_{N(f)}
\left(fT\right) = K_\waff^{\tau'}\left(\ul{\wtlat}^\tau \times 
\ul{\frt} \right)$. 
\item \textup{(ii)} If $\tau$ is regular, this is also $K_\waff^
{\tau'-\sigma(\ul\frt)-\dim\ul\frt}\left(\ul{\wtlat}^\tau\right)$, 
and is free, of finite rank over $\bZ$. 
\end{trivlist}
\end{theorem} 
\begin{proof}
The first part is Lemma \ref{keylemma}. Provided that all stabilisers 
of $\waff$ on $\ul{\wtlat}^\tau$ are finite, part (ii) follows from (i) by 
integration along $\ul\frt$, and $\sigma(\ul\frt)$ is the twisting of 
the equivariant Thom class. 

Now, $\ul\Pi\subset\waff$ has finite index, and acts on $\ul{\wtlat}^\tau$ 
by translation, via the linear map $\kappa^\tau: \ul\Pi\to \ul\wtlat$, 
defined by restricting $[\tau]$ to $H^3_{\ul T}(\ul T)$. Topological 
regularity of $\tau$ implies finiteness of the quotient $\ul{\wtlat}^\tau/
\waff$ and of all stabilisers.
\end{proof}

\begin{remark}\label{otherpicture}
\begin{trivlist}\itemsep0ex
\item (i) Considering the action of $N(f)\ltimes\ul\frt$ on $fT$ 
leads to the presentation 
\[
fT \cong N(f)\ltimes\ul\frt /\nafff
\]
where the stabiliser $\nafff$ of $f$ fits now in an extension 
$1 \to \ul T\to \nafff \to \waff \to 1$. 
\item (ii) Without Lemma \ref{keylemma}, the isomorphism \eqref{homog} 
identifies $K^\tau_{N(f)}(fT)$ with $K^\tau_{\nafff}(\ul{\frt})$; the 
right-hand side has a sensible topological interpretation, because the group 
action is proper, resulting in a stack of the type studied in \cite{fht1}. 
It is tempting to integrate along $\ul\frt$, to land in the 
$\nafff$-equivariant twisted $K$-theory of a point. However, no 
\textit{topological} definition of $K$-theory that we know allows this 
operation (cf.\ Remark~\ref{integration}); this could perhaps be done 
by $C^*$-algebra methods.
\item (iii) For a loop group interpretation 
of $\nafff$, $\waff$ and its action on $\ul{\wtlat}^\tau$, see Remark 
\ref{philoop} below.
\end{trivlist}\end{remark}

\subsection{Induction from conjugacy classes.} \label{inductorb}
The following result (with Theorem~\ref{split} in the next section) is 
the basis for our original construction \cite{freedicm} of twisted 
$K$-classes. For each element of the natural basis of Theorem~\ref 
{toruskcompute}.ii, it selects a distinguished $N(f)$-conjugacy class 
in $fT$.\footnote{Up to an overall ambiguity, see Remark \ref{ambiguity}.} 
We shall revisit this when discussing the Dirac families in Chapter~\ref 
{diracfam}.   

\begin{proposition} \label{orbits}
If $\tau$ is regular, $K^\tau_{N(f)}(fT)$ is spanned by classes 
supported on single $N(f)$-orbits.
\end{proposition}

\begin{proof} An affine action of $\waff$ on $\ul\frt$ is inherited from 
the conjugation$\ltimes$translation action of the ambient group $\left(
N(f)/\ul{T}\right)\ltimes \ul\frt$. There is also a $\waff$-action on 
the affine copy $\ul{\wtlat}^\tau\otimes\bR$ of $\ul\frt^*$, defined above 
by $[\tau]$. Such actions are classified by the groups
\[
H_\waff^1(\ul{\frt}) \cong 
		\mathrm{Hom}_{\widetilde{W}^f}(\ul\Pi,\ul{\frt})
\quad\mbox{and}\quad 
H_\waff^1(\ul {\frt}^*) \cong 
		\mathrm{Hom}_{\widetilde{W}^f}(\ul\Pi, \ul{\frt}^*),
\] 
respectively; $\widetilde{W}^f$ acts by conjugation. The first class is 
the natural map $\ul\Pi\to \ul\frt$; the second, the map $\kappa^\tau \otimes 
\bR$. Hence, the two actions of $\waff$ are isomorphic by some 
translate $\kappa^\tau_\nu: \ul\frt \to \ul\frt^*$ of $\kappa^\tau\otimes\bR$.

A class in $K^{\tau'-\sigma(\ul\frt)}_\waff(\ul{\wtlat}^\tau)$ can 
be pushed forward to $K^{\tau'}_{\waff}(\ul{\wtlat}^\tau \times 
\ul{\frt})$ using the graph of the inverse map $(\kappa^\tau_\nu)^{-1}$. 
Under (\ref{toruskcompute}.i), its image in $K^\tau_{N(f)}(fT)$ is supported 
on a single conjugacy class, if the original lived on a single $\waff
$-orbit. 
\end{proof} 

\begin{remark}\label{ambiguity} 
\begin{trivlist}\itemsep0ex 
\item (i) $\kappa^\tau_\nu$ descends to an affine isogeny from $f\ul T$ to 
$\ul{\wtlat}^\tau\otimes\bT$, preserving the actions of $\widetilde{W}^f$.
The quotient spaces $f\ul T/\widetilde{W}^f$ and $fT/N(f)$ are isomorphic,
and the conjugacy classes in Prop.~\ref{orbits} lie in the fibre of this 
isogeny over the base-point $\ul{\wtlat}^\tau$ of the second torus.\footnote
{Note that $f\ul T$ covers $fT/T$, which in turn surjects onto the conjugacy 
classes.} Specifically, a class in the (twisted) $K_\waff(\ul{\wtlat}^\tau)$ 
supported on a $\waff$-orbit $\Omega$ corresponds to one in $K^\tau_{N(f)}
(fT)$ with support at the single $\widetilde{W}^f$-orbit $f\exp\left(
(\kappa^\tau_\nu)^{-1}\Omega\right)$.
\item(ii) An ambiguity in the set of orbits results from our freedom in 
identifying the affine spaces $\ul\frt$ and $\ul\frt^*$: we are free to 
translate by the $W^f$-invariant part of $\ul\frt$. 
\end{trivlist}\end{remark}

\subsection{Relation to loop groups.}\label{philoop}
The isomorphism \eqref{homog} implies that $\waff$ is $\pi_1$ of the 
homotopy quotient of $fT$ by $N(f)$. The latter is equivalent to the 
classifying space $BL_fN$: this is best revealed by the gauge action 
of $L_fN$ on the contractible space of connections on the principal 
$N$-bundle over the circle, with holonomies in $[fT]$; fixing the fibre 
over a base-point, the space of holonomies becomes $fT$, while the 
residual symmetry group is $N(f)$. All in all, $\waff = \pi_0L_fN$.

Every component of $L_fN$ contains loops of minimal length, so the 
subgroup $\Gamma_fN \subset L_fN$ of $f$-twisted \textit{geodesic loops} 
is an extension of $\waff$ by $\ul T$. This $\Gamma_f N$ is in fact 
isomorphic to $\nafff$: to equate them, interpret the presentation of 
$fT$ in Remark~\ref{otherpicture}.i as the quotient of $N(f)\times\ul\frt$, 
the set of flat bundles over the interval based at the endpoints and with 
\textit{constant} connection forms, under the gauge action of $\Gamma_fN$. 

The action of $\waff$ on $\ul{\wtlat}^\tau$ and its twisting $\tau'$ also 
have a loop group description. The connection picture above gives an 
equivalence between the smooth groupoids associated to the actions of 
$N(f)$ on $fT$ and of $\nafff$ on $\ul\frt$ (via $\waff$). We are 
interested in twistings coming from $\bT$-central extensions of 
$\nafff$. In that case, the action of $\waff$ on $\ul{\wtlat}^\tau$ 
arises from the conjugation action of $\nafff$ on $\ul T$, in the 
central extension. The subgroup of $(\naff^\tau)_\lambda$ stabilising 
a weight $\lambda\in \ul{\wtlat}^\tau$ is an extension of $(\waff)_\lambda$ 
by $\ul{T}^\tau$; pushing out via $\lambda$ gives a $\bT$-central extension 
of $(\waff)_\lambda$, and these extensions assemble to the twisting $\tau'$. 

\section{General compact groups}
\label{topkcompact}

For any compact $G$, we will describe $K^\tau_G(G)$ in terms of the 
maximal torus $T$ of $G$ and its normaliser $N$. We must first recall 
some facts about disconnected groups; readers focusing on the connected 
case may skip ahead to \S\ref{weylmap}. We keep the notations of 
\S\ref{topktorus}.

\subsection{Diagram automorphisms.} Choose a set of simple root vectors 
in $\frg_\bC$, satisfying, along with their conjugates and the simple 
co-roots, the standard $\mathfrak{sl}_2$ relations. 
 
\begin{proposition} \label{diagaut}
Every outer automorphism of $\frg$ has a distinguished implementation,
called \textit{diagram automorphism}, which preserves $\frt$ and its 
dominant chamber and permutes the simple root vectors. 
\end{proposition}

\begin{proof} The variety $G/N$ of Cartan sub-algebras in $\frg$ has the 
rational cohomology of a point, so any automorphism of $\frg$ fixes a 
Cartan sub-algebra, by the Lefschetz theorem. Composing with a suitable 
inner automorphism ensures that we preserve $\frt$ and the dominant 
chamber. Conjugation by $T$ provides the freedom needed to permute the 
simple root vectors without scaling.
\end{proof}

\begin{remark}
The use of Lefschetz's theorem is justified, because any automorphism  
of $\frg$ lifts to the orientable double cover of $G/N$, which is 
the variety of oriented Cartan sub-algebras.
\end{remark}
 
\begin{corollary} \label{reducetocenter}
$G$ can be reduced to an extension of $\pi_0G$ by the centre of the 
identity component $G_1$. 
\end{corollary}

\begin{proof} The subgroup of $G$-elements whose Ad-action on $\frg$ is 
a diagram automorphisms meets every component of $G$, and meets $G_1$ in 
its centre. This is our reduction.
\end{proof}

\subsection{Conjugacy classes in $G$.}\label{specialf} The push-out of 
\eqref{reducetocenter} to a maximal torus $T$ is called a \textit
{quasi-torus} $Q_T \subset G$; it meets every component of $G$ in a 
translate of $T$. $Q_T$ depends on $T$ and a choice of dominant chamber. 
Choose $f\in Q_T$; its Ad-action on the dominant chamber must fix some 
interior points, so $\ul\frt=\frt^f$ contains $\frg$-regular elements. 
The identity component $\ul T$ of the invariant subgroup $T^f$ is then 
a maximal torus of the centraliser $G^f$ of $f$. 

Call $W = N/T$, $W_1 = (N\cap G_1)/T$ the Weyl groups of $G$ and $G_1$; 
we have $W = \pi_0G \ltimes W_1$, by \eqref{reducetocenter}. Call $[f]$ 
the image of $f$ in the quotient $fT/T$ by $T$-conjugation. Conjugation 
by $N(f)$, the subgroup of $N$ preserving the component $fT$, descends 
to an action of the group $W^f=\pi_0N(f)$ on $fT/T$. Let $\ul W:=
W^f\cap W_1$, and $\ul{\widetilde W}$ its extension by $[T/\ul T]^f$ 
restricted from the $\widetilde{W}^f$ of \eqref{wfext}.

\begin{lemma} \label{gentwistedconj}
\begin{trivlist}\itemsep0ex
\item (i) The space of conjugacy classes $fG_1/G_1$ is $(fT/T)\left/
\right.$. 
\item (ii) The Weyl group of $G_1^f$ is an extension by $\pi_0T^f$ of 
the $\ul W$-stabiliser of $[f]$. 
\end{trivlist}
\end{lemma}

\begin{proof}
Part (i) reformulates Theorem ?.? of \cite{td}: indeed, $fT/T$ is the 
quotient of $f\ul T$ under conjugation by $[T/\ul T]^f$, whence it follows 
that $(fT/T)\left/\ul W\right. \cong f\ul T/\ul {\widetilde{W}}$. That 
is the description in \cite{td}.
 
The normaliser of $\ul T$ in $G_1^f$ is $N\cap G_1^f$, by regularity,  
and exactness of $1 \to T^f \to N\cap G_1^f \to \ul W$ implies (ii). 
\end{proof}

\begin{remark}
Translation by $f$ identifies $fT/T$ with the co-invariant torus $T_f$ 
(quotient of $T$ by the sub-torus $\{xfx^{-1}f^{-1}\,|\, x\in T\}$). The 
$\ul W$-action on $fT/T$ is \textit{affine} under the quotient $\ul W$-action 
on $T_f$. However, the two $\ul W$-actions \textit{agree} when $f$ is a 
diagram automorphism $\vep$: $\ul W$ indeed isomorphic to the Weyl group 
of $\frg^\vep$ (Appendix \ref{afftwist}).
\end{remark} 

\subsection{The Weyl map $\omega$.} \label{weylmap}
Decompose $K^\tau_G(G) = \bigoplus\nolimits_f K^\tau_{G(f)}\left
(fG_1 \right)$ over representatives $f\in Q_T$ of conjugacy classes in 
$\pi_0$; $G(f)$ denotes the stabiliser of the component $fG_1$. The 
$G(f)$-equivariant map 
\[
\omega : G(f)\times_{N(f)}fT \to fG_1,\qquad 
g\times ft\mapsto g\cdot ft\cdot g^{-1} 
\]
induces two morphisms in twisted $K$-theory, \textit{restriction} 
$\omega^*$ and \textit{induction}\footnote{The names are justified in 
\S\ref{semiinf}.} $\omega_*$:
\[
K^\tau_{N(f)}\left({fT} \right) \cong K^\tau_{G(f)}\left( 
{G(f)\times_{N(f)} fT}\right) \begin{array}{c}{\raisebox{-.7ex}{
$\stackrel{\omega_*}{\displaystyle\longrightarrow}$}}
\\ {\stackrel{\displaystyle\longleftarrow}
{\scriptstyle\omega^*}} \end{array} K^\tau_{G(f)}\left({fG_1}\right). 
\]

\begin{theorem} \label{split}
The composition $\omega_*\circ \omega^*$ is the identity.
\end{theorem}
\noindent Consequently, $K^\tau_{G(f)}(fG_1)$ is a summand in $K^\tau
_{N(f)}(fT)$, split as an $R(G)$-module. To identify it, we will call 
a weight in $\ul{\wtlat}^\tau$ \textit{regular} if it corresponds to a 
regular conjugacy class in $fG_1$, under the isomorphism $\kappa^\tau_\nu: 
\ul\frt \to \ul\frt^*$ in the proof of \eqref{orbits}. Clearly, this 
condition is preserved by $\waff$.

\begin{theorem} \label{kcompute}
$K^\tau_{G(f)}(fG_1)$ is the summand in $K^\tau_{N(f)}(fT)$ corresponding 
to the regular weights: 
\[
K^{\tau}_{G(f)}(fG_1) = K_\waff^{\tau'-\sigma(\ul\frt)-\dim T^f}
\left(\ul{\wtlat}^\tau_{reg}\right).
\]
\end{theorem} 

\begin{remark}\label{waffdef}
$\waff$ is called the \textit{$f$-twisted, extended affine Weyl group} of 
$G$. Regular weights are those \textit{not} fixed by any \textit{affine 
reflection} of $L_f\frg$ (\S\ref{morewaff}, \S\ref{afroot}). Reaching 
ahead a bit, this can be seen as follows. The action of $\naff$ 
on $\ul{\wtlat}^\tau$ (Remark \ref{philoop}) is part of the co-adjoint 
action of $L_fG$ on an affine slice (at level $[\tau]$) in $L_f(\frg^*)
^\tau$ (\S\ref{coadjoint}). The latter is isomorphic the space of connections 
over the circle with holonomy in $fG_1$, with the gauge action. That, 
in turn, is equivalent to the conjugation $G(f)$-action on $fG_1$. Regular 
are those weights whose co-adjoint Lie algebra stabiliser is minimal, 
namely $\ul\frt$. Singular weights will be fixed by some Weyl reflection 
in their $L_f\frg$-stabiliser, which is an affine Weyl reflection. Some 
facts about simple algebras are recalled in Appendix~\ref{afftwist}.
\end{remark}

\begin{proof}[Proof of (\ref{split})] 
The quotient spaces ${fT}\left/N(f)\right.$ and ${fG_1}\left/G(f)\right.$ 
are isomorphic under $\omega$ (Lemma \ref{gentwistedconj}). We shall show 
that $\omega_*\circ\omega^*$ is the identity on small neighbourhoods of 
conjugacy classes: a Mayer-Vietoris argument then implies that the map is 
a global isomorphism. However, $K^\tau_{N(f)}(fT)$ is spanned by classes 
induced from single orbits (Prop.~\ref{orbits}). Their $\omega_*$-images 
are fixed by $\omega_*\circ\omega^*$, so the theorem follows.

We need a local model for the Weyl map. We work near $f$, which was 
arbitrary in $Q_T$. Because $\ul T$ contains regular elements, $N^f:= 
N \cap G^f$ is the normaliser of $\ul T$ in $G^f$. Now, the translate 
$f\cdot\exp(\frg^f)$ is a local slice for $G_1$-conjugation near $f$, 
while $f\cdot\exp(\ul\frt)$ is one for $T$-conjugation in $Q_T$. Therefore, 
a local, $G^f$-equivariant model for $\omega$ is the Dirac induction 
map of \S\ref{dirind},
\begin{equation} \label{localweyl}
G^f \times_{N^f}\ul\frt \to \frg^f,
\end{equation}
and our claim reduces to Proposition \ref{diracbwb}.
\end{proof}

\begin{proof}[Proof of (\ref{kcompute})]
We use the construction \eqref{orbits} of $K$-classes from single 
conjugacy classes. Let $f\in Q_T$ and observe, from $f\ul T \cong f\cdot
\ul\frt/\ul\Pi$ and Lemma \ref{gentwistedconj}, that the Weyl group of 
$G^f$ is identified with the stabiliser in $\waff$ of the associated weight. 
Singular weights are then fixed by the Weyl reflection in some $\mathfrak
{sl}_2$ centralising $f$, and their $K$-classes are killed by the local 
induction \eqref{localweyl}. Near a regular $f$, on the other hand, the 
local model for $\omega$ is an isomorphism, so regular weights contribute 
non-zero generators in $K^\tau_G(G)$.
\end{proof}

\part{Loop groups and admissible representations}
In this chapter, we summarise some basic facts about loop groups, 
twisted loop groups and their Lie algebras, as well as the classification 
of admissible representations in terms of the action on affine regular 
weights of the extended affine Weyl group. The key result, Thm.~\ref
{repclassif}, is certainly known, but does not seem to appear in the 
literature in this form (but see \cite{val} for simple groups). 
This combines the theorem of the lowest weight with Mackey's 
irreducibility criterion.

We shall need to distinguish between representations of the 
polynomial loop algebras and their Hilbert space completions, 
and we convene to mark \textit{uncompleted} spaces by a prime. 

\section{Refresher on affine algebras}
\label{affnot}
\subsection{Affine algebras.} \label{basicext} We use the notations of 
\S\ref{finitedirac}; in particular, $\frg$ is now simple. The \textit
{Fourier polynomial loop algebra} $\plg_\bC$ has the Fourier basis
$\xi_a(m) = z^m\xi_a$. Its \textit{basic central extension} 
$\tildeplg := \mi\bR K \oplus \plg$, with central generator $K$, is
defined by the 2-cocycle sending $\xi\wedge\eta\in \Lambda^2\plg$ to 
$K\cdot\mathrm{Res}_{z=0}\langle d\xi |\eta \rangle$. The \textit{affine 
Lie algebra} $\hatplg = \tildeplg\oplus\mi\bR E$ arises by adjoining a 
new element $\mi E$, where the \textit{energy} $E$ satisfies $\left[ 
E,K\right] = 0$ and $\left[E,\xi(n)\right] = n \xi(n)$, for any $\xi 
\in \frg$. Unlike $\tildeplg$, $\hatplg$ carries an ad-invariant bilinear 
form, extending the basic one on $\frg$: 
\begin{equation} \label{bil}
\langle k_1K+\xi_1+e_1E\,|\,k_2K+\xi_2+e_2E \rangle
\mapsto \frac{1}{2\pi}\int_0^{2\pi} {\langle \xi_1(t)|
\xi_2(t)\rangle dt} + k_1e_2 + k_2e_1.
\end{equation}
\subsection{Lowest-weight modules.}\label{affdefs} 
A projective representation of $\plg$ has \textit{level $k$} if it 
extends to a strict representation of $\widetilde{L}'\frg$ in which $K$ 
acts as the scalar $k$. This means that we can choose the action $R_a(m)$ 
of $\xi_a(m)$ so that
\[
[R_a(m),R_b(n)]=f_{ab}^cR_c(m+n)+km\delta _{ab}\delta _{m,-n}.
\]
Call $\frh := \mi\bR K \oplus\frt\oplus\mi\bR E$ a \textit{Cartan 
sub-algebra} of $\hatplg$, and let $\frN := \bigoplus_{n>0} z^n\frg_
\bC\oplus\frn \subset \plg_\bC$. A \textit{lowest weight vector} in an 
$\hatplg$-module $\philb$ is an $\frh$-eigenvector killed by $\ol{\frN}$. 
Call $\philb$ a \textit{lowest weight module}, with lowest weight $\left
({k,-\lambda,m} \right)$, if it is generated by a lowest weight vector 
$\bfv$ of that $(K,\frt,E)$-weight. The factorisation $U(\hatplg) 
= U(\frN)\otimes U(\frh)\otimes U(\ol{\frN})$ shows that $\philb$ is 
generated by $\frN$ from $\bfv$. Defining the \textit{positive alcove} 
$\fra\subset\frt$ as the subset of dominant elements $\xi$ satisfying 
$\theta(\xi) \le 1$, we have the following:

\begin{proposition}\label{ineq}
In a $\left(k,-\lambda,m \right)$-lowest-weight module, the weight 
$\left(k,\omega,n \right)$ of any other $\frh$-eigenvector satisfies 
$n-m \in \bZ$ and $\left(\omega +\lambda\right)(\xi) + n > m$, for $\xi$ 
inside $\fra$. 
\end{proposition}

\begin{proof} All weights of $\frN$ verify these conditions, with 
$\lambda, m = 0$.  
\end{proof}

\subsection{Integrable modules.}\label{integrable} A lowest-weight 
module is \textit{integrable} if the action exponentiates to the 
associated simply connected loop group.\footnote{A precise definition 
is a bit delicate, but there are some simple equivalent Lie algebra 
conditions \cite[III]{kac}; for instance, it suffices that the action 
should exponentiate on all root $\mathfrak{sl}_2 $-subgroups \cite[VII]
{psloop}.} Integrable representations are unitarisable, completely 
reducible, and the irreducible ones are parametrised by their lowest 
weights $\left(k,-\lambda,m \right)$, in which $k$ must be a non-negative 
integer and $\lambda$ a dominant $T$-weight satisfying $\lambda \cdot 
\theta \le k$, in the basic inner product. These weights correspond 
to points of the scaled alcove $k\cdot\fra$.

\subsection{Spinors.} The \textit{complex Clifford algebra} $\cliff
(\plg^*)$ is generated by the odd elements $\left\{\psi^a(m)\right\}$ 
dual to $\left\{\xi_a(m) \right\}$, satisfying
\[
\psi^a(m)\psi^b(n) + \psi^b(n)\psi^a(m) = 2\delta_{ab}\delta_{m-n}.
\]
Choose an irreducible, $\bZ/2$-graded, positive energy module $\pspin$ 
of $\cliff(\plg^* )$. As a vector space, this can be identified with 
the graded tensor product $\spin(0) \otimes \Lambda^\bullet\left(z\frg_
\bC[z]\right)$, for an irreducible, graded spin module $\spin(0)$ of 
$\cliff(\frg^*)$. $\pspin$ carries a hermitian metric, in which 
$\psi^a(n)^* = \psi^a(-n)$; so $\psi (\mu)$ is self-adjoint for $\mu 
\in \plg^*$. There are obvious actions of $\frg$ and $E$ on $\pspin$, 
intertwining with $\cliff(\plg^*)$. The lowest $E$-eigenvalue is 0, 
achieved on $\spin(0) \otimes 1$. Setting  
\begin{equation} \label{quad}
K \mapsto h^\vee, \qquad \xi_a(m) \mapsto \sigma_a(m) := 
-\frac{1}{4}\sum\nolimits_{p+q=m}{f_{bc}^a\,\psi ^b(p)\psi ^c(q)}  
\end{equation}
extends them to an action of $\hatplg$, with intertwining relation 
$\left[\sigma_a(m),\psi^b(n)\right] = f_{ca}^b\psi ^c(m+n)$. One 
derives \eqref{quad} by considering the adjoint representation $\plg$ 
to the orthogonal Lie algebra $\mathfrak{so}_\mathrm{res}(L\frg)$, 
``restricted" as in \cite{psloop} with respect to the splitting 
$\plg_\bC = z\frg_\bC[z] \oplus \frg_\bC[z^{-1}]$. Formula \eqref{quad} 
is then the quadratic expression of the spin representation of 
$\mathfrak{so}_\mathrm{res}$ in terms of Clifford generators 
\cite{koststern}. 

The following key result follows from the Kac character formula.
It is part of affine algebra lore; but see \cite{fht2} for a proof.

\begin{proposition} \label{spinrep}
As a representation of $\hatplg$, $\pspin$ is a sum of copies of the
integrable irreducible representation of level $h^\vee$ and lowest
weight $(-\rho)$. The lowest weight space, which is isomorphic to
the multiplicity space, is also the $\frg$-lowest-weight space in 
$\spin(0)$, and is a graded, irreducible $\cliff(\frt^*)$-module.
\end{proposition} 

\section{Twisted affine algebras} \label{twistedaff}
The loop algebras $L\frg$ have \textit{twisted versions}, arising from
the automorphisms of non-trivial principal $G$-bundles over the circle.
They are closely related to the outer automorphisms of $\frg$ and to the 
\textit{twisted simple affine algebras} in Tables (Aff 2,3) of \cite{kac}: 
each of those is a central extension of a twisted loop algebra,\footnote
{More precisely, the basic central extension of the sub-algebra of 
twisted, $\frg$-valued Fourier polynomials.} plus an outer derivation $E$. 

\subsection{} \label{81}
The algebra $L_\vep\frg$ of loops in $\frg$ twisted by an automorphism 
$\vep$ depends, up to isomorphism, only on the conjugacy class of $\vep$ 
in the outer automorphism group of $\frg$. Thanks to Proposition \ref
{diagaut}, we may as well assume that $\vep$ is a diagram automorphism. 
When $\frg$ is simple, this will have order $1,2$ or $3$; in general, 
we insist that the order $r$ should be finite. This leads to an attractive 
algebraic model for $L_\vep \frg$ as the invariant part of a copy of 
$L\frg$, based on the $r$-fold cover $\rS$ of the unit circle $S^1$, 
under the Galois automorphism which rotates the cover by ${2\pi}/r$ and 
applies $\vep$ point-wise. To find the geometric meaning of this 
construction, let $G_1$ be the simply connected group with Lie algebra 
$\frg$ and $G = \bZ/r \ltimes_\vep G_1$. The quotient of the trivial 
bundle $\rS\times G$ under the action of $\bZ/r$ which rotates the 
circle and left-translates the fibres $G$ is a principal $G$-bundle $P$ 
over $S^1$, and its Lie algebra of gauge transformations is $L_\vep\frg$.

\subsection{Standard form.}\label{structwist} 
Let $\frg$ be simple, for the rest of this section. The (smooth) twisted 
affine algebra $\widehat{L}_\vep\frg$ is the invariant part of $\widehat{L}
\frg$, in our construction above. Its structure is described in \cite[VI, 
VIII]{kac}. Inherited from the ambient $\widehat{L}\frg$ is a linear 
decomposition $\widehat{L}_\vep\frg = \mi\bR K \oplus L\frg\oplus\mi\bR E$. 
We now \textit{rescale} $K$ and $E$ to $r\times$, resp.\ $1/r\times$ the 
originals. Then, $E$ is the natural generator for the rotation of the unit 
(downstairs) circle, while the bilinear form \eqref{bil} is still 
ad-invariant. Using the \textit{standard connection} $\nabla_0$ on $P$, 
descended from the trivial connection on $\rS$, the $2$-cocycle of the 
central extension $\widetilde{L}_\vep \frg := \mi\bR K \oplus L_\vep \frg$ 
is again expressed as an integral over the unit circle, and the Lie bracket 
takes the following form:
\begin{equation}\label{twistcoc}
\left[\xi,\eta\right](t) = \left[\xi(t),\eta(t)\right]
+ \frac{K}{2\pi\mi} \oint {{\langle\nabla _0\xi |\eta \rangle}}.
\end{equation}

\subsection{Lowest-weight modules.} \label{posroot}
The r\^oles of $\frt$, $\frh$ and $\frN$ are taken over by their Galois 
invariants within the ambient $L'\frg$; we denote them by underlines. 
The structure of weights and roots parallel the untwisted case;\footnote
{Corollary \ref{constext} below imposes a small distinction for the weight 
lattice for twisted $\mathrm{SU}(2\ell+1)$; see (\ref{weights}).} details 
are summarised in Appendix \ref{afftwist}. 

In a \textit{level $k$} representation of $\hatplge$, $K$ acts as the 
scalar $k$. A \textit{lowest weight vector} is an $\frh$-eigenvector killed 
by $\ol{\ul\frN}$, and a \textit{lowest-weight module} is one generated by 
a lowest weight vector. We call such a module \textit{integrable} if the 
action of all the root $\mathfrak{sl}_2$ sub-algebras of $\hatplge$ is so; 
in that case, the module is unitarisable, and the Lie algebra action 
exponentiates to one of $\widehat{L}_\vep G$ on the Hilbert space completion. 
Integrable representations are semi-simple, and the irreducible ones are 
parametrised by their level $k$ and their lowest weight $(k,-\lambda)$, 
in which $\lambda$ is dominant and satisfies $\ul\theta \cdot\lambda\le 
k/r$.

The underlined Doppelg\"angers for $\rho$, $\theta$ and $\fra$ require 
a comment: $\ul\rho$ has the obvious meaning, the half-sum of positive 
roots for $\ul{\frg} := \frg^\vep$, but $\ul\theta$, which cuts out $\ul\fra$ 
from the dominant chamber of $\ul\frt$ by the relation $\ul\theta(\xi)\le 1/r$, 
is \textit{not} the highest root of $\ul\frg$, but rather the highest 
weight of $\frg/\ul\frg$. Therewith, the analogue of Proposition \ref{ineq} 
holds true.

A geometric sense in which $\ul\fra$ plays tho r\^ole of $\fra$ is the 
following. Let $\cA$ denote the space of smooth connections on the bundle 
$P$; the quotients $\cA\left/L_\vep G_1\right.$ (by gauge transformations) 
and $\vep G_1\left/{G_1}\right.$ (by conjugation) are isomorphic by the 
holonomy map. The classification \eqref{twistedconj} of twisted conjugacy 
classes gives the following.
 
\begin{proposition}\label{gauge} Every smooth connection on $P$ is a 
smooth gauge transform of $\nabla_0 + \xi dt$, for a \textit{unique} $\xi 
\in\ul\fra$. That is, $\ul\fra$ is a global slice for $L_\vep G$: 
$\ul\fra \cong \cA\left/L_\vep G\right.$. \hfill\qed
\end{proposition}

\subsection{The Clifford algebra.} \label{twistcliff}
A basis for $\hatplge_\bC$ suited to calculations arises from a complex 
orthonormal $\vep$-eigen-basis $\{\xi_a\}$ of $\frg_\bC$, so chosen that 
the indexing set carries an involution $a\leftrightarrow\bar a$ with 
$\xi_{\bar a} = -\xi_a^*$. If $\vep(a) \in \bZ/r$ corresponds to the $\vep
$-eigenvalue of $\xi_a$, then $\{\xi_a(m)\}$ forms a basis of $\frg$, as 
$m+\vep(a)/r$ ranges over $\bZ$. Raising and lowering indexes involves 
a bar; for instance, the relations in the complex Clifford algebra of 
$\plge^*$ are $\left[\psi^a(m),\psi^b(n)\right] = 2\delta_{a\ol{b}} 
\cdot\delta_{m-n}$. 

A positive energy, graded spin module $\pspin$ can be identified, as a 
vector space, with $\spin(0) \otimes \Lambda^\bullet\left(\ul\frN \right)$, 
for a graded spin module $\spin(0)$ of $\cliff(\ul \frg^*)$. As in (\ref
{quad}), the obvious actions of $\ul\frg$ and $E$ extend to a lowest-weight 
representation of $\hatplge$, with a bar in the raised index $a$, but, 
remarkably, with the same $h^\vee$. As in Prop.~\ref{spinrep}, the 
representation can be identified using the Kac character formula; its 
lowest-weight space is the lowest $\ul\frg$-weight space in $\spin(0)$, 
has pure weight $(-\ul\rho)$ and is a graded irreducible $\cliff(\ul
{\frt}^*)$-module. 

\subsection{The loop group.} The extensions in \S\ref{basicext} and 
\eqref{twistcoc} are so normalised as to generate all central extensions 
of $LG$ by the circle group $\bT$. Call $L_\vep G_1$ the twisted loop 
group of $G_1$.

\begin{proposition} \label{basictwext}
$\widetilde{L}_\vep \frg$ is the Lie algebra of a \emph{basic central 
extension} $\widetilde{L}_\vep G$ of $L_\vep G$, with central circle 
parametrised by $\left\{\left.z^K\right|\lvert z\rvert =1\right\}$, 
and whose Chern class generates $H^2\left(L_\vep G_1;\bZ\right)=\bZ$.
\end{proposition} 

\begin{proof} The untwisted case is handled in \cite{psloop}, so we 
focus on $r>1$. Being the space of sections of a $G_1$-bundle over $S^1$, 
$L_\vep G_1$ is connected and simply connected. Further, $\pi_2L_\vep G_1 
= H^1\left(S^1;\pi_3G_1\right) = \bZ$, and Hurewicz gives us $H^2(L_\vep 
G_1;\bZ) = \bZ$. Since $\pi_2LG_1 = H^1(\rS;\pi_3G_1)$, the restriction 
$H^2(LG_1) \to H^2(L_\vep G_1)$ has index $r$. Our extension of $L_\vep 
G_1$ will be the $r$-th root of the restriction of $\widetilde{L}G_1$, 
the basic extension of the ambient, untwisted loop group. Having fixed 
the cocycle \eqref{twistcoc}, the obstructions to existence and uniqueness 
of this root are topological, living in $H^2$ and $H^1$ of $L_\vep G_1$ 
with $\bZ/r$-coefficients, respectively; and they vanish as seen. Finally, 
we have a semi-direct decomposition $L_\vep G\cong\bZ/r\ltimes L_\vep G_1$, 
and the $\vep$-action on $L_\vep G_1$ preserves the cocycle \eqref{twistcoc}, 
so it lifts to an automorphism action on the central extension (again, 
by vanishing of the topological obstructions). We let $\widetilde{L}_\vep G
=\bZ/r\ltimes_\vep L_\vep G_1$. 
\end{proof}

\begin{corollary} \label{constext}
The basic extension $\widetilde{L}_\vep G$ restricts trivially to the 
constant subgroup $G_1^\vep$, except when $G_1=\mathrm{SU}(2\ell +1)$ 
and $r=2$, in which case $G_1^\vep = \mathrm{SO}(2\ell +1)$, and we 
obtain the $Spin^c$-extension. 
\end{corollary}  

\begin{proof} The flag variety $L_\vep G_1/G_1^\vep$ is simply connected, 
with no $H^3$. (This follows, for instance, from its Bruhat stratification
by even-dimensional cells.) The Leray sequence for the fibre bundle 
$L_\vep G_1 \twoheadrightarrow L_\vep G_1/G_1^\vep$ shows that $H^2
(L_\vep G_1)$ surjects onto $H^2(G_1^\vep)$. However, $G_1^\vep$ is 
simply connected, save in the cases listed, whence the result.
\end{proof}

\section{Representations of $L_fG$} \label{repclassifsect}
We now classify the admissible representations of the loop groups at levels 
$\tau-\ul\sigma$ for which $\tau$ is regular, in terms of the 
affine Weyl action on regular weights.

\subsection{Notational refresher.}\label{refresh} Let $f$ be an element 
of the quasi-torus $Q_T$ and call $L_fG$ is the $f$-twisted smooth loop group 
of $G$ (\S\ref{twistedgroups}), $\tau$ a regular central extension and $\ul\sigma$ the extension defined by the spin module $\spin$ of $L_f\frg^*$ 
(\S\ref{morita}). Gradings are incorporated into our twistings. The extended 
affine Weyl group $\waff = \pi_0L_fN$ acts on $\ul{\wtlat}^\tau$ by conjugating 
the central extension of $\ul T$, and a tautological twisting $\tau'$ is 
defined for this action, wherein each $\tau$-affine weight defines a 
$\bT$-central extension of its stabiliser in $\waff$ (\S\ref{philoop}). 
We now restate Theorem~4 without Clifford algebras; it is the 
lowest-weight classification of representations, enhanced to track the 
action of the components of $L_fG$.
 
\begin{theorem} \label{repclassif}
\begin{trivlist}
\item\textup{(i)} The category of admissible representations of $L_fG$ 
of level $\tau -\ul\sigma$ is equivalent to that of finite-dimensional, 
$\waff$-equivariant, $\tau'-\sigma(\ul\frt)$-twisted vector bundles over 
$\ul{\wtlat}^\tau_{reg}$. 
\item\textup{(ii)} The $K$-groups of graded admissible representations 
are naturally isomorphic to the twisted equivariant $K$-theories 
$K_{\waff}^{\tau'-\sigma(\ul\frt)+*}(\ul{\wtlat}^\tau_{reg})$. 
\end{trivlist}
\end{theorem}

\noindent The reader may wish to consult the simple Example \ref{example}, 
where $G=N$. In general, the equivalence in (i) arises as follows. A 
regular weight $\mu$ defines a polarisation of $L_f\frg$, which selects, 
for each admissible representation $\hilb$, a lowest-weight space in 
$\hilb\otimes\spin$ with respect to $L_fG\ltimes\cliff(L_f\frg^*)$. The 
$(-\mu)$-eigen-component under $\ul T$ of this lowest weight-space is 
a $\cliff(\ul\frt)$-module, and factoring out the spinors on $\frt$ gives 
the fibre of our vector bundle at $\mu\in \ul{\wtlat}^\tau$. 

Dirac induction provides the inverse equivalence. Each $\mu\in\ul{\wtlat}
^\tau_{reg}$ determines a regular co-adjoint orbit $\frO_\mu\subset L_f
(\frg^*)^\tau$, over which a twisted representation of the $\waff$-stabiliser 
defines a $(\tau-\sigma(\ul\frt))$-twisted, $L_fG$-equivariant vector 
bundle. The Dirac index of this bundle along $\frO_\mu$, coupled to the \textit{highest-weight} spinors, is the desired representation of $L_fG$. 
Its level $(\tau-\ul\sigma)$ arises from the shift by the level $\sigma
(\ul\frt)-\ul\sigma$ of the highest-weight spinors on $L_f\frg/\ul\frt$. 

Dirac induction in infinite dimensions is only a heuristic notion, but 
can be realised in this case by the \textit{Borel-Weil construction}, as 
a space of holomorphic sections \cite{psloop}. We will review that in 
\S\ref{fusion}, where it is needed, but we will make no use of it this 
section. 

\vskip2ex
Proving \eqref{repclassif} requires some preparation. Split $\frg$ 
into its centre $\frz$ and derived sub-algebra $\frg'$.  

\begin{proposition}\label{splitlooplie} $L_f\frg'$ splits canonically 
into a sum of simple, possibly twisted loop algebras. Central extensions 
of $L_f\frg$ are sums of extensions of $L_f\frz$ and of the simple 
summands.
\end{proposition} 
\begin{proof} In the decomposition of $\frg'$ into simple ideals, $f
$-conjugation permutes isomorphic factors. To a cycle $C$ of length 
$\ell(C)$ in this permutation, we assign one copy of the underlying 
simple summand $\frg(C)$ and the automorphism $\vep(C):= \Ad(f)^
\ell$. This is a diagram automorphism of $\frg(C)$, whose fixed-point 
sub-algebra is isomorphic to that of $\Ad(f)$ on the summand 
$\frg(C)^{\oplus\ell}$ in $\frg'$. Then, $L_f\frg'$ is isomorphic to the 
sum of loop algebras $L_{\vep(C)}\frg(C)$, with the loops parametrised 
by the $\ell(C)$-fold cover of the unit circle. The splitting arises from 
the eigenspace decomposition of $\Ad(f)$ on $\frg(C)^{\oplus\ell}$. 
As the summands are simple ideals, uniqueness is clear. The splitting 
of the extension follows from the absence of one-dimensional characters 
of the simple summands.
\end{proof}

\subsection{More on $\waff$.} \label{morewaff}
The proposition splits $\ul\frt$ into $\ul\frz:=\frz^f$ and the sum of 
the Cartan sub-algebras $\ul\frt(C)$. Call $\tau\cdot\ul\fra\in\ul\frt$ the 
product of $\ul\frz$ and the positive alcoves (\S\ref{affdefs}, 
\S\ref{structwist}) in the $\ul\frt(C)$, scaled by the simple components 
of the level $[\tau]$, and let $\tau\cdot\ul\fra^*$ be its counterpart in 
$\frt^*$ in the basic inner product on $\frg'$. Reflection about the 
walls of $\tau\cdot\ul\fra^*$ generate a normal subgroup $\waff(\frg,f) 
\subset\waff$, under whose action the transforms of the alcove are 
distinct and tessellate $\ul\frt^*$ \eqref{scwaff}. The two groups agree 
when $G$ is simply connected, but in general we have an exact sequence, 
split by the inclusion of $\pi$ in $\waff$ as the stabiliser of $\tau\cdot\fra$,
\begin{equation}\label{wext}
1 \to \waff(\frg,f) \to \waff\to \pi:=\pi_0L_fG\to 1.
\end{equation}

Regular are those weights not lying on any alcove wall (Remark~\ref
{waffdef}). The alcoves correspond to positive root systems on $L_f\frg$ 
which are conjugate to the standard one (\S\ref{posroot}), the simple 
roots being the outward normals to the walls. The positive root spaces 
span a polarisation of $L_f\frg'$; the various polarisations, plus the 
original one on $L_f\frz$, are conjugate under $\Gamma_fN\subset L_fG$, 
so they define the same class of admissible representations. 

\subsection{Mackey decomposition in $K$-theory.}\label{normal}
Let $H$ be a group, $M$ a normal subgroup, $\upsilon$ a central extension 
of $H$. Conjugation leads to an action of $H/M$ on isomorphism classes 
of $\upsilon$-representations of $M$. Let $Y$ be a family of isomorphism 
classes, satisfying the conditions
\begin{enumerate}\itemsep0ex
\item $Y$ is stable under $H/M$;
\item Every point in $Y$ has finite stabiliser in $H/M$;
\item The $M$-automorphisms of any representation in $Y$ are scalars.
\end{enumerate}
\noindent There is a tautological projective vector bundle $\bP R$ 
over $Y$, whose fibre $\bP R_y$ at $y\in Y$ is the projective space on a 
representation of isomorphism type $y$. Its uniqueness up to canonical 
isomorphism, and hence $H$-equivariance, follow from condition (iii). 
The bundle defines a $\bT$-central extension of the action groupoid of 
$H^\upsilon$ on $Y$. This central extension is split over $M^\upsilon$, 
so dividing out by the latter gives a central extension, or twisting, 
$\upsilon'$ for the $H/M$-action on $Y$.

Call an $H$-representation \textit{$Y$-admissible} if its restriction 
to $M$ is a finite-multiplicity sum of terms of type in $Y$, with only 
finitely many $H/M$-orbit types. For instance, this includes all induced 
representations $\mathrm{Ind}_M^K(R_y)$. The same construction as in 
Lemma~\ref{keylemma} establishes the following:
\begin{proposition} \label{repextension}
The category of $Y$-admissible representations of $H$ is equivalent to 
that of $\upsilon'$-twisted, $H/M$-equivariant vector bundles over $Y$, 
supported on finitely many orbits. \qed
\end{proposition} 
\noindent In this equivalence, a $M$-representation $\hilb$ 
is sent to the bundle whose fibre at $y$ is $\mathrm{Hom}^M(R_y,\hilb)$. 
Conversely, to a bundle over $Y$ we associate its space of sections. 
The relation to Construction \ref{constr} can be made explicit by 
choosing a representation $\hilb$ of $H^\upsilon$ containing all elements 
of $Y$. The projective bundle $\bP\mathrm{Hom}^M(R;\hilb)$ over $Y$ gives 
a model for the twisting $\upsilon'$ of the $H/M$-action.

\begin{proof}[Proof of (\ref{repclassif}).] The unitary lowest-weight 
representations of the Lie algebra correspond to the admissible ones of 
the simply connected cover of the identity component $(L_fG)_1$. For 
the simple summands, integrable representations are classified by 
lowest-weights \cite{kac}. Analytic regularity of $\tau$ on the centre 
$L_f\frz \cong \ul\frz \oplus {L_f\frz}\left/\ul\frz\right.$ means that 
the second summand has a unique irreducible lowest-weight representation. 
Unitary irreducibles of $\ul\frz$ are labelled by the points of 
the $\tau$-affine dual space. Descent of representations to $(L_fG)_1$ 
is controlled by an integrality constraint imposed by $\ul T$: 
parametrising the admissible irreducibles of $(L_fG)_1$ by their lowest 
weights $(-\lambda)$, the shifted weights $(\lambda+\ul\rho)$ range over 
$\ul{\wtlat}^\tau_{reg+} := \ul{\wtlat}^\tau_{reg}\cap \tau\cdot\ul\fra^*$. 

As $\waff(\frg,f)$ acts freely on $\ul \wtlat^\tau_{reg}$, and the orbits 
are in bijection with the points in $\ul\wtlat^\tau_{reg+}$, we get an 
identification 
\begin{equation}\label{repcompute}
K^{\tau'-\sigma(\ul\frt)}_{\waff}(\ul{\wtlat}^\tau_{reg}) = 
	K^{\tau'-\sigma(\ul\frt)}_\pi(\ul{\wtlat}^\tau_{reg+}).
\end{equation} 
We apply Proposition \ref{repextension} to $H = L_fG$, $M = (L_fG)_1$, 
$\upsilon = \tau-\ul\sigma$, $Y=\ul{\wtlat}^\tau_{reg+}$. The actions of 
$\pi$ described in \S\ref{refresh} and \S\ref{repextension} do match, because 
the (sign-reversed) lowest weight $(\ul\sigma,\ul\rho)$ of $\spin$ is 
$\pi$-invariant. To conclude the proof, it remains to identify the 
$\pi$-twistings $\upsilon'$ and $\tau'-\sigma(\ul\frt)$. 

The subgroup of $\naff$ lying over $\pi$ preserves the lowest-weight 
space in any $\upsilon$-representation $\hilb$ of $L_fG$, and so the 
projective action of $\pi$ on the resulting lowest-weight bundle over 
$Y$ represents $\upsilon'$. Similarly, a model for $\tau'$ arises from the 
action of $\pi$ on the lowest-weight space in $\hilb\otimes\spin$, 
distributed over the (sign-reversed) eigenvalues in $\ul{\wtlat}^\tau_
{reg+}$. The second bundles differs from the first by a factor 
of $\spin(\ul\frt)$, and this represents the twisting $\sigma(\ul\frt)$.  
\end{proof}

\subsection{Example: $G=N$.}\label{example}  
Let $V := L_f\frt \ominus \ul{\frt}$ and $L_fN \cong \nafff \ltimes 
\exp(V)$, as in \S\ref{anregsect}. Regularity of $\tau$ confines us to 
sums of Heisenberg extensions of $V$ and topologically regular extensions 
$\Gamma_fN^\tau$. The lowest-weight module $\fock$ of $\exp(V)$ carries 
a (projective) intertwining action of $\pi_0\nafff$. An admissible 
representation $\hilb$ of $L_fN$ factors then as $\fock \otimes\mathrm
{Hom}^V(\fock;\hilb)$, where the second factor is (the $\ell^2$ completion 
of) a \textit{weight module} of $\naff$, which means that it is 
$\ul T$-semi-simple, of finite type. Our classification now becomes the 
following, more precise 

\begin{proposition}
Global sections give an equivalence from the category of $\waff
$-equivariant, $\tau'$-twisted vector bundles on $\ul{\wtlat}^\tau$
with that of weight $\tau$-modules of $\nafff$.
\end{proposition}

\noindent It is understood here that $\ul{T}^\tau$ acts with weight 
$\lambda$ on the fibre at $\lambda \in \ul{\wtlat}^\tau$. The proposition 
follows directly from Prop.~\ref{repextension}. Weight modules split 
into irreducibles, which are induced from stabilisers of single weights. 

\part{From representations to $K$-theory}\label{diracfam}
To an admissible representation $\hilb$ of $LG$ at fixed level 
$\tau-\sigma$, we assign a family of Fredholm operators parametrised 
by an affine copy of $L\frg^*$, equivariant for the \textit{affine 
action} of the loop group $LG$ at the shifted level $\tau$. The 
underlying space of the family is $\hilb \otimes\spin$, and the 
operator family is the analogue of the one in \S\ref{finitedirac}, 
but is based on the \textit{Dirac-Ramond operator}. We recall 
this operator in \S\ref{affdiracsect}, and reproduce the calculation 
\cite{land}, \cite{taub} of its Laplacian, which we extend to twisted 
algebras. Our family defines an $LG$-equivariant twisted $K$-theory 
class over $L\frg^*$, which we identify, when $\hilb$ is irreducible, 
with the Thom push-forward of the natural line bundle on a single, 
integral co-adjoint orbit. The passage from representation to orbit 
and line bundle is an inverse of Kirillov's quantisation of co-adjoint 
orbits. The affine copy of $L\frg^*$ carrying our family can be 
identified with the space of $\frg$-connections over the circle with 
the gauge action, leading to an interpretation of our family as a 
cocycle for $K^\tau_G(G)$.

\section{The affine Dirac operator and its square} 
\label{affdiracsect}
Let $\frg$ be simple and let $\philb$ be a lowest weight module for 
$\hatplg$, with lowest weight $(k,-\lambda,0)$. Consider the following 
formally skew-adjoint operator on $\philb \otimes \pspin$: 
\begin{equation}
\label{affdirac}
\dirac = \dirac_0 := R_a(m) \otimes \psi^a(-m)
+ \frac{1}{3} \cdot \sigma_a(m)\psi^a(-m).
\end{equation}
This is known to physicists as the Dirac-Ramond operator \cite{mick}; 
in the mathematical literature, it may have been first considered by 
Taubes \cite{taub}, and, more recently, studied in detail by Landweber 
\cite{land}, based on Kostant's compact group analogue. Denote by 
$T_a(m)$ the total action $R_a(m)+\sigma _a(m)$ of $\xi _a(m)$ on 
$\philb \otimes \pspin$, and let $k^\vee :=k+h^\vee $.

\begin{proposition} \label{affdiracrels} 
$\left[\dirac,\psi^b(n) \right] = 2 T_b(n)$, \ 
$\left[\dirac,T_b(n)\right] = - n k^\vee\cdot \psi^b(n)$. 
\end{proposition}

We postpone the proof for a moment and explore the consequences. 
Clearly, the commutation action of the Dirac Laplacian $\dirac^2$ 
on the $T_\bullet$ and the $\psi$ agrees with that of $-2k^\vee E$. 
Normalise the total energy operator $E$ on $\philb \otimes \pspin$ 
to make it vanish on its lowest eigenspace $\hilb(0)\otimes \spin(0)$. 
This last space is $\dirac$-invariant, and the only terms in \eqref
{affdirac} to survive on it are those with $m=0$. These sum to the 
Dirac operator for $\frg$, acting on its representation $\hilb(0)$. 
The latter squares to $-(\lambda +\rho)^2$. Since $\philb \otimes
\pspin$ is generated by the actions of the $T_\bullet$ and the 
$\psi$ on $\hilb(0)\otimes \spin(0)$, the following formula for the 
Dirac Laplacian results:

\begin{equation}
\label{affdiracsquare}
\dirac^2=-2k^\vee E-\left(\lambda +\rho\right)^2.
\end{equation}
In particular, $\dirac$ is invertible, with discrete, finite multiplicity 
spectrum.

\begin{remark} \label{comarg}
Because the $\sigma$ are expressible in terms of the $\psi$, the 
Dirac operator \eqref{affdirac} is expressible in terms of the operators
$T_\bullet$ and $\psi$ alone. Define the \textit{level $k^\vee$ 
universal enveloping algebra} of $\plg$, $U_{k^\vee}(\plg):= 
{U(\tildeplg)}/(K-k^\vee)$. Then, $\dirac$ is an odd element in a 
certain completion\footnote{The most natural completion is that 
containing infinite sums of normal-ordered monomials, of bounded degree 
and energy; this acts on all lowest weight modules of $\plg\ltimes \psi(\plg^*)$} of the ``semi-direct tensor product'' of $\mathrm
{Cliff}(\plg^*)$ by $U_{k^\vee}(\plg)$, acting via $\ad$. The first 
equation in \eqref{affdiracrels} determines $\dirac$ uniquely, because 
no odd elements of the completed algebra commute with all the $\psi$. 
However, a definite lifting $T_\bullet$ of $L\frg$ into $U_{k^\vee}
(L\frg)$ has been chosen here. This shows up more clearly in the next 
section, where we consider the family of $\dirac$'s parametrised by 
all possible linear splittings of the central extension 
$\widetilde{L}\frg$ (cf.~also \S\ref{coordfree}). 
\end{remark}

\begin{proof}[Proof of (\ref{affdiracrels}).] The first identity
follows by adding the two lines below, in which summation over $m\in \bZ$
is implied, in addition to the Einstein convention: 
\begin{eqnarray*}
\left[ {R_a(m)\otimes \psi^a(-m),\psi^b(-n)} \right] 
		&=& 2 R_b(-n),\\
\left[ {\sigma_a(m)\psi^a(-m),\psi ^b(-n)} \right] &=& 
	2\sigma_b(-n)+ f_{ac}^b \psi^c(m-n)\psi^a(-m) 
	= 6\sigma_b(-n). 
\end{eqnarray*}
The second identity in \eqref{affdiracrels} follows from the first. 
Indeed: 
\begin{eqnarray*}
\left[\left[ \dirac,T_b(n) \right],\psi ^c(p) \right] 
&=& \left[ \dirac,\left[ T_b(n),\psi^c(p) \right]\right]
-\left[ T_b(n),\left[ \dirac,\psi^c(p) \right] \right]\\
&=& f_{db}^c \left[ \dirac,\psi^d(p+n)\right]
	- 2 \left[T_b(n),T_c(p)\right]\\
&=& 2 f_{db}^c T_d(p+n) - 2 f_{bc}^d T_d(p+n) - 2 nk^\vee\cdot 
\delta_{bc} \delta_{n+p}\\
&=& -2 nk^\vee \cdot \delta _{bc}\delta _{n+p}\\
&=& -nk^\vee\left[\psi ^b(n),\psi ^c(p) \right],
\end{eqnarray*} 
whence we conclude that the odd operator $\left[\dirac,T_b(n)\right] 
+ nk^\vee \psi^b(n)$ commutes with all the $\psi$; hence it 
is zero, as explained in Remark \eqref{comarg}. 
\end{proof}

\subsection{The twisted case.} \label{twistdirac} 
With the same notation and the same definition \eqref{affdirac} 
of $\dirac$, we have
\begin{equation}\label{twistdiracrels}
\left[\dirac,\psi^b(n) \right] = 2 T_{\bar b}(n),\quad
\left[\dirac,T_b(n)\right] = - nk^\vee\cdot \psi^{\bar b}(n); 
\end{equation}
and we obtain, as before, the formula for the Dirac Laplacian: 
\begin{equation} \label{twistdiracsquare}
\dirac_0^2 = -2k^\vee E - (\lambda + \ul{\rho})^2.
\end{equation}

\vspace{-2ex}
\subsection{The affine Dirac operator.} \label{coadjoint}
The relation between the finite and affine Dirac Laplacians \eqref
{diracsquare} and \eqref{affdiracsquare} becomes more transparent 
if we use spinors on the full Kac-Moody algebra. Let $\widehat{L}\frg^* 
= \mi\bR\kstar\oplus L\frg^* \oplus\mi\bR\delta$, where $\kstar$ 
(denoted $\Lambda$ in \cite{kac}) is dual to $K$ and $\delta$ to $E$. 
Identifying it with $\widehat{L}\frg$ by the bilinear form \eqref{bil}, 
the co-adjoint action of $\widehat{\xi} = (k,\xi,e)$ on $\widehat{L}\frg^*$ 
becomes
\begin{eqnarray} \label{coad}
\kstar &\mapsto& \mi\:d\xi/dt = -\left[E,\xi\right],\quad
\delta \mapsto 0, \nonumber \\
\mu \in L\frg^* &\mapsto& \ad^\vee_{\xi(t)}\mu(t) + 
e\cdot\mu'(t) + \mi\delta\cdot\oint\mu\xi'dt.
\end{eqnarray}
The Spin module for $\cliff(\widehat{L}\frg^*)$ is $\widehat{\spin} 
= \spin \oplus \psi^{\kstar}\cdot\spin$. The corresponding Dirac operator, 
\[
\widehat{\dirac} := \dirac + E\psi^\delta + K\psi^{\kstar},
\]
commutes with the (new) total action $T_\bullet$ of $\widehat{L}\frg$ 
and satisfies a simpler formula $\widehat{\dirac}{}^2 = -\left(\lambda 
+ \rho\right)^2$, whose verification we leave to the reader.

\section{The Dirac family on a simple affine algebra}
\label{diracfamsect}
We now assume the representation $\philb$ of $\hatplg$ to be integrable; 
it is then unitarisable, and its Hilbert space completion $\hilb$ carries
an action of the smooth loop group $LG$. Furthermore, $k^\vee > 0$.

\subsection{The level hyperplanes.}\label{gaugeaff} The co-adjoint 
action \eqref{coad} preserves the fixed-level hyperplanes $\mi 
k^\vee\kstar + \widetilde{L}\frg^* \subset\widehat{L}\frg^*$. Ignoring 
$\delta$ leads to the \textit{affine action at level $k^\vee$} on 
$L\frg^*$. The correspondence 
\[
\mi k^\vee\kstar + \mu \quad\leftrightarrow\quad d/dt+\mu/k^\vee
\] 
identifies this action with the gauge action on the space 
$\cA$ of $\frg$-valued connections on the circle.

\begin{proposition} \label{equivar}
The assignment $\mu \mapsto \dirac_\mu:=\dirac + \mi\psi(\mu)$, from
$L\frg^*$ to $\mathrm{End}(\philb\otimes\pspin)$, intertwines the 
affine action at level $k^\vee$ with the commutator action. 
\end{proposition}

\begin{proof}
$\left[T(\xi),\dirac_\mu \right] =  k^\vee \psi\left([E,\xi]
\right) + \mi\left[\sigma(\xi),\psi(\mu) \right] = \mi\psi\left(
- k^\vee d\xi/dt + \ad^\vee_\xi(\mu)\right)$, as desired. 
\end{proof}

\subsection{The Laplacian.} \label{ident} Formulae \eqref{affdiracrels} 
and \eqref{affdiracsquare} give
\begin{eqnarray} \label{famsquare}
\dirac_\mu ^2 &=& \dirac^2 +\mi\left[\dirac,\psi(\mu)\right]
	-\psi(\mu)^2\nonumber \\
&=& -2k^\vee E-(\lambda +\rho)^2 + 2\mi\langle T|\mu\rangle - \mu^2 
\\
&=& -2\left(k^\vee E - \mi\langle T\,|\mu\rangle +\langle\lambda +
\rho\,| \mu\rangle\right) - (\lambda +\rho - \mu)^2.
\nonumber
\end{eqnarray}
When $\mu\in\frt^*$, we can view this formula as a generalisation of 
\eqref{affdiracsquare}, as follows. The first term in \eqref{famsquare} 
is $-2k^\vee E_\mu$, with a \textit{modified energy operator} 
\[
E_\mu = E - \mi\langle T\,|\mu/k^\vee \rangle +\langle\lambda +
\rho\,| \mu/k^\vee\rangle. 
\]
This is associated to the connection 
$d/dt + \mu/k^\vee$ in the same way that $E$ is associated to the trivial 
connection: they intertwine correctly with the action of $L\frg$. 
Furthermore, $E_\mu$ is additively normalised so as to vanish on the 
$-(\lambda +\rho )$-weight space within $\hilb(0) \otimes \spin(0)$. 
As we are about to see, when $\mu/k^\vee \in\fra^*$, that weight space 
is the lowest eigenspace for the Dirac Laplacian on $\hilb\otimes\spin$. 

\subsection{The Dirac kernels.} To study a general $\dirac_\mu$, we 
conjugate by a suitable loop group element to bring $\mu$ into 
$k^\vee\fra^*$. As $\dirac_\mu$ now commutes with $\frt$ and $E$, 
we can evaluate (\ref{famsquare}) on a weight space of type $(\omega,n)$, where 
$T(\mu)  = \mi\left\langle \omega | \mu\right\rangle $, and obtain 
\begin{equation} \label{onweightspace}
\dirac_\mu ^2=-2\left(k^\vee n + \langle\omega +\lambda +\rho |\mu\rangle
\right) - (\lambda +\rho - \mu)^2
\end{equation}
Now, a weight of $\hilb \otimes \spin$ splits as $\left(\omega,n\right) 
= \left(\omega_1,n_2\right)+\left(\omega_2,n_2\right)$, into weights of 
$\hilb$ and $\spin$. Proposition (\ref{ineq}) asserts that $(\omega_i+
\lambda)\cdot\mu + k^\vee n_i \ge 0$,with equality only if $\mu/k^\vee$ 
is on the boundary of $\fra^*$, or else if $\omega = -(\lambda +\rho )$ and 
$n=0$. But then, (\ref{onweightspace}) can only vanish if, additionally, 
$\mu=\lambda+\rho$. Since that lies in the interior of $k^\vee\fra^*$, 
we obtain the following.

\begin{theorem} \label{kernel}
The kernel of $\dirac_\mu $ is nil, unless $\mu$ is in the affine 
co-adjoint orbit of $(\lambda+\rho)$ at level $k^\vee$. If so, $\ker
\dirac_\mu $ is the image, under the same transformation, of the 
$-(\lambda+\rho)$-weight space in $\hilb(0)\otimes \spin(0)$.\qed
\end{theorem} 

\noindent The last space is the product of the lowest-weight space 
$\bC\bfv$ of $\hilb(0)$ with that of $\spin(0)$; the latter is a graded, 
irreducible $\cliff(\frt)$-module. As in finite dimensions, the more 
canonical statement is that the kernels of the $\dirac_\mu$ on the 
``critical'' co-adjoint orbit $\frO$ of $\lambda +\rho$ in $\mi 
k^\vee\kstar + L\frg^*$ assemble to a vector bundle isomorphic to 
$\spin(\cN)(-\lambda-\rho)$, the normal spinor bundle twisted by 
the natural line bundle on $\frO$. This vector bundle has a natural 
continuation to a neighbourhood of $\frO$ as the lowest eigen-bundle 
of $\dirac_\mu$. We can describe the action of $\dirac_\mu$ there, 
when $\mu$ moves a bit off $\frO$. 

\begin{theorem} \label{thom}
Let $\mu \in \frO$, $\nu \in \cN_\mu $ a normal vector to $\frO$ at 
$\mu$ in $\cA$. The Dirac operator $\dirac_{\mu +\nu}$ preserves 
$\ker (\dirac_\mu)$, on which it acts as Clifford multiplication 
$\mi\psi(\nu)$. \qed
\end{theorem}  

\subsection{Twisted $K$-theory class. }
Proposition \ref{equivar} shows that our constructions are preserved 
by the action of $LG$, so the Fredholm bundle $\left(\hilb \otimes 
\spin, \dirac_\mu \right)$ over $L\frg^*$ defines a twisted, $LG
$-equivariant $K$-theory class supported on $\frO$.  Formula \eqref
{onweightspace} bounds the complementary spectrum of $\dirac_\mu$ 
away from zero, so the embedding of the lowest eigenbundle induces 
an equivalence of twisted, $LG$-equivariant $K$-theory classes in some 
neighbourhood of $\frO$. Proposition \ref{thom} identifies the $K$-class 
with the Thom push-forward of the line bundle $\cO(-\lambda -\rho)$, 
from $\frO$ to $L\frg^*$. Finally, identifying the level $k^\vee$ 
hyperplane in $\widetilde{L}\frg$ with $\cA$ as in \S\ref{gaugeaff} and 
using the holonomy map from $\cA$ to $G$ interprets our Dirac family as 
a class in $K^\tau_G(G)$, in degree $\dim\frg\pmod{2}$.

\subsection{Twisted affine algebras.} The results extend verbatim to 
twisted affine algebras, if we use the presentation $L_\vep\frg$ of 
\S\ref{twistedaff}. Let $\cA_\vep$ be the space of smooth connections 
on the $G$-bundle of type $\vep$ and recall the distinguished 
connection $\nabla_0$ of \S\ref{structwist}. 
\begin{proposition} 
\begin{trivlist}
\item\textup{(i)} The identification of the affine hyperplane 
$\mi k^\vee\kstar+L_\vep\frg^* \subset$ with $\cA_\vep$ sending $\mu$ to 
$\nabla_0+\mu/k^\vee$ is equivariant for the action of $L_\vep\frg$.
\item\textup{(ii)} The assignment $\mu \mapsto \dirac +\mi\psi(\mu)$ 
intertwines the affine co-adjoint and commutator actions.
\item\textup{(iii)} Formula (\ref{famsquare}) for $\dirac_\mu^2$,
and its consequences (\ref{kernel}) and (\ref{thom}),
carry over, with $\rho$ replaced by $\ul\rho$.\qed
\end{trivlist}
\end{proposition}

\section{Arbitrary compact groups}\label{arb}

We now extend the construction of the Dirac family, and the resulting 
map from representations to twisted $K$-classes, to the space $\cA_P$ 
of connections on a principal bundle $P$ over the circle, with arbitrary
compact structure group $G$. The Lie algebra $L_P\frg$ of the loop 
group $L_PG$ of gauge transformations splits into a sum of abelian and 
simple loop algebras, and the central extension preserves the splitting
(Prop.~\ref{splitlooplie}). To assemble the families  for the individual 
summands, we must still discuss the abelian case and settle their 
equivariance under the non-trivial components of $L_PG$. 

\subsection{The Abelian case.} Assume, as in Def.~\ref{anreg},
that the the central extension $\widetilde{L}\frz$ takes the form 
$\left[\xi,\eta\right] = b(S\xi,\eta) \cdot K$, for the $L^2$ 
pairing in an inner product on the abelian Lie algebra $\frz$. 
Letting $\widehat{L}\frz := \mi\bR K\oplus L\frz\oplus\mi\bR S$, the 
discussion of \S\ref{affdiracsect} and \S\ref{diracfamsect} carries 
over, with the difference that $\fra = \frz$, $L\frz$ acts trivially 
on the spin module $\spin$, and $\rho$ and $h^\vee$ are null. For 
instance, $\dirac := R_a \otimes \psi^a$, summing over a basis of 
$L\frz$, and relations \eqref{coordfree} in $U_k(L\frz)\otimes\cliff$ 
(cf.~Remark \ref{comarg}) are obvious. 

An admissible irreducible representation of $\widehat{L}\frz$ has 
the form $\fock\otimes\bC_{-\lambda}$, for the Fock representation 
$\fock$ of $\widehat{L}\frz/\frz$ and a $\tau$-affine weight $\lambda$ 
of $\frz$, and we obtain
\[
\dirac^2 = -2S -\lambda^2, \qquad 
\dirac_\mu^2 = -2S -(\lambda-\mu)^2.
\]
\noindent The kernel is identified as before: it is supported on the 
affine subspace $\mi\kstar + \lambda + L\frz^*\ominus\frz^*$ of 
$\widetilde{L}\frz^*$. This is a single co-adjoint orbit of the identity 
component of $LZ$, and the family represents the Thom push-forward of the 
$LZ$-equivariant line bundle $\cO(-\lambda)$, from that orbit to the 
ambient space.

\subsection{Spectral flow over $Z$.}\label{quiso} The positive polarisation 
$\frU \subset L\frz_\bC\ominus\frz_\bC$ of \S\ref{loweightrep} leads to 
vector space identifications $\pspin\cong\spin(0)\otimes \Lambda^\bullet
(\frU)$ and $\fock'\cong \mathrm{Sym}(\frU)$. Decomposing $\dirac_\mu = 
\dirac^\frz_\mu + \dirac^{L\frz/\frz}$ into zero-modes and $\frU$-modes, 
we recognise in the first term is the Dirac family of \S\ref{spectorus}, 
lifted to $\frz^*$ and restricted to the single summand $\bC_{-\lambda}
\subset\fock_{[-\lambda]}$; whereas $\dirac^{L\frz/\frz} = \partial+
\partial^*$, for the Koszul differential
\[
\partial:\mathrm{Sym}^p(\frU)\otimes\Lambda^q(\frU)
\to \mathrm{Sym}^{p+1}(\frU)\otimes\Lambda^{q-1}(\frU).
\]
Thus, $\dirac_\mu$ is quasi-isomorphic to the finite-dimensional family 
$(\bC_{-\lambda},\dirac^\frz_\mu)$ over $\frz^*$. The induced $LZ$-module 
will have the form $\fock'\otimes\fock_{[-\lambda]}$, and dropping the 
factor $\Lambda^\bullet(\frU)\otimes \fock'$, which is equivalent to $\bC$, 
recovers our spectral flow family of \S\ref{fock}.

\subsection{Characterisation of}\label{coordfree}$\dirac_\mu$.
Proposition \ref{equivar} ensures the equivariance of our Dirac family 
for the connected part of the loop group. When $G$ is not simply connected, 
we must extend this to the other components. This is accomplished by an 
intrinsic characterisation of $\dirac_\mu$. We restate relations \eqref
{affdiracrels} and \eqref{twistdiracrels} in a ``coordinate-free'' way:
\[
\left[\dirac_\mu,\psi(\nu) \right] = 2\langle T\,|\, \nu\rangle 
	+ 2\mi\langle \mu\,|\, \nu\rangle \qquad
\left[\dirac_\mu,T(\xi)\right] = \psi\left(\ad_\xi^\vee
( k^\vee\kstar -\mi \mu)\right),
\]
where the bracket in the first equation is contraction
in the bilinear form (\ref{bil}). Observe now that the first 
formula expresses the total action of $\nu$ on $\hilb\otimes\spin$, 
\textit{in the lifting of $L_P\frg$ to $\widetilde{L}_P\frg$ defined 
by the line} $\mi k^\vee\kstar+\mu\subset\widetilde{L}_P\frg^*$. 
In the second formula, we have used the co-adjoint action of \S\ref
{coadjoint}. As explained in Remark \ref{comarg}, the first relation 
uniquely determines $\dirac_\mu$, and we conclude
\begin{proposition} 
The assignment $\mu \mapsto \dirac_\mu$ is equivariant under all 
compatible automorphisms of $L_P\frg$, $\hilb$ and $\spin$ which 
preserve the bilinear form on $L_P\frg$. \qed 
\end{proposition}

\subsection{Coupling to representations.} The Dirac family 
$\dirac_\mu$ lives on an affine copy of $L_P\frg^*$, namely the 
hyperplane over $\mi\in\mi\bR$ in the projection $(L_P\frg^*)^\tau \to 
\mi\bR$, dual to the central extension \eqref{exactLie}. We transport 
it to $\cA_P$ by identifying the two as $L_PG$-affine spaces. For the simple 
factors, this is described in \S\ref{diracfamsect}; but on the abelian 
part, there is an ambiguity: we can translate by the Lie algebra of the 
centre of $L_PG$. Under the holonomy map, this ambiguity matches the 
one encountered in (\ref{ambiguity}.ii), when we identified $\tau\cdot\fra^*$ 
with the space of holonomies. Note, however, that the regularity and 
singularity of the affine weights matches the one of the underlying 
(twisted) conjugacy classes in $G$, irrespective of the chosen 
identification.

Coupling $\dirac_\mu$ to graded, admissible representations results in
twisted $K$-classes on $\cA_P$, equivariant under $L_PG$. This is also 
an Ad-equivariant twisted $K$-classes over $G$, supported on the components 
which carry the holonomies of $P$. 
\begin{proposition}
The isomorphism of Theorem 3 is induced by the Dirac family map, 
from admissible representations to $K$-classes. 
\end{proposition}
\begin{proof}
This follows by comparing the Dirac kernels to the classification of 
irreducibles by their lowest-weight spaces in \S\ref{repclassifsect}, 
and again with the basis of $K^\tau_G(G)$ described in Proposition \ref
{orbits}.
\end{proof}

\part{Variations and Complements}\label{var}

This chapter exploits the correspondence between representations and 
$K$-classes to produce analogues of known constructions in representation 
theory in purely topological terms.

\section{Semi-infinite cohomology} \label{semiinf}

In this section, we give alternative formulae \eqref{alt1}, \eqref{alt2}
for the Dirac operator $\dirac$. With the Lie algebra cohomology results 
of Bott \cite{bott} and Kostant \cite{kost1} and with Garland's loop group 
analogues \cite{gar}, the new formulae explain the magical appearance of 
the kernel on the correct orbit. The relative Dirac operators of \cite
{kost2} and \cite{land} allow us to interpret the morphisms $\omega^*$ 
and $\omega_*$ of \S\ref{topkcompact} in terms of well-known constructions 
for affine algebras, namely \textit{semi-infinite cohomology} and \textit
{semi-infinite induction} \cite{ff}. 

We work here with polynomial loop algebras and lowest-weight modules; for 
simplicity, we omit $f$-twist, underlines and the primes from the notation. 
We shall also use $\ad^\vee$ to denote the co-adjoint action of a Lie 
algebra on its dual, reserving the ``$*$" for hermitian adjoints.

\subsection{Lie algebra cohomology.} The triangular decomposition 
$L\frg_\bC = \frN\oplus\frt_\bC \oplus \ol{\frN}$ factors the spin module 
as $\spin = \spin(\frt^*) \otimes \Lambda^{\bullet}\ol{\frN}^*$. The action 
of $\ol{\frN}$ on a lowest-weight module $\hilb$ leads to a Chevalley 
differential on the Lie algebra cohomology complex, 
\begin{equation}\begin{split} \label{kosdiff}
\bar\partial &: \hilb \otimes \Lambda^k\ol{\frN}^*
\to \hilb \otimes \Lambda^{k+1}\ol{\frN}^*, \\ 
\bar\partial &= R_{-\alpha}\otimes \psi^\alpha + 
	\frac{1}{2}\psi^\alpha\cdot\ad^\vee_{-\alpha}
\end{split}\end{equation}
where we have used a root basis of $\ol{\frN}$ and its dual basis 
$\psi^\alpha$ of Clifford generators. Let $\bar\partial^*$ be the 
hermitian adjoint of $\bar\partial$, and denote by $\dirac^\frt_{-\rho}$ 
the $\frt$-Dirac operator with coefficients in the representation $\hilb 
\otimes \Lambda^{\bullet} \ol{\frN}^*\otimes \bC_{-\rho}$ of $T$. 

\begin{proposition} \label{alt1}
$\dirac = \bar{\partial}+\bar{\partial}^*+\dirac^\frt_{-\rho}$,
and $\dirac^\frt_{-\rho}$ commutes with $\bar{\partial} + 
\bar{\partial}^*$. 
\end{proposition}

\begin{proof} Commutation is obvious. It is also clear that the $R$-terms 
on the two sides agree; so, it remains to compare the Dirac $(\sigma 
\psi)/3$-term in \eqref{affdirac} with ${{\psi\cdot\ad^\vee}/2} 
+ {{(\psi\cdot\ad^\vee )^*}/2}$, plus the ad-term in $\dirac^\frt$. Now, 
all three terms have cubic expressions in the Clifford generators, and 
we will check their agreement. We have
\begin{eqnarray*}
\frac{1}{2}\psi^\alpha\cdot\ad_{-\alpha}^\vee &=& \frac{1}{4}
\sum\nolimits_{\ofrac{\alpha,\beta > 0}{\gamma < 0}} 
{f_{\alpha\beta\gamma}\psi^\alpha\psi^\beta\psi
^\gamma}, \\
\frac{1}{2}\left( {\psi^\alpha\cdot\ad_{-\alpha}^\vee}\right)^*
&=& \frac{1}{4}\sum\nolimits_{\ofrac{\alpha,\beta>0}{\gamma<0}} 
{\bar f_{\alpha\beta\gamma}\psi^{-\gamma}
\psi^{-\beta}\psi^{-\alpha}}.
\end{eqnarray*}
\noindent Disregarding the order of the generators, their difference 
contains precisely the terms in ${\sigma\psi}/3$ involving two positive 
roots and a negative one, respectively two negative roots and a positive 
one; whereas the $\ad$-term in $\dirac^\frt$ similarly collects the 
${\sigma\psi}/3$-terms involving exactly one $\frt^*$-element. Clearly, 
this accounts for all terms in ${\sigma\psi}/3$. We have thus shown 
that the symbols of these operators agree in (a completion of) 
$\Lambda^3(L\frg^*)$.

The difference between the two must then be a linear $\psi$-term.
However, both operators commute with the maximal torus $T$ and with
the energy $E$; so the difference is $\psi(\mu)$, for some $\mu 
\in \frt_\bC^*$. A quick computation gives, for $\nu\in\frt^*$,
\begin{eqnarray*}
\left[\bar{\partial},\psi(\nu)\right] &=& \left[\bar{\partial}^*,
	\psi(\nu)\right] = 0,\\ 
\left[\dirac^\frt_{-\rho},\psi(\nu)\right] &=& 2 T(\nu) = 
\left[\dirac,\psi(\nu)\right];
\end{eqnarray*} 
so $\left[\psi(\mu),\psi(\nu)\right] = 0$ for all $\nu$, and it follows 
that $\mu =0$. 
\end{proof}
 
\subsection{The Dirac kernels.}\label{altproof}
Proposition \ref{alt1} gives a new explanation for the location of 
$\ker\dirac_\mu$. If $\hilb$ is irreducible with lowest weight $(-\lambda)$, 
we have, on $\hilb\otimes \Lambda^q \ol{\frN}^*\otimes\spin(\frt^*)$,
\begin{equation}\label{liecohom}
\ker\left(\bar{\partial} + \bar{\partial}^*\right) \cong 
	H^q \left(\ol{\frN};\hilb\right)\otimes\spin(\frt^*) =
	\bigoplus_{\ell(w) = q} 
	\bC_{w(-\lambda -\rho) +\rho}\otimes\spin(\frt^*),
\end{equation}
\noindent embedded in the Lie algebra complex as harmonic co-cycles; the 
sum ranges over the elements of length $q$ in the affine Weyl group of 
$\frg$. If $\mu \in\frt^*$, then $\dirac_{\mu -\rho}^\frt = \dirac^\frt_
{-\rho} +\mi\psi(\mu)$ commutes with $\left(\bar{\partial} + \bar{\partial}
^*\right)$, so $\ker\dirac_\mu$ is also the kernel of $\dirac_{\mu-\rho}^
\frt$ on \eqref{liecohom}. Clearly, the latter is non-zero precisely when 
$\mu$ is one of the $w(\lambda+\rho)$; otherwise, it follows that the 
highest eigenvalue of $\dirac_\mu^2$ is the negative squared distance 
to the nearest such point, in agreement with \S\ref{ident}.

\subsection{Semi-infinite cohomology.} A similar construction applies 
to a decomposition of a rather different kind. Splitting $L\frg_{\bC} 
= L\frn\oplus L\frt_{\bC} \oplus L\bar{\frn}$ gives a factorisation
\[
\spin(L\frg^*) = \spin(L\frt^*) \otimes \Lambda^{\infty/2+\bullet}
(L\frn^*), 
\]
where the right-most factor is the exterior algebra on the non-negative 
Fourier modes in $\frn^*$ and the duals of the negative ones, the latter 
carrying degree $(-1)$ \cite{fgz,ff}. A formula similar to \ref{kosdiff} 
defines a differential $\partial$ for semi-infinite Lie algebra 
cohomology, acting on $\hilb \otimes \Lambda^{\infty/2} \left(L\frn^*
\right)$. With the same $\hilb$, the semi-infinite cohomology can be 
expressed as a sum of positive energy Fock spaces $\fock \otimes \bC_\mu$ 
for ${L\frt}\left/\frt\right.$, on which $T$ acts with weight $\mu$: 
\begin{equation}
H^{\infty/2+q}\left(L\frn;\hilb \right) = \bigoplus_{\ell(w) = q} 
	\fock \otimes \bC_{w(-\lambda -\rho) +\rho}. 
\end{equation}
Because the splitting of $L\frg_\bC$ was $LT$-equivariant, $LT$ act on 
$\Lambda^{\infty/2}\left(L\frn^*\right)$; it commutes with $\partial$, so 
acts on the cohomology; but the non-trivial components shift the degree. 
Passing to Euler characteristics, we can collect terms into irreducible 
representations $\fock \otimes \fock_{[\mu]}$ of $LT$ (\S\ref{spectorus}) 
and obtain a sum over the finite Weyl group\footnote{In the $f$-twisted 
case, this is the extension $\widetilde{W}^f$ of \eqref{wfext}, and the 
$\fock_{[\mu]}$ are the irreducible $\tau$-modules of $\ul\Pi\times\ul T$.} 
\begin{equation} \label{infindex}
\sum\nolimits_q {(-1)^q H^{\infty/2 + q}\left(L\frn;\hilb\right)} = 
\sum\nolimits_{w\in W} {\vep(w)\cdot \fock \otimes \fock_
	{[w(-\lambda-\rho)+\rho]}}. 
\end{equation} 

\subsection{Relative Dirac operator.}\label{reldir} 
Define $\dirac^{L\frg/L\frt}:=\partial +\partial^*$; its index is given 
by \eqref{infindex}. 

\begin{proposition} \label{alt2}
$\dirac = \dirac^{L\frt}+\dirac^{{L\frg}/{L\frt}}$, and the three
operators commute. \qed
\end{proposition}
\noindent The proof is very similar to the one of Prop.~\ref{alt1}; 
see \cite{land} for more help. Similarly, we have $\dirac_\mu^{L\frg} 
= \dirac^{L\frt}_\mu + \dirac^{{L\frg}/{L\frt}}$, and the three operators 
commute when $\mu \in L\frt^*$. As in \S\ref{altproof}, it follows that 
the restriction to $L\frt^*$ of our Dirac family on $\hilb \otimes \spin$ 
is stably equivalent to $\dirac_\mu^{L\frt}$, acting on the alternating 
sum in \eqref{infindex}. Comparing this with the construction \eqref
{orbits} of $K$-classes from conjugacy classes and with the local model 
of the Weyl map \eqref{weylmap}, we obtain the following

\begin{theorem} \label{seminfres}
Under the Dirac family construction Chapter \ref{diracfam}, the 
semi-infinite $L\frn$-Euler characteristic, from $R^{\tau-\sigma}(LG)$ 
to $R^\tau(LT)$, corresponds to the Weyl restriction 
$\omega ^*:K^\tau_G(G)\to K^\tau_T(T)$. \qed
\end{theorem}

\begin{remark} \begin{trivlist}\itemsep0ex
\item (i) In the twisted case, this applies to the restriction 
$\omega ^*:K^\tau_{G(f)}(fG_1)\to K^\tau_T(fT)$.
\item (ii) We have used $LT$ for simplicity, but the result applies to 
$LN$, which preserves the relative Dirac $\dirac^{L\frg/L\frt}$ (though 
not the semi-infinite differential $\partial$). We then detect the 
restriction to $K^\tau_{N(f)}(fT)$.
\end{trivlist}\end{remark}

\section{Loop rotation, energy and the Kac numerator}\label{energysect}

\subsection{Conditions for rotation-equivariance.}
The admissible loop group representations of greatest interest 
admit a circle action intertwining with the loop rotations (\S\ref
{looprot}). This will be the case iff the following two conditions 
are met:
\begin{trivlist}\itemsep0ex
\item(i) The loop rotation action lifts to the central extension $LG^\tau$,
\item(ii) The polarisation used in defining admissibility is 
rotation-invariant (\S\ref{admissibility}). 
\end{trivlist}
A lifting in (i) defines a semi-direct product $\bT\ltimes LG^\tau$. Subject 
to condition (ii), the Borel-Weil construction \cite{psloop} shows that all 
of admissible representations carry actions of the identity component of 
this product, and the $\bT$-action is determined up to an overall shift on 
each irreducible. The action can be extended to the entire loop group as 
in \S\ref{repclassifsect}, and this leads to the same classification 
of irreducibles, but with the extra choice of normalisation for the 
circle action.

With respect to condition (i), it is convenient to allow \textit{fractional} 
circle actions: that is, we allow the circle of loop rotations to be 
replaced by some finite cover. A lifting of the rotation action to 
$LG^\tau$ refines the level $[\tau]$ to a class in $H^3\left(B(\bT
\ltimes LG)\right)$. The obstruction to such a refinement is the 
differential $\delta_2: H^3_G(G_1)\to H^2(B\bT)\otimes H^2_G(G_1)$ in 
the Leray sequence for the projection to $B\bT$. All torsion obstruction 
vanish when $\bT$ is replaced by a suitable finite cover. Rationally, 
$H^*_G(G_1)$ is the invariant part of $H^*_T(T)$ under the Weyl group $W$ 
of $G$, and for the torus we have the following.

\begin{lemma}\label{symcond}
A class in $H^3(T\times BT)$ lifts to a rotation-equivariant one
iff its component in $H^1(T)\otimes H^2_T$ is symmetric.
\end{lemma}
\begin{proof} The differential $\delta_2$ vanishes on the $H^*(T)$ 
factor, and is determined its effect on $H^2_T$: this is mapped 
isomorphically onto $H^2(B\bT)\otimes H^1(T)$. On $H^3(T\times BT)$, 
this becomes the anti-symmetrisation map $H^1(T)\otimes H^2_T \to H^2(T) 
=\Lambda^2H^1(T)$.
\end{proof}
\begin{remark} \label{ssfine}
For semi-simple $G$, symmetry is ensured by Weyl invariance. 
\end{remark}

Adding loop rotations to the landscape leads to the quotient stack of 
the space $\cA$ of smooth connections by the action of $\bT\ltimes LG$. 
This is a smooth stack, with compact quotient and proper stabiliser, 
but which, unlike the quotient stack $G_G$ of $G$ by its own $\Ad
$-action, cannot be presented as a quotient of a manifold by a compact 
group. The $K$-theory of such stacks is discussed in \cite{fht1}. Let 
$\widehat{\wtlat}^\tau = \wtlat^\tau \oplus\bZ\delta$ be the level $\tau$ 
slice of the affine weight lattice \eqref{afroot}.

\begin{proposition}
We have isomorphisms $R^{\tau-\sigma} (\bT\ltimes LG) \cong 
K^{\tau'-\sigma(\frt)}_{\waff}(\widehat{\wtlat}^\tau)\cong 
K^{\tau+\dim\frg}_\bT(G_G)$, obtained by tracking the loop rotation 
in Thm.~\ref{kcompute} and in the Dirac family. 
\end{proposition}

\noindent The middle group is a free $R_\bT$-module, with the generator 
acting on $\widehat{\wtlat}^\tau$ by $\delta$-translation. Killing the 
augmentation ideal forgets the circle action in the outer groups and 
$\delta$ in the middle group, and recovers the isomorphisms in Theorems~3 
and 4.

\begin{proof}The argument is a repetition of \eqref{keylemma}, \eqref
{toruskcompute} and \eqref{kcompute}, with the extra $\bT$-action. 
The main difference is that we are now dealing with the $K$-theories 
of some smooth, proper stacks, which are no longer global quotients, 
but only locally so. However, the proofs of \eqref{keylemma} and \eqref
{kcompute} proceed via a local step, which continues to apply, globalised 
using the Mayer-Vietoris principle.
\end{proof}

\subsection{Positive energy.}\label{poserg} 
The natural choices for the Fredholm operator $S$ defining the Lie algebra 
cocycle in \eqref{anreg} are multiples of the derivative $-\mi d/dt$; the 
polarisation $\frP$ is then the semi-positive Fourier part of $L\frg_{\bC}$. 
With those choices, lowest-weight modules of $L\frg$ carry a bounded-below 
energy operator $E$, unique up to additive normalisation, generating the 
intertwining loop rotation action. If the restriction to $H^1(T)\otimes 
H^2_T$ of $[\tau]$ is symmetric, loop rotations lift fractionally to 
$LG^\tau$; and, if that same bilinear form is positive, $E$ is bounded below 
on admissible $\tau$-representations of the group.

This generalises easily to the \textit{twisted} loop groups $L_PG$ of 
gauge transformations of a principal bundle $P$ over $S^1$. The 
diffeomorphisms of the bundle $P$ which cover the loop rotation form an 
extension of the rotation group $\bT$ by $L_PG$; any connection on $P$ 
whose holonomy has finite order gives a fractional splitting of this 
extension. This group replaces the $\bT\ltimes LG$ of the trivial bundle 
case. The topological constraint for rotation-equivariance of a extension 
$\tau$ is now the symmetry of the map $\kappa^\tau$ in \S\ref{affaction}.  

\subsection{The Kac numerator.} \label{kachar}
For the remainder of this section, we make the simplifying assumption that 
$G$ is connected, with $\pi_1G$ free. Positive energy representations of 
$LG$ are then determined by their restriction to the subgroup $\bT\times G$ 
of circle rotations and constant loops; moreover, loop rotations extend to 
a trace-class action of the semi-group $\{q\in\bC^\times | |q|<1\}$. 
If $\hilb$ is irreducible with lowest-weight $(-\lambda)$, the value of 
its character at $q\in\bC$ and $g\in G$ is given by the \textit{Kac 
formula} \cite{kac}
\begin{equation}\label{charform}
\mathrm{Tr}\left(qg|\hilb\right)= \frac{
	\sum\nolimits_{\mu}\vep(\mu)\cdot q^{\|\mu\|^2/2}\cdot 
		\mathrm{Tr}\left(g | V_{\rho-\mu}\right)}{\Delta(g;q)},
\end{equation}
where $\mu$ ranges over the dominant regular affine Weyl transforms of 
$(\lambda+\rho)$ at level $[\tau]$, $\vep(\mu)$ is the signature of the 
transforming affine Weyl element, $V_\mu$ the $G$-representation with 
lowest weight $\mu$, $\|\mu\|^2 := \langle(\kappa^\tau)^{-1}(\mu) |\mu\rangle$ 
defined by the level $[\tau]$, and the \textit{Kac denominator} for 
$(L\frg,\frg)$
\[
\Delta(g;q) = \prod_{n>0}\det\left(1-q^n\cdot\ad(g)\right)
\] 
independent of $\lambda$ and $\tau$, representing the (super)character 
of spinors on $L\frg/\frg$. We shall now see how \eqref{charform} is 
detected by our $K_\bT$-group.

Including the identity $e\in G$ defines a Gysin map
\[
\mathrm{Ind}: R^{\tau-\sigma(\frg)}(\bT\times G) \to 
		K^{\tau+\dim G}_\bT(G_G),
\]
with $\tau$ on the left denoting the restricted twisting and $\sigma(\frg)$ 
the Thom twist of the adjoint representation. Dualising over $R_\bT$, while 
using the bases of irreducible representations to identify $K^\tau_\bT(G_G)$ 
with its $R_\bT$-dual, leads to an $R_\bT$-module map 
\[
\mathrm{Ind}^*: K^{\tau+\dim G}_\bT(G_G)\to \mathrm{Hom}_\bZ\left(
		 R^{\tau-\sigma(\frg)}(G); R(\bT)\right); 
\]
the right-hand side is the $R(\bT)$-module of formal sums of (twisted) 
$G$-irreducibles with Laurent polynomial coefficients. The choice of basis 
gives an indeterminacy of an overall power of $q$ for each irreducible, 
which must be adjusted to give an exact match in the following theorem. 
Let $[\hilb]$ be the $K^\tau_\bT(G_G)$-class corresponding to $\hilb$. 

\begin{theorem}
$\mathrm{Ind}^*[\hilb]$ is the Kac numerator in \eqref{charform}.
\end{theorem}
\begin{proof}
The theorem is a consequence of two facts. First is the relation 
\begin{equation}\label{indrel}
q^{\|\lambda+\rho\|^2/2}\cdot\mathrm{Ind}\left(V_{-\lambda}\right) = 
	\vep(\mu)\cdot q^{\|\mu\|^2/2}\cdot\mathrm{Ind}\left(V_{\rho-\mu}\right),
\end{equation}
holding for any $\mu$ in the affine Weyl orbit of $(\lambda+\rho)$. Second 
is the fact that, with our simplifying assumption that $G$ is connected 
with free $\pi_1$, the twisted $K$-class $\mathrm{Ind}(V_{-\lambda})$ 
corresponds to an irreducible representation of $LG^\tau$. (There are no 
affine Weyl stabilisers of regular weights). 

We can check \eqref{indrel} by restriction to the maximal torus $T$. 
The Weyl denominator is the Euler class of the inclusion $T\subset G$; 
multiplying by it while using the Weyl character formula converts the Kac 
numerator for $(L\frg,\frg)$ to that of $(L\frt,\frt)$, and we are reduced 
to verifying the theorem for the torus (with $\vep(\mu) =1$ and without 
$\rho$-shifts, as the affine Weyl group is now the lattice $\pi_1T$). 

The twisting $\tau$ defines a line bundle $L$ over the representation ring 
of the stabiliser over $T$.\footnote{This $L$ is the free $\bZ$-module over
the covering $Y$ in Construction \eqref{constr}.} The stabiliser itself is 
a bundle of groups with fibre $\bT\times T$ and holonomy around a loop $\gamma\in\Pi$ given by the automorphism 
\[
q^m t^{\mi\lambda} 
	\mapsto q^{m+\langle\lambda |\gamma\rangle}t^{\mi\lambda},
\] 
where $q\in\bT$, $t\in T$ and $\lambda:\pi_1T\to\bZ$ is an integral weight. 
For $L$, this must vary by multiplication by a unit $q^{\phi(\lambda,
\gamma)}\cdot t^{\mi\kappa^\tau(\gamma)}$. (The exponent $\kappa^\tau
(\gamma)$ of $t$ is detected by restriction to $q=1$.) We claim that the 
only option, up to automorphism, is $\phi(\lambda,\gamma) = \langle\kappa^\tau(\gamma)| \gamma\rangle/2$, resulting in the holonomy 
\[
q^{m+\|\lambda\|^2/2} t^{\mi\lambda} \mapsto 
	q^{m+\|\lambda +\kappa^\tau(\gamma)\|^2/2}\cdot 
		t^{\mi(\lambda+\kappa^\tau(\gamma))};
\] 
travelling around $\gamma$ shows that induction from the characters  
$q^{\|\lambda\|^2/2} t^{\mi\lambda}$ and $q^{\|\lambda +\kappa^\tau
(\gamma)\|^2/2}\cdot t^{\mi(\lambda+\kappa^\tau(\gamma))}$ of the stabiliser 
must lead to the same twisted $K$-class, which proves \eqref{indrel} 
for the torus and hence our theorem. To check the claim, note the two 
relations
\[\begin{split}
\phi(\lambda+\mu,\gamma) &= \phi(\mu, \gamma), \\
\phi(\lambda,\gamma+\gamma') &= \phi(\lambda,\gamma) + 
		\phi(\lambda,\gamma') + \langle\kappa^\tau(\gamma)|\gamma'\rangle
\end{split}
\]
the first, by computing the holonomy of $t^{\mi(\lambda+\mu)} = 
t^{\mi\lambda} t^{\mi\mu}$ in two different ways (using the module structure 
of $L$) and the second, from the homomorphism condition. These imply that 
$\phi(\lambda,\gamma) = \langle\kappa^\tau(\gamma)| \gamma\rangle/2$,  modulo 
a linear $\gamma$-term; but the latter can be absorbed by a shift 
$t^{\mi\lambda}\mapsto t^{\mi(\lambda+\nu)}$ in $T$-characters, representing 
an automorphism of $L$.
\end{proof}
\begin{remark}
This can be generalised to twisted loop groups and their disconnected 
versions, but to determine a representation uniquely, we must restrict 
it to a larger subgroup of the loop group, meeting every torsion 
component in a translate of the maximal torus. We expect in that case 
to recover the extension of the Kac character due to Wendt \cite{wend}.
\end{remark}

\section{Fusion with $G$-representations}
\label{fusion}
For positive energy representations, the \textit{fusion product} of 
conformal field theory defines an operation $*: R(G)\otimes R^\tau(LG) 
\to R^\tau(LG)$. We will now recall its construction and prove its 
agreement with the $R(G)$-action on $K^\tau_G(G)$. For notational 
clarity, we treat the untwisted loop groups, the twisted result 
following by judicious insertion of underlines and $f$-subscripts. 

\subsection{Example: $G_1$ is a torus.} \label{torusfusion}
Recall from \S\ref{anregsect} that, when $G=N$, $LN\cong \Gamma N\ltimes
\exp(L\frt\ominus\frt)$, where $\Gamma N = \naff$ is the subgroup of 
geodesic loops. Evaluating geodesic loops at a point $x$ in the circle 
gives a homomorphism $E_x:LN \to N$. If $V$ is a finite-dimensional 
$N$-representation, the pull-back $E_x^*V$ is an admissible $LN
$-representation, and fusing with $V$ is simply tensoring with $E_x^*V$. 

Note that $h$ is not the ``evaluation at $x$" homomorphism on the whole of 
$LN$; the latter would \emph{not} lead to admissible representations. 
For non-abelian $G_1$, we need the more complicated definition that follows, 
essentially moving the base-point $x$ inside the disk.  

\subsection{Segal's holomorphic induction.} Let $\hilb$ be a positive energy 
admissible $\tau$-representation of $LG$, and $V$ a $G$-module whose 
$\rho$-shifted highest weights lie in the alcove $\tau\cdot\fra^*$ (\S\ref
{morewaff}). Such $G$-modules are called \textit{small}. Let also $A$ 
be a complex annulus, with an interior base-point $x$. The obvious group 
$\mathrm{Hol}(A;G_\bC)$ of holomorphic maps with smooth boundary values 
acts on $\hilb$, by restriction to the inner boundary, on $V$ by evaluation 
at $x$, and maps into a copy of $LG_\bC$ by restriction to the outer 
boundary. G.~Segal defines the fusion of $\hilb$ with $V$ along $A$ as 
the \textit{holomorphic induction}
\begin{equation}\label{holind}
\hilb * V_x := 
	\mathrm{Ind}_{\mathrm{Hol}(A;G_\bC)}^{LG_\bC}\left(
					\hilb\otimes V_x\right),
\end{equation}
by which we mean the space of right $\mathrm{Hol}(A;G_\bC)$-invariant 
holomorphic maps from $LG_\bC$ to $\hilb\otimes V$. Conjecturally, this 
is a completion of an admissible representation.

The rigorous implementations of this construction that we know are algebraic. 
The direct product $\widehat\hilb$ of energy eigenspaces in $\hilb$ is a 
representation of the \textit{Laurent polynomial loop group} $L'G_\bC := 
G_\bC[z,z^{-1}]$. After evaluation at $z=x$, $L'G_\bC$ also acts on $V$. 
The completion of $L'G_\bC$ at $z=\infty$ is the group of formal Laurent loops 
$G_\bC(\!(w)\!)$ ($w=z^{-1}$). Its algebraic, positive energy $\tau-\sigma
$-modules are completely reducible, and the irreducibles are precisely 
the direct sums $\philb$ of energy eigenspaces in irreducible admissible 
representations $\hilb$ of $LG$.\footnote{Experts will know that, when $G$ 
is not semi-simple, these algebraic loop groups are highly non-reduced group 
(ind)-schemes, and their formal part must be included in the discussion.} 
Constructing the induced representation now from \textit{algebraic} functions, 
the following important lemma permits the subsequent definition. 

\begin{lemma}\label{posind}
$\mathrm{Ind}_{L'G_\bC} ^{G_\bC(\!(w)\!)}\left(\widehat{\hilb}\otimes 
V\right)$ is a finitely reducible, positive-energy representation 
of $G_\bC(\!(w)\!)$. 
\end{lemma} 

\begin{definition} 
The fusion product $\hilb * V_x$ is $\mathrm{Ind}_{L'G_\bC}
^{G_\bC(\!(w)\!)}\left( \widehat{\hilb}\otimes V\right)$.
\end{definition}  

\noindent Using brackets to denote the associated $K$-classes, the 
fusion is identified by the following

\begin{theorem}\label{tensorfusion}
In $K^\tau_G(G)$ with its topological $R(G)$-action, $[\hilb * V_x] = 
[\hilb]\otimes[V]$.
\end{theorem}

\noindent The proof of this theorem requires some preliminary constructions.

\subsection{Borel-Weil construction.} 
We need to review the construction of $\hilb$ by algebraic induction 
from a Borel-like subgroup (called the Iwahori subgroup), but minding 
the group $\pi_0LG$ of components.\footnote{To see the problem, recall 
that every representation of a \textit{connected} compact Lie group is 
holomorphically induced from a Borel subgroup $B$. This fails in the 
disconnected case, where induction from the quasi-Borel $Q_T\cdot B$ is 
required instead.} This is neatly accomplished by a construction due to 
Beilinson and Bernstein. 

The \textit{quasi-iwahori subgroup} $\QB\subset L'G_\bC$ is the normaliser 
of $\frN$; it meets every component of $L'G_\bC$ in a translate of the 
standard Iwahori subgroup. Killing $\exp(\frN)$ converts $\QB$ into a 
subgroup $Q_L\subset(\naff)_\bC$, which plays the r\^ole of a (complex) 
quasi-torus for the loop group. $Q_L$ is the normaliser of $\frN$ in 
$(\naff)_\bC$. There is a Cartesian square
\[
\begin{array}{ccc}
Q_L & \longrightarrow & \naff\\
\downarrow & & \downarrow \\
\pi_0LG & \longrightarrow &\waff\:,
\end{array}\] 
where the bottom horizontal arrow is the splitting of \eqref{wext} defined
by the positive alcove.

Call $\cU$ the algebraic vector bundle over the full flag variety $X':= 
L'G_\bC/\QB$ whose fibre at a coset $\gamma\QB$ is the space $\philb/
\frN^\gamma\philb$ of co-invariants in $\philb$, with respect to the 
conjugated nilpotent $\frN^\gamma:=\gamma\frN\gamma^{-1}$. (This fibre 
is isomorphic to the lowest-weight space for the opposite polarisation.) 
Then, $\widehat\hilb$ is the space of algebraic sections of $\cU$ over $X'$. 
A result of Kumar \cite{kumar} ensures the vanishing of higher cohomologies 
of this bundle.

\begin{remark}\begin{trivlist}\itemsep0ex
\item (i) $\QB$ acts (projectively) on the space $U:=\philb/\frN\philb$, 
which defines a projective $L'G_\bC$-vector bundle over $X'$; ``unprojectivising" this bundle at level $\tau-\sigma$ results in $\cU$.
\item (ii) The same prescription defines $\cU$ over the ``thicker" flag 
variety $X:= G_\bC(\!(w)\!)/\QB$, and its sections there lead to the ``thin" 
version $\philb$ of the same representation.
\end{trivlist}\end{remark}

\subsection{Derived induction.} 
The fibre of $\cU$ at $1$ is a representation of $\QB$ which factors 
through $Q_L$, and whose highest weights are in $\wtlat^\tau_{reg+}$, as 
discussed in \S\ref{repclassifsect}. We now study the ``derived 
induction" $\mathbf{R}\mathrm{Ind}$ from $Q_L$-modules to $LG$-modules, 
by which we mean the Euler characteristic over $X$ of a vector bundle 
associated to a general $(\tau-\sigma)$-module of $Q_L$. By \S\ref
{repclassifsect} again,
\[
{}^{\tau-\sigma} R(Q_L)\cong 
	{}^{\tau'-\sigma(\frt)} K_\pi(\wtlat^\tau),
\]
with the action and twistings defined there, and we claim that $\mathbf{R} 
\mathrm{Ind}$ is the result of the direct image map, followed by restriction 
to the regular part:
\begin{equation}\label{Rind}
{}^{\tau'-\sigma(\frt)} K_\pi(\wtlat^\tau) \to
{}^{\tau'-\sigma(\frt)} K_{\waff} (\wtlat^\tau) \to
{}^{\tau'-\sigma(\frt)} K_{\waff} (\wtlat^\tau_{reg}). 
\end{equation}
From \S\ref{repclassifsect} and the vanishing of higher cohomology, this 
is known for weights in $\wtlat^\tau_{reg+}$. Because $\waff\cong 
\pi\ltimes\waff(\frg)$ and $\tau\cdot\fra^*$ is a fundamental domain for 
$\waff(\frg)$, it suffices to show that $\mathbf{R}\mathrm{Ind}$ is 
anti-symmetric under this last group and that weights on the walls of 
$\tau\cdot\fra^*$ induce $0$. Both statements follow from Bott's reflection 
argument \cite{bott} applied to the simple affine reflections.

\begin{proof}[Proof of (\ref{posind})]
$\QB$ acts on on $V$ by evaluation at $z=x$; calling $\cV_x$ the associated 
vector bundle over $X$, transitivity of induction shows that
\[
\mathrm{Ind}_{L'G_\bC} ^{G_\bC(\!(w)\!)}
		\left(\widehat{\hilb}\otimes V\right) \cong
	\Gamma\left(X; \cU\otimes \cV_x\right). 
\]
and the Lemma now follows from Theorem 4 of \cite {tel}.
\end{proof}

\begin{proof}[Proof of (\ref{tensorfusion})]
Theorem 4 of \cite {tel} also ensures the vanishing of higher cohomologies 
when $V$ is small. We will identify $\hilb * V_x$ by deforming $\cV_x$. 
Scaling $x\mapsto 0$ deforms the action of $\QB$ on $V_x$ into the 
representation $V_0$, pulled back from the quotient map $\QB\to Q_L$. 
More precisely, any point-wise evaluation $LG \to G$ embeds $Q_L$ into 
$N$, and $V_0$ is obtained from $V$ under $Q_L\to N\subset G$. The Euler 
characteristic of the bundle $\cU\otimes\cV_x$ is unchanged under 
deformation, because of the rigidity of admissible representations of 
$G(\!(w)\!)$,\footnote{And the techniques of \cite{tel}, which reduce this 
to a finite type problem.} and we conclude that 
\[
\hilb * V_x \cong \mathbf{R}\mathrm{Ind}(\cU\otimes V_0). 
\]
To prove the theorem, we must show that $\mathbf{R}\mathrm{Ind}: R(Q_L) 
\to K^\tau_G(G)$ is an $R(G)$-module map, under the inclusion $Q_L\subset G$. 
Factoring $\mathbf{R}\mathrm{Ind}$ as in \eqref{Rind}, this property is 
clear for the second step, restriction to $\wtlat^{reg}$, since that is 
nothing but the map $\omega^*$ of \S\ref{kcompute}. A different description 
makes the same obvious for the first step, the direct image. Indeed, 
\[
{}^{\tau'-\sigma(\frt)} K_{\waff} (\wtlat^\tau)\cong 
	K^\tau_{\naff}(\frt),\qquad
{}^{\tau'-\sigma(\ul\frt)} K_\pi(\wtlat^\tau) \cong
	K^\tau_{Q_L}(\frt),
\]
as in Remark \ref{otherpicture}. The direct image map becomes now 
induction along the inclusion $Q_L\subset\naff$, and this is clearly 
a module homomorphism under $R(N)$, to which $\naff$ maps by evaluation. 
\end{proof}

\section{Topological Peter-Weyl theorem}
\label{peterweylsection}

We now describe a topological version of the Peter-Weyl theorem for loop 
groups; beyond its entertainment value, the result can be used to confirm 
that the bilinear form in the Frobenius algebra $K^\tau_G(G)$ of \cite
{fht1} agrees with the natural duality pairing in the Verlinde ring, as 
claimed in \cite[\S8]{fht3}.\footnote{This interpretation is only 
available for twistings that are transgressed from $BG$.} 

\subsection{Compact groups.} \label{pwfinite}
One version of the Peter-Weyl theorem asserts that the two-sided regular 
representation of a compact group $G$ --- the space of functions on $G$, 
under its translation actions on the left and on the right --- is (a 
topological completion of) the direct sum $\bigoplus V\otimes V^*$, ranging 
over the irreducible finite-dimensional modules $V$.\footnote{For a Lie 
group, the direct sum describes the polynomial functions.} A variation 
of this, for a central extension $G^\tau$ by $\bT$, describes the space 
of sections of the associated line bundle over $G$ as the corresponding 
sum over irreducible $\tau$-representations. 

Qua $G\times G$-module, the regular representation of $G$ is induced from 
the trivial $G$-module, under the diagonal inclusion $G\subset G\times G$. 
For \textit{finite} $G$, the result can be expressed in terms of equivariant 
$K$-theory: it asserts that the trivial representation $[1] \in R(G)$ maps, 
under diagonal inclusion, to the class $\sum [V\otimes V^*] \in R(G\times 
G)$. To see this topological direct image more clearly, identify $R(G)$ 
with $K_{G\times G}(G)$, with the left$\times$right action, and push 
forward to a point with $G\times G$ action. In the presence of a twisting 
$\tau$ for $R(G)$, we map $[1]\in R(G)$ instead to $R^{\tau\times(-\tau)}
(G\times G)$. In constructing this last push-forward, we have used the 
natural trivialisation of the sum of a central extension $\tau$ of $G$ 
with its opposite, so that the required twisting on $K_{G\times G}
(G)$ is canonically zero. 

\begin{remark}
When $\tau$ is graded, our formulation of Peter-Weyl conceals a finer point. 
The module $R^\tau(G)$ of graded representations has now an odd component 
$R^{\tau+1}(G)$, defined from the \textit{super-symmetric representations} 
\cite[\S4]{fht3}. These are graded $G$-modules with a commuting action 
of the rank one Clifford algebra $\cliff(1)$. The contribution of such a 
super-symmetric representation $V$ to the Peter-Weyl sum is the (graded) 
tensor product $V\otimes_{\cliff(1)} V^*$ over $\cliff(1)$, and not over 
$\bC$. However, this \textit{does} match the cup-product 
\[
R^{\tau+1}(G)\otimes R^{-\tau+1}(G) \to R^{\tau\times(-\tau)+0}
	(G\times G);
\]
indeed, the (graded) tensor product $V\otimes_\bC V^*$ has a commuting 
$\cliff(2)$ action, and defines an element of $K^2$, which is indeed where 
the cup-product initially lands \cite{mich}. Tensoring over $\cliff(1)$ 
instead of $\bC$ is the Morita identification of complex $\cliff(2)$-modules 
with vector spaces, which implements the Bott isomorphism to $K^0$. 
\end{remark}

\subsection{Loop groups.} 
Before discussing the loop group analogue of this, let us recall the algebraic 
Peter-Weyl theorem for loop groups, a special case of the Borel-Weil theorem 
of \cite{tel}. As in the preceding section, denote by $G_\bC(\!(z)\!)$ and 
$G_\bC(\!(w)\!)$ be the two Laurent completions of the loop group $LG$ at 
the points $0$ and $\infty$ on the Riemann sphere. The Laurent polynomial 
loop group $L'G_\bC= G[z,z^{-1}]$ embeds in both (with $w =z^{-1}$). The 
quotient variety $Y:= G(\!(w)\!)\times_{L'G_\bC} G(\!(z)\!)$ for the diagonal 
action is a homogeneous space for the product of the two loop groups, which 
should be thought regarded as a generalised flag variety. For any twisting 
$\tau$, the product $\cO(\tau-\sigma)\boxtimes\cO(\sigma-\tau)$ of the 
opposite line bundles on the two factors carries an action of $L'G_\bC$, 
so it descends to an (algebraic) line bundle on $Y$. A special case of the 
Borel-Weil-Bott theorem of \cite{tel} asserts that, as a representation 
of $G(\!(w)\!)\times G(\!(z)\!)$,
\[
\Gamma\left(Y; \cO(\tau-\sigma)\boxtimes\cO(\sigma-\tau)\right) \cong 
\bigoplus\nolimits_\hilb \hilb'\otimes\overline{\hilb}',
\] 
with the sum ranging over the lowest-weight representations $\hilb$ for 
$G_\bC(\!(w)\!)$ at level $\tau-\sigma$.

\subsection{Topological interpretation.}
The topological construction in \S\ref{pwfinite} breaks down for infinite 
compact groups, but remarkably, it does carry over to loop groups. To start 
with, the diagonal self-embedding of $G$ leads to a Gysin map 
\[
\iota_*: K_G(G) \to K_{G\times G}(G\times G),
\] 
with the $\Ad$-action in both cases. This models the diagonal $LG\to LG 
\times LG$ when $G$ is \textit{connected}. In general, the restriction 
of $\iota_*$ to $K_G(G_1)$ corresponds to $LG$, whereas its restriction 
to $K_{G(f)} (fG_1)$, as in \S\ref{weylmap}, captures the diagonal 
embedding for the twisted loop group $L_fG$. Similarly, for any $\tau$, 
we get a map 
\begin{equation}\label{iota}
\iota^\tau_*: K_G(G) \to K^{\tau\times(-\tau)}_{G\times G}(G\times G), 
\end{equation}
cancelling the pulled-back twisting by the same observation as before: 
the sum of extensions $\tau + (-\tau)$ is trivial on the diagonal $LG$.
 
\begin{remark}
To construct $\iota_*$, we replace $K_G(G)$ with the isomorphic group 
$K^*_{G\times G}(G\times G)$, the action being now 
\[
(g_1,g_2).(x,y) = (g_1xg_1^{-1},g_1yg_2^{-1}).
\]
The isomorphism with $K^*_G(G)$ arises by restriction to the diagonal $G$'s. 
The map $G\times G\to G\times G$ inducing $\iota_*$ sends $(x,y)$ to 
$(x, y^{-1}xy)$. Note that the relative tangent bundle of this map is
(stably) equivariantly trivial, and there is a preferred relative orientation, 
if we use the same dual pair of Spin modules on each pair of $\frg$'s.
\end{remark}

\begin{theorem}[Peter-Weyl for Loop Groups] \label{pw}
When $\tau$ is regular, we have
\[
	\iota^\tau_*(1) = \sum\nolimits_\hilb [\hilb\otimes\hilb^*], 
\]
summing over the irreducible admissible representations $\hilb$ of 
$LG$, in the correspondence of Theorem~3. The analogue holds for 
each twisted loop group of $G$.
\end{theorem}

\noindent Without using Theorem~3, we can assert that $\iota_*(1)$ has a 
diagonal decomposition in the basis of $K^\tau_G(G)$ produced from 
regular affine Weyl orbits and irreducible representations of the 
centralisers (Theorem~\ref{kcompute}), and the complex-conjugate basis 
for $K^{-\tau}_G(G)$. The two formulations are of course related 
by Theorem~\ref{repclassif}. To state this more precisely, we use the 
``anti-diagonal" class $[\Delta^-]$ on $\ul{\wtlat}^\tau \times
\ul{\wtlat}^{-\tau}$, which is identically $1$ on pairs $(\lambda,-\lambda)$ 
and null elsewhere. It is equivariant for the diagonal $\waff$-action. Also 
let $\tau'' = \tau'-\sigma(\ul\frt)$. 

\begin{lemma}\label{diagweights}
The sum in (\ref{pw}) corresponds to the direct image of $[\Delta^-]$ under 
the direct image map  
\[
K_{\waff}\left(\ul{\wtlat}^\tau_{reg} 
		\times\ul\wtlat^{-\tau}_{reg}\right) 
\longrightarrow K^{\tau"\times(-\tau")}_{\waff\times\waff}
		\left(\ul{\wtlat}^\tau_{reg} 
		\times\ul\wtlat^{-\tau}_{reg}\right).
\] 
\end{lemma}

\begin{proof} Replacing both sides with the sets of orbits, represented 
by weights $\mu\in\ul{\wtlat}^\tau_+$ and stabilisers $\pi_\mu\subset\waff$, 
we get the direct sum over $\mu$ of the diagonal push-forwards 
\[
R(\pi_\mu) \to R^{\tau"}(\pi_\mu) \otimes R^{-\tau"}(\pi_\mu), 
\]
and apply the topological Peter-Weyl theorem to each $\pi_\mu$. 
\end{proof}

For a torus $T$, the representation categories of $LT^\tau$ and $\Gamma
^\tau =(\Pi\times T)^\tau$ are equivalent, and $\iota^\tau_*$ captures the 
Peter-Weyl theorem for $\Gamma^\tau$: diagonal induction of the trivial 
representation to $\Gamma^\tau\times\Gamma^{-\tau}$ leads to the sum 
in \eqref{pw}. This result generalises to every group $\naff$ of 
($f$-twisted) geodesic loops in $N$, and is the basis for the general 
proof. To convert it into a topological statement, we will factor both the algebraic and the topological induction (direct image) maps into two steps, 
with the second step being described by Lemma~\ref{diagweights}. Agreement 
of the other, first step is then verified by a Dirac family construction 
akin to \S\ref{spectorus}. As the general case may be obscured by the 
requisite notational clutter, we handle the torus first. 

\subsection{Example: $G = T$.}  
Let $\ell = \dim T$ and factor $\iota^\tau_*$ into the direct images 
\begin{equation}\label{178}
K^0_T(T) \xrightarrow{B\mathrm{diag}_*} 
	K^{\tau\times(-\tau)-\ell}_{T\times T}(T) 
	\xrightarrow{\mathrm{diag}_*}
	K^{\tau\times(-\tau)+0}_{T\times T}(T\times T),
\end{equation}
along the obvious diagonal morphisms. Describing $\mathrm{diag}_*$ is easy. 
Double use the Key Lemma~\ref{keylemma}, with the \textit{same} group $T^2$, 
followed by direct images (along $\frt$ and $\frt^2$) lead to isomorphisms
\begin{equation}\label{179}\begin{split}
K^{\tau\times(-\tau)-\ell}_{T\times T}(T) 
	&\cong K^0(\wtlat^\tau\times_\Pi\wtlat^{-\tau}), \\
K^{\tau\times(-\tau)+0}_{T\times T}(T\times T)
	&\cong K^0(\wtlat^\tau/\Pi\times\wtlat^{-\tau}/\Pi).
\end{split}\end{equation}
Moreover, $\mathrm{diag}_*$ becomes the direct image between the groups 
on the right, which is the map in Lemma \ref{diagweights}.

In view of Lemma~\ref{diagweights}, we must check that $B\mathrm{diag}_*[1]$ 
in the middle group of \eqref{178} is the anti-diagonal $[\Delta^-]$. 
We have a commutative square 
\[\begin{array}{ccc}
K^0_T(T) & \xrightarrow{B\mathrm{diag}_*} & 
		{}^{\tau\times(-\tau)}K^{-\ell}_{T\times T}(T) \\
 \uparrow & \Box  & \uparrow  \\
K^0(T) & \xrightarrow{p_*}  & K^{\tau-\ell}_T(T)  
\end{array}\]
with the vertical arrows being the pull-backs, along the projection of 
$BT$ to a point and the multiplication $BT^2\to BT$. Our anti-diagonal 
class, in the upper right, is the pull-back of the sum of the irreducible 
classes in $K^{\tau-\ell}_T(T)$. But we identified this in \S\ref{dirimg} 
with $p_*[1]$, as desired.  

\begin{proof}[Proof of (\ref{pw}).] 
Fix a twisting element $f$ in the quasi-torus; we prove the theorem for 
$L_fG$. We use the notation of \S\ref{topktorus} and \S\ref{topkcompact}, 
except that we write $G$ for $G(f)$, $N$ for $N(f)$, $W$ for $W^f$ for 
simplicity. 

\subsection*{Step 1.} In view of the commutative square, in which 
$\omega_*(1) = 1$,   
\[
\begin{array}{ccc}
 K_N(fT) & \xrightarrow{\;\iota_*\;} & K^{\tau\times(-\tau)}
 								_{N\times N}(fT\times fT)  \\
\vspace{-5pt}\\
\big\downarrow\scriptstyle{\omega_*} &  & 
		\quad\big\downarrow\scriptstyle{\omega_*}  \\
\vspace{-5pt} \\
 K_G(fG_1) & \xrightarrow{\;\iota_*\;}  &  K^{\tau\times(-\tau)}
 					_{G\times G}(fG_1\times fG_1) 
\end{array}
\]
it suffices to settle the upper $\iota_*$: that is, we may assume $G=N$. 

\subsection*{Step 2.} 
Let $\delta(N)$ be the left equaliser of the two projections $N^2
\rightrightarrows N/\ul T$. Its Ad-action on $fT^2$ preserves the 
diagonal cp[y of $fT$. With $\ul\ell = \dim\ul T$, we can factor 
$\iota^\tau_*$ as 
\begin{equation}\label{1712}
K^0_N(fT) \xrightarrow{B\mathrm{diag}_*} 
	K^{\tau\times(-\tau)-\ul\ell}_{\delta(N)}(fT) 
		\xrightarrow{\mathrm{diag}_*}
	K^{\tau\times(-\tau)+0}_{N^2}(fT\times fT).
\end{equation}
Moreover, we have the ``key Lemma" isomorphisms for $M=\ul T^2$ in 
$\delta(N)$ and $N^2$, 
\begin{equation}\label{1713}\begin{split}
K^{\tau\times(-\tau)-\ul\ell}_{\delta(N)}(fT) &\cong 
	K^{\tau"\times(-\tau")+0}_{\waff}
	\left(\ul{\wtlat}^\tau\times \ul{\wtlat}^{-\tau}\right), \\
K^{\tau\times(-\tau)+0}_{N^2}(fT\times fT) &\cong 
	K^{\tau"\times(-\tau")+0}_{\waff^2}
	\left(\ul{\wtlat}^\tau\times\ul{\wtlat}^{-\tau}\right).
\end{split}\end{equation}
and, as in \eqref{idcomp}, $\mathrm{diag}_*$ is the push-forward from 
upper to lower $K$-groups. We are reduced to showing that $B\mathrm
{diag}_*[1] \in K^{\tau\times(-\tau)-\ul\ell}_{\delta(N)}(fT)$ 
is the anti-diagonal class in the upper right group. 

\subsection*{Step 3.} 
Call $\delta(\nafff)$ the left equaliser of the projections $\nafff 
\times \nafff\rightrightarrows \waff$. The presentation \eqref{homog} 
of $fT$ as a homogeneous space for $N\ltimes\ul\frt$ leads to the isomorphisms 
\begin{equation}
\begin{split}\label{????}
K_N(fT) &\cong K^{\tau\times(-\tau)}_{\nafff}(\ul\frt), \\
K^{\tau\times(-\tau)}_{\delta(N)}(fT) &\cong 
	K^{\tau\times(-\tau)}_{\delta(\nafff)}(\ul\frt),
\end{split}
\end{equation}
as flagged in Remark \ref{otherpicture}. The twisting $\tau\times(-\tau)$ 
is null on the diagonal copy of $\nafff$ in $\delta(\nafff)$, but 
trivialising it in relation to the other twistings is the key step in 
finding $B\mathrm{diag}_*$. 

\subsection*{Step 4.} 
Call $\cO(\tau)$ the line bundle over $\ul T \cong \delta(\nafff)/\nafff$ 
descended from the line bundle of the extension $\tau\times(-\tau)$ of 
$\delta(\nafff)$. This $\cO(\tau)$ carries a projective action of 
$\delta(\nafff)$, by left translations, and its space of sections over 
$\ul T$ is, by definition, the representation $\mathrm{Ind}[1]$ induced 
from $\bC$ under the embedding $\naff \subset \delta(\naff)^{\tau
\times(-\tau)}$. This is the sought-after class $[\Delta^-]$ in \eqref{1713}. 

\subsection*{Step 5.} 
Finally, we show that, under the standard trivialisation of the extension 
$\tau\times(-\tau)$ over $\naff$, the direct image of $[1]$ along the 
topological induction
\[
K^{^{\tau\times(-\tau)}+0}_{\naff}(\ul\frt) 
	\xrightarrow{B\mathrm{diag}_*}
	K^{\tau\times(-\tau)-\ul\ell}_{\delta(\naff)}(\ul\frt)
\]
is represented by the Dirac family on $\ul\frt$ coupled to $\mathrm{Ind}[1]$. 
This implies its agreement with $[\Delta^-]$. The argument merely repeats 
the discussion in \S\ref{dirimg}, after observing that $[1]$ corresponds 
to the class of $\cO(\tau)$ in the chain of isomorphisms  
\[
[1]\in K_{\naff}(\ul\frt) \cong 
	K^{\tau\times(-\tau)}_{\naff}(\ul\frt) \cong
	K^{\tau\times(-\tau)}_{\delta(\nafff)}(\ul T\times \ul\frt).
\]   
\end{proof}


\part*{Appendix}
\addtocontents{toc}{\small\bfseries Appendix}
\appendix

\section{Affine roots and weights in the twisted case}
\label{afftwist}
We recall here the properties of diagram automorphisms, which lead to 
a concrete description of the twisted affine algebras in terms of simple, 
finite-dimensional ones. The connection between the two questions is 
due to Kac, to which we refer for a complete discussion \cite[\S7.9 
and \S7.10]{kac}; but we reformulate the basic facts more conveniently 
for us.

\subsection{}\label{fromkac}
When $\frg$ is simple, the order of a diagram automorphism $\vep$ is 
$r=1,2$ or 3. We assume that $\vep \ne 1$; $\frg$ must then be simply 
laced. We summarise the relevant results from \cite{kac}. 

\begin{itemize} 
\item The invariant sub-algebra $\ul \frg := \frg^\vep$ is simple, with 
Cartan sub-algebra $\ul \frt := \frt^\vep$ and Weyl group $\ul{W} := 
W^\vep$. 
\item The simple roots are the restrictions to $\ul\frt$ of those of 
$\frg$. 
\item The ratio of long to short root square-lengths in $\ul \frg$ is 
$r$, save for $\frg =\mathfrak{su}(3)$, when $\ul{\frg}= \mathfrak{su}(2)$. 
\item The $\vep$-eigenspaces are irreducible $\ul\frg$-modules. The two 
$\vep\ne 1$-eigenspaces are isomorphic when $\frg =\mathfrak{so}(8)$ 
and $r=3$. 
\end{itemize}

\subsection{The weight $\ul\theta$.} \label{newtheta}
Denote by $\ul \theta$ the highest weight of $\frg \left/\ul \frg\right.$, 
and let $a_0=2$ when $\frg =\mathfrak{su}(2\ell+1)$ and $r=2$; else, 
let $a_0=1$. Then, ${\ul\theta}\left/{a_0}\right.$ is the short dominant 
root of $\ul\frg$. (When $\frg = \mathfrak{su}(2\ell+1)$, $\ul{\frg}= 
\mathfrak{so}(2\ell+1)$ and $\frg \left/\ul \frg\right.$ is ${\mathrm{Sym}^2
\bR^ {2\ell+1}}\left/\bR\right.$, whose highest-weight is twice the short 
root.) The basic inner product on $\frg$ restricts to $a_0$ times the 
one on $\ul\frg$; so $\ul{\theta}^2 = {2a_0}\left/r\right.$. 

\begin{remark} With reference to \cite[VI]{kac}, we have $\ul \theta 
= \sum a_i\alpha _i - a_s\alpha _s$, where $s=0$, except when $\frg = 
\mathfrak{su}(2\ell+1)$, in which case $s=2\ell$: if so, our $\ul\theta$ 
differs from $\theta$ in \textit{loco citato}. 
\end{remark}

\subsection{Twisted affine Weyl group.} Denote by $\ul \fra$ the simplex 
of dominant elements $\xi \in \ul \frt$ satisfying $\ul\theta (\xi)\le 1/r$. 
The \textit{$\vep$-twisted affine Weyl group} $\waff(\frg,\vep)$ is 
generated by the reflections about the walls of $\ul \fra$. Let 
$\ul{R}^\prime\subset \ul \frt$ correspond to the root lattice 
$\ul R$ in $\ul \frt^*$ under the $\frg$-basic inner product. 

\begin{proposition}[{\cite[Props.~6.5 and 6.6]{kac}}] \label{scwaff}
$\waff(\frg,\vep)$ is the semi-direct product of the $\ul{R}'$-translation 
group by $\ul W$. Its action on $\ul \frt$ has $\ul\fra$ as fundamental 
domain. The $\waff(\frg,\vep)$-stabiliser of any point in $\ul\fra$ is 
generated by the reflections about the walls containing it. 
\end{proposition}

\begin{proof} This follows from the analogous result for the untwisted
affine algebra based on the Langlands dual to $\ul\frg$, in which
$\ul\fra$ is the fundamental alcove and $\ul{R}'$ the co-root lattice.  
\end{proof}

\subsection{Twisted conjugacy classes.} When $G$ is simply connected, 
the points of $\fra$ parametrise the conjugacy classes in $G$. The 
alcove $\ul \fra$ fulfils the same role for the $\vep$-twisted 
conjugation $g: h\mapsto g\cdot h \cdot\vep(g)^{-1}$. 

\begin{proposition} \label{twistedconj}
If G is simply connected, every $\vep$-twisted conjugacy class in has 
a representative $\exp (2\pi\xi)$, for a unique $\xi \in \ul \fra$. 
The twisted centraliser of $\exp (2\pi \xi)$ in G is connected, and its 
Weyl group is isomorphic to the stabiliser of $\xi$ in $\waff(\frg,\vep)$. 
\end{proposition}

\begin{proof} For the first part, we must show, given (\ref{scwaff})
and (\ref{gentwistedconj}), that the integer lattice of $T_\vep$ in 
$\frt_\vep \cong \ul\frt$ is $\ul R^\prime$, and that the $\vep$-twisted 
action of $\ul W$ on $T_\vep $ is the obvious one. Now, the first lattice 
is the image, in the quotient $\frt_\vep$ of $\frt$, of the integer 
(co-root) lattice $R^\vee \subset \frt$ of $T$. As $\frg$ is simply 
laced, $R^\vee$ is identified with the root lattice $R$ of $\frg$ in 
$\frt^*$ by the basic inner product, so the integer lattice of $T_\vep$ 
is also the image of $R$ in $\ul\frt^*$. But, by (\ref{fromkac}), this 
agrees with the root lattice $\ul R$ of $\ul\frg$. Concerning $\ul W$, 
since that is the Weyl group of $G^\vep $, we can find $\vep$-invariant 
representatives for its elements, and their $\vep$-action coincides with 
the usual one. 

Connectedness of twisted centralisers, for simply connected $G$, is due 
to Borel \cite{borel}. Moreover, because maximal tori are maximal abelian 
subgroups, $T^\vep$ is connected as well; and \eqref{gentwistedconj} 
identifies the Weyl groups of centralisers as desired. 
\end{proof}

\begin{remark} Connectedness of $T^\vep$ can also be seen directly, 
as follows. Clearly, the $\vep$-fixed point set  $\exp\left(\fra^\vep
\right)$ in the simplex $\exp(\fra)$, is connected. By regularity of 
$\ul \frt$, every component of $T^\vep$ contains a regular element. 
This must be conjugate to some $a \in \exp(\fra)$, hence of the form 
$w(a)$, with $w\in W$ and $a\in \exp(\fra)$. Invariance under $\vep$ 
implies $w(a)=\vep(w(a)) = \vep (w)(\vep (a))$. As $a$ and $\vep(a)$ 
are both in $\fra$ and regular, it follows that $w=\vep (w)$ and $a=
\vep(a)$, so $w(a)$ is in the $\ul W$-image of $\exp(\ul{\frt})$, hence 
in $\ul T$.
\end{remark}

\subsection{Affine roots and weights.} \label{afroot} The sub-algebra
$\ul{\mathfrak h} = \mi\bR K\oplus\ul\frt\oplus\mi\bR E$ plays the role 
of a \textit{Cartan sub-algebra} of $\hatplge$. The \textit{affine roots}, 
living in $\ul{\mathfrak h}^*$, are the $\frh$-eigenvalues of the adjoint 
action on $\hatplge$. Define the elements $\delta$ and $\kstar$ of $\ul 
{\frh}^*$ by $\delta(E) = 1/r$, $\kstar(K)=1$, $\delta(K) = \delta(\ul 
\frt) = \kstar(\ul \frt) = \kstar(E)=0$. The \textit{simple affine roots} 
are the simple roots of $\ul\frg$, plus $\delta -\ul\theta$; their $\bZ$-span 
is the \textit{affine root lattice} $R_\mathrm{aff}$. The \textit{positive 
roots} are sums of simple roots. The \textit{standard nilpotent sub-algebra} 
$\ul\frN$ is the sum of the positive root spaces, and a triangular 
decomposition $\hatplge_\bC = \ol{\ul\frN} \oplus \ul{\frh}_\bC \oplus 
\ul\frN$ is inherited from $\widehat{L}'\frg_\bC$. 

The \textit{simple co-roots} are those of $\ul\frg$, plus $(K-\beta^\vee)
\left/a_0\right.$, where $\beta^\vee$, the long dominant co-root of $\ul
\frg$, satisfies $\lambda(\beta^\vee) = {\langle\lambda|\ul\theta\rangle}/
r$.\footnote{Recall from \ref{newtheta} that ${\ul\theta}\left/{a_0}\right.$ 
is a short root, and $\ul{\theta}^2= {2a_0}\left/r\right.$.} The restriction 
$\ul{\widetilde T}$ of the basic central extension \eqref{basictwext} to 
$\ul T$ is the quotient of $\mi\bR K\oplus\ul\frt$ by the affine co-root 
lattice $R_\mathrm{aff}^\vee$.
 
The \textit{weight lattice} $\widetilde\wtlat$ of $\widetilde{L}_\vep G$, 
in $\ul{\frh}^*$, is the integral dual of $R_\mathrm{aff}^\vee$, and 
comprises the characters of $\ul{\widetilde T}$. Calling $\ul\wtlat$ 
the (simply connected) weight lattice of $\ul\frg$, we have

\begin{equation} \label{weights}
\widetilde{\wtlat} = 
\begin{cases}
\bZ\kstar \oplus\ul\wtlat, &\text{if } \frg 
\ne \mathfrak{su}(2\ell+1),\\ 
2\bZ\kstar \oplus\ul{\wtlat}^+ \cup (2\bZ+1)\kstar \oplus
\ul{\wtlat}^- &\text{if } \frg = \mathfrak{su}(2\ell+1)
\end{cases}
\end{equation}
the superscript indicating the value of the character on the central 
element of $\mathrm{Spin}(2\ell +1)$. The \textit{affine weight lattice} 
$\widehat\wtlat$ includes, in addition, the multiples of $\delta$, giving 
the energy eigenvalue.

The \textit{dominant weights} pair non-negatively with the simple
co-roots; this means that $(k,\lambda,x)$ is dominant iff  $\lambda$ 
is $\ul \frg$-dominant and $\lambda\cdot\ul\theta \le k/r$. 
The affine Weyl group $\waff(\frg,\vep )$ preserves the constant 
level hyperplanes, and its lattice part $\ul R$ (\ref{scwaff}) acts 
by $k$-fold translation at level $k$. Every positive-level weight has 
a unique dominant affine Weyl transform. {Regular weights} are those 
not fixed by any reflection in $\waff(\frg,\vep)$. The important 
identity $\langle \ul{\rho} |\ul{\theta} \rangle + {\ul \theta}^2/2 = 
{h^\vee}/r$ \cite[VI]{kac} implies that an \textit{integral} weight 
$(k,\lambda,x)$ is dominant iff the shifted weight $(k+h^\vee, 
\lambda +\ul \rho,x)$ is dominant regular.


\vspace{1.5cm}
\small{
\noindent\textsc{D.S.~Freed:} Department of Mathematics, UT Austin,
Austin, TX 78712, USA. \texttt{dafr@math.utexas.edu}\\
\textsc{M.J.~Hopkins:} Department of Mathematics, MIT, Cambridge, MA 
02139, USA. \texttt{mjh@math.mit.edu}\\
\textsc{C.~Teleman:} DPMMS, CMS, Wilberforce Road, Cambridge CB2 1TP, UK.
\texttt{teleman@dpmms.cam.ac.uk}}

\end{document}